%% file: main.tex
\documentclass{article}
\usepackage[english]{babel}

\usepackage[letterpaper,top=2cm,bottom=2cm,left=2cm,right=2cm,marginparwidth=1.75cm]{geometry}

\usepackage{float}
\usepackage{calrsfs}
\usepackage{amssymb}
\usepackage{amsmath}
\usepackage{amsthm}
\usepackage{subcaption}
\usepackage{multirow}
\usepackage{graphicx}
\usepackage{amsfonts}
\usepackage{mathtools}
\usepackage{mathrsfs}
\usepackage{placeins}
\usepackage{comment}
\usepackage{bbm}
\usepackage{dsfont}
\usepackage[export]{adjustbox}
\usepackage{overpic}

\usepackage{algorithm}
\usepackage{algpseudocode}

\usepackage{booktabs}

\usepackage{empheq}
\usepackage[mathscr]{euscript}
\usepackage[dvipsnames,x11names]{xcolor}
\usepackage[colorlinks=true, allcolors=OliveGreen]{hyperref}
\usepackage{amsmath}   
\usepackage{tikz}
\usepackage{pgfplots}
\usetikzlibrary{external}
\usepackage{authblk}
\usepackage{array}
\newcolumntype{C}{>{\centering\arraybackslash}p{2.5cm}}
\usepackage{nicematrix}
\usepackage{caption}
\usepackage{subcaption}
\usepackage{tikz}

\usetikzlibrary{patterns,tikzmark}
\usetikzlibrary{matrix,decorations.pathreplacing,calc}
\usepackage{hf-tikz}

\usepackage{orcidlink}
\usepackage{csvsimple}
\usepackage{siunitx}
\usepackage{tikz} 

\usetikzlibrary{arrows.meta, positioning, shapes.geometric}
\usepackage{pgfplots}
\pgfplotsset{compat=newest} 
\usepgfplotslibrary{units} 
\usepackage{tikz}
\usetikzlibrary{matrix,decorations.pathreplacing}
\pgfkeys{/pgfplots/legend entry/.code=\addlegendentry{#1}}

\usepackage{pdfpages}

\newcommand{\averagel}{\{\!\!\{}
\newcommand{\averager}{\}\!\!\}}

\newcommand{\jumpl}{[\![}
\newcommand{\jumpr}{]\!]}

\newcommand{\Napl}{\mathrm{Na}^+}

\newcommand{\Kpl}{\mathrm{K}^+}

\usepackage[export]{adjustbox}
\usepackage{lipsum}                     
\usepackage{xargs}                      
\usepackage[colorinlistoftodos,prependcaption]{todonotes}
\newcommandx{\unsure}[2][1=]{\todo[linecolor=red,backgroundcolor=red!25,bordercolor=red,#1]{#2}}
\newcommandx{\change}[2][1=]{\todo[linecolor=blue,backgroundcolor=blue!25,bordercolor=blue,#1]{#2}}
\newcommandx{\info}[2][1=]{\todo[linecolor=OliveGreen,backgroundcolor=OliveGreen!25,bordercolor=OliveGreen,#1]{#2}}
\newcommandx{\improvement}[2][1=]{\todo[linecolor=Plum,backgroundcolor=Plum!25,bordercolor=Plum,#1]{#2}}
\newcommandx{\thiswillnotshow}[2][1=]{\todo[disable,#1]{#2}}

\DeclareMathAlphabet{\mathcalligra}{T1}{calligra}{m}{n}

\graphicspath{{Images/}}

\tikzset{  font={\fontsize{15pt}{12}\selectfont}}

\title{Mathematical and numerical modeling of coupled oxygen dynamics
and neuronal electrophysiology\footnote{\textbf{Funding}: This work is partially funded by the European Union (ERC SyG, NEMESIS, project number 101115663). Views and opinions expressed are, however, those of the authors only and do not necessarily reflect those of the European Union or the European Research Council Executive Agency. Neither the European Union nor the granting authority can be held responsible for them. FD has been partially funded by the National Recovery and Resilience Plan (NRRP) Mission 4 - Component 2 - Investment 1.2, funded by European Union. CBLS has been funded by the National Recovery and Resilience Plan (NRRP), Mission 4, Component 1 – Investment 3.4 and Investment 4.1, funded by the European Union. 
The present research is part of the activities of the Dipartimento di Eccellenza 2023-2027 grant, funded by MUR. PFA, SP and CBLS are members of INdAM-GNCS. }}

\author[1]{Francesco Daniele\orcidlink{0009-0009-1473-3524}}

\author[1]{Caterina B. Leimer Saglio\orcidlink{0009-0007-7887-919X}}

\author[1]{Stefano Pagani\orcidlink{0000-0002-6662-3433}}

\author[1]{Paola F. Antonietti\orcidlink{0000-0002-2138-3878}}

\affil[1]{MOX-Dipartimento di Matematica, Politecnico di Milano, Piazza Leonardo da Vinci 32, Milan, 20133, Italy}

\begin{document}
\maketitle

\begin{abstract}
Modeling and simulating how oxygen supply shapes neuronal excitability is crucial for advancing the understanding of brain function in pathological scenarios, such as ischemia. This condition is caused by a reduced blood supply, leading to the deprivation of oxygen and other metabolites; this energy deficit impairs ionic pumps and causes cellular swelling. In this work, this phenomenon is modeled through a volumetric variation law that links cell swelling to local oxygen concentration and the percentage of blood flow reduction. The swelling law links volume changes to local oxygen and the degree of blood-flow depression, providing a simple mechanistic pathway from hypoxia to tortuosity-driven transport impairment. The interplay between oxygen supply and excitability in brain tissue is described by coupling the monodomain model with specific neuronal ionic and metabolic models that characterize ion and metabolite concentration dynamics. 
The numerical approximation of this coupled multiscale problem is particularly challenging, owing to the presence of sharp and fast-propagating wavefronts and complex geometrical domains.
To address these challenges, suitable space- and time-adaptive schemes are employed to capture the action potential dynamics accurately.
This multiscale model is discretized in space with the high-order $p$-adaptive discontinuous Galerkin method on polygonal and polyhedral grids (PolyDG) and integrated in time with a Crank-Nicolson scheme. We numerically investigate different pathological scenarios on a two-dimensional idealized square domain and on a realistic geometry, both discretized with a polygonal grid, analyzing how subclinical and severe ischemia can affect brain electrophysiology and metabolic concentrations.
\end{abstract}

\section{Introduction}\label{sec:1}
\input{introduction}

\section{The mathematical model}
\label{sec:mathematicalmodel}

\input{mathematicalmodel}

\section{PolyDG formulation}
\label{sec:PolyDG}
\input{PolyDG}

\section{Numerical results}
\label{sec:NumericalResults}
\input{Numerical_result0D}

\input{Num2D}
\input{NumericalresultBrain}

\section{Conclusions}
\label{sec:10}
\noindent This work introduces an electro–metabolic framework that couples oxygen availability and ATP production to ionic homeostasis and membrane excitability under ischemic stress. The coupling is modeled by taking into account ATP-dependent Na$^+$/K$^+$ pumping and glial $\text{K}^+$ uptake. At the same time, metabolism provides oxygen/ATP dynamics, as well as a swelling law that ties volume to local hypoxia and blood-flow depression. The system is coupled with the monodomain model and discretized with a high-order, $p$-adaptive PolyDG discretization and semi-implicit time integration, allowing the simultaneous resolution of fast spikes and slow ionic–metabolic changes. The numerical results demonstrate how energy limitations reshape excitability. As ischemic severity increases, the model transitions from near-physiological behavior to high-frequency, burst-like discharges, accompanied by extracellular potassium accumulation and intracellular sodium rise. This establishes a quantifiable transition from oxygen shortage to high-frequency activity, marking the onset of epileptiform activity. A complementary sensitivity analysis reveals that extracellular potassium transport acts as an amplifier: decreasing the diffusion coefficient sustains the positive feedback between firing and further depolarization, thereby pushing the system into pathological regimes. Beyond homogeneous settings, a localized ischemic subregion of grey matter allows spikes to propagate into surrounding tissue, showing how local energy failure can trigger organ-scale activity when oxygen transport is impaired. Future developments will focus on biological and computational contexts. On the physiological side, we will calibrate parameters to incorporate explicit perfusion and vascular reactivity \cite{koppl20203d}, and account for dynamic tortuosity. From the numerical point of view, we will pursue scalable implementations for 3D simulations to connect the model to realistic measurements.

\bibliographystyle{hieeetr}
\bibliography{sample.bib}

\end{document}

%% file: introduction.tex
The human brain is a complex and highly organized system that accounts for approximately 20$\%$ of the energy consumption of the body, despite only representing 2$\%$ of its weight \cite{kety1957metabolism_energy, clarke1999_energy}. To reduce energy consumption and conduction delays, local and long-term connections of the brain are organized into anatomically distinct regions: white matter and grey matter. White matter is primarily composed of axons and oligodendrocytes, whereas the grey matter consists of neurons and unmyelinated fibers. Glial cells are present in both grey and white matter; among them, astrocytes are found with a higher density in grey matter, where they play a key role in synaptic regulation and metabolic coupling \cite{beard2022astrocytes}. Neuronal electrical activity strongly depends on continuous energy supply, and at the same time, metabolic dynamics are regulated by neural activation. Modeling this two-way coupling enables an accurate representation of brain function, which is regulated by the balance between energy demand and energy production. This is a crucial point to understand both physiological and pathological conditions in which energy availability and neuronal activity are simultaneously impaired \cite{goriely2015mechanics}.

In this work, we focus on the modeling of ischemic events in realistic brain geometries, encompassing both grey and white matter tissues, where ischemic processes predominantly originate in the former, and analyzing both biochemical reactions and the electrical response of neural cells. Within this framework, we consider neurons and astrocytes, which differ in both structure and function \cite{BONVENTO20211546_neuron_astrocyte, allaman2011astrocyte_neuron}: neurons are excitable cells specialized in transmitting electrical signals, while astrocytes are non-excitable glial cells involved in metabolic processes and modulation of synaptic activity. Once a neuron is activated, sodium ions flow into the cell and potassium ions flow out; this imbalance is restored under physiological conditions by Na$^+$/K$^+$-ATPase pumps, sustained by ATP consumption \cite{CAPORANGEL2020110093}. Astrocytes play a key role in $\text{K}^+$ clearance, which is fundamental to prevent abnormal firing and epileptic patterns \cite{Potassio_buffer_Contreras2021}.
Glial cells contribute to restoring ionic homeostasis, and pathological increases of extracellular potassium can intensify excitability and trigger spontaneous activity. Cerebral metabolism refers instead to the ensemble of biochemical processes through which the brain utilizes glucose and oxygen to sustain neuronal activity, maintain ionic gradients, and support cellular homeostasis. Glucose is metabolized through glycolysis and oxidative phosphorylation to produce ATP, required to fuel membrane pumps and other bioenergetic processes \cite{attwell2001energy}. Astrocytes regulate glucose uptake, neurotransmitter recycling, and ion buffering, thereby coupling metabolism with electrophysiology.  

From a modeling perspective, coupling these processes is challenging, since electrophysiology evolves on millisecond timescales while metabolism typically evolves over seconds to minutes. In this work, we develop a multiscale model, extending the Barreto–Cressman ionic model \cite{hodgkin1952propagation,cressman2009_BC1,barreto2011ion} with astrocytic contributions, and combining it with a reduced metabolic model inspired by \cite{CALVETTI2018, idumah2023spatial_calvetti2}. The aim is to capture the interplay between neuronal excitability, energy supply, and astrocytic regulation, with a particular focus on pathological conditions such as cerebral ischemia and cytotoxic edema.

Ischemia results from reduced blood flow and oxygen availability, impairing Na$^+$/K$^+$ pumps and leading to sodium accumulation inside cells \cite{muller}. Due to osmotic effects, this drives water influx and swelling, especially in astrocytes, which further reduces extracellular space and alters ionic diffusion \cite{rungta2015cellular_edema}. These mechanisms contribute to seizures, spreading depression, and eventually neuronal death \cite{cell_volume_seizures_ullah}. Ischemic stroke can promote the onset of pathological high-frequency electrical activity and epileptic seizures by altering the brain normal electrical homeostasis \cite{lee2000brain, chen2022poststroke}. The sudden deprivation of oxygen and nutrients destabilizes neuronal membranes, leading to uncontrolled ion movements and increased excitability of neurons. The massive release of glutamate aggravates this effect, as the main excitatory neurotransmitter overstimulates its receptors and promotes the influx of sodium and calcium. These changes disrupt normal communication between neurons and support cells, causing abnormal connections and persistent hyperexcitability.

At the tissue and organ scale, ionic models are coupled with the monodomain or bidomain equations \cite{mardal2022mathematical,jaeger2022deriving,leimer2024highorder} to describe the propagation of the transmembrane electric potential.  
The numerical simulation of these processes is computationally demanding: the sharp and fast wavefronts of the transmembrane potential require fine spatial resolutions and very small timesteps, which makes standard finite element approaches extremely computationally expensive \cite{HamBathe2012_FEM_wave,kager2007seizure,somjen2008computer,stefanescu2012computational,mardal2022mathematical,schreiner2022simulating,saetra2023neural}.  
To overcome these limitations, we employ a discontinuous Galerkin scheme on polytopal grids (polyDG \cite{Antonietti2013A1417,Bassi201245,cangiani_hp-version_2017,cangiani2022hp,antonietti2025lymph}), which provides flexibility for complex geometries and heterogeneous tissues, efficiently resolves steep gradients, and naturally supports $h$-, $p$-, and $hp$-adaptivity strategies due to the local nature of the discretization spaces. Temporal adaptivity is introduced to handle the multiscale dynamics in the neuronal model, while spatial $p$-adaptivity \cite{leimer2025p} selectively increases the polynomial degree in regions with sharp wavefronts for spatial multiscale simulations. These techniques improve accuracy while significantly reducing the number of degrees of freedom of the system and the overall computational costs. 

This work presents a mathematical model of the coupling between metabolism and electrophysiology in both physiological and ischemic conditions, together with a computational study showing increased sensitivity of ischemic regions to external stimuli, potentially leading to spontaneous pathological spiking activity.
Results also provide new insights into the interaction between energy supply and neuronal dynamics, as well as the large-scale effects of localized ischemic regions.  

The remainder of the manuscript is organized as follows:
in Section~\ref{sec:mathematicalmodel}, we introduce the mathematical model of brain electrophysiology and cerebral metabolism under investigation. Section~\ref{sec:polyDG} presents the semi-discrete formulation within the PolyDG framework and the fully discrete method based on employing the Crank–Nicolson time marching scheme. In Section~\ref{sec:NumericalResults}, we report several sensitivity analyses on the proposed neuronal model, as well as numerical experiments in simplified and realistic two-dimensional settings to investigate the influence of ischemic regions on the onset and propagation of pathological neuronal activity.

%% file: mathematicalmodel.tex
In this section, we present the coupled mathematical model describing cerebral hemodynamics, metabolism, and neuronal electrophysiology. The metabolic and hemodynamic component of the model is based on the formulation introduced in \cite{CALVETTI2018}, while neuronal electrophysiology is described following \cite{cressman2009_BC1}. These models are then spatially and temporally coupled to analyze the influence of oxygen supply on epileptic events. 

\begin{figure}[h]
    \centering
    \includegraphics[width=\textwidth]{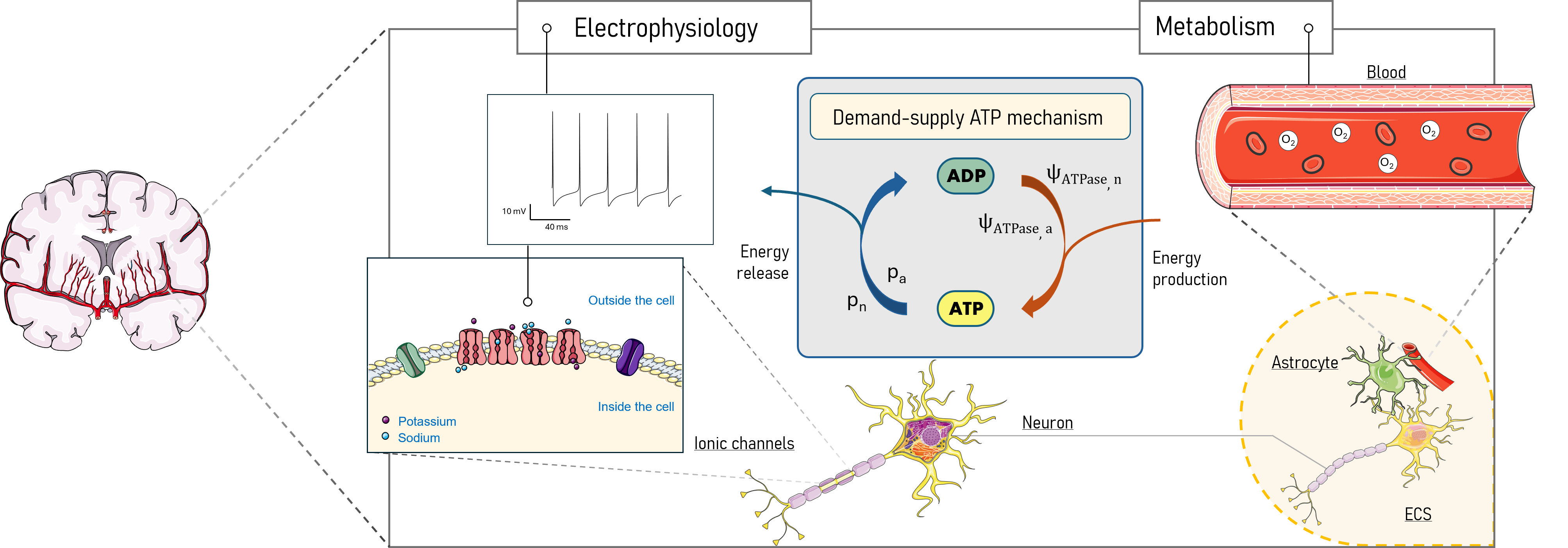}
    \caption{Schematic representation of electro-metabolic model. The equilibrium between metabolic and electrophysiological dynamics is preserved by the demand/supply mechanism of ATP production. $p_n$ and $p_a$ denote phosphorylation states of neuron and astrocyte, respectively, while outputs of ATP demand are the ATP dephosphorylation fluxes, which are additional inputs of the metabolic processes. Figure adapted from Servier Medical Art (https://smart.servier.com), licensed under CC BY 4.0 (https://creativecommons.org/licenses/by/4.0/).}
    \label{fig:schematic_model}
\end{figure}

The complex cerebral metabolism can be described using a multi-compartment well-mixed model with four separate compartments \cite{CALVETTI2018}: blood, extracellular space (ECS), neurons, and astrocytes. Cellular communication and biochemical processes occur through the ECS and are driven by their energetic needs, which are satisfied by ATP production through metabolic processes. The coupling between the two models is regulated by the availability of ATP, which ensures the proper functioning of the ion pumps. Inputs of the metabolic system are the blood flow $q(t)$ and arterial concentration $q_a(t)$. In contrast, the one for the electrophysiology model is a time-constant function $\xi(t)$, representing the external stimulus. The overall system, highlighting energetic connections, is illustrated in Figure \ref{fig:schematic_model}.

The generation of action potentials (APs) relies not only on ionic fluxes but also on the continuous supply of metabolic energy. Neurons maintain steep electrochemical gradients of sodium and potassium through the action of $\Napl$, $\Kpl$, ATPases, which consume large amounts of ATP. Astrocytes support this process by buffering extracellular potassium and participating in glutamate recycling, thereby contributing to neuronal excitability. For this reason, the description of neural activity cannot be limited to pure electrophysiology, but must include the tight link with metabolism; the two processes evolve simultaneously on different time scales within a coupled modeling framework \cite{CALVETTI2018}.

The electrophysiological component of the model is based on a modified Barreto–Cressman conductance model \cite{barreto2011ion, cressman2009_BC1}. It tracks the membrane potential, intra- and extracellular ion concentrations, and gating variables, while explicitly including astrocytic potassium uptake. Ionic currents depend not only on channel dynamics but also on the energetic state of the system: the $\Napl$/$\Kpl$ pumps and glial uptake fluxes are scaled by the availability of ATP in neurons and astrocytes. The coupling mechanism is realized through ATP dephosphorylation fluxes. The electrophysiology model prescribes the energetic demand required to sustain ionic currents, while the metabolic model supplies ATP according to oxygen availability. This bidirectional interaction reflects the fundamental balance between electrical activity and energy metabolism: intense firing requires higher ATP consumption, which, in turn, accelerates oxygen depletion in cellular compartments.
\subsection{Metabolic multi-compartment model}
\label{sec:model0D}
Metabolites are carried by the blood, first entering the ECS and subsequently the neurons and astrocytes. Biochemical reactions occur only in these two cellular compartments \cite{CALVETTI2018}. To mathematically describe all the metabolic biochemical processes, several physical laws are incorporated into the model. Among the various processes, we focus on oxygen dynamics, assuming a simplifying hypothesis that considers glucose and lactate concentrations to be constant, thereby neglecting their temporal evolution.
The oxygen concentration in blood $[O_2]_b$ satisfies the following differential equation
\begin{align}
    \eta_b \frac{d[O_2]_b}{dt} &= \frac{q}{F} \left( {q_a} - [O_2]_b \right) - J_{O_2},
\end{align}
where $q=q(t)$ is the time dependent blood flow, $q_a$ represents the arterial concentration, $F\in(0,1)$ is the mixing ratio of arterial and venous blood and $J_{O_2}>0$ is the net flux between blood and ECS.
The equation describes the balance between oxygen delivery and consumption within the system. Oxygen can be present either dissolved freely in plasma or bound to hemoglobin. The total oxygen concentration is defined according to Hill's equation
\begin{align}
    [O_2]_b &= [O_2]_{b,\text{free}} + 4 \text{Hct} [\text{Hb}] 
    \frac{[O_2]_{b,\text{free}}^n}{K_H^n + [O_2]_{b,\text{free}}^n} = H([O_2]_{b,\text{free}}),
    \label{eq:O2free}
\end{align}
\noindent where Hct is the hematocrit, $[\text{Hb}]$ is the concentration of hemoglobin in plasma, $K_H$ is the affinity constant, and $n=5/2$. The transport flux obeys a modified Fick’s law for the free oxygen concentration
\begin{align}
    J_{O_2}
    = \lambda_{b,O_2} \left( H^{-1} ([O_2]_b) - [O_2]_{\text{ECS}} \right)^{k}, \qquad k = 0.1.
\end{align}
\noindent In the extracellular space no chemical reactions take place and consequently the concentrations change only through exchanges of the substances with adjacent compartments. The oxygen dynamics here follows
\begin{align}
      \eta_{\text{ECS}} \frac{d[O_2]_{\text{ECS}}}{dt} &= J_{O_2} - j_{O_2}^{n} - j_{O_2}^{a},
\end{align}

\noindent where $j_{O_2}^{n}$ and $j_{O_2}^{a}$ are the oxygen fluxes from ECS to neuron and to astrocyte, respectively. Also in these cell compartments the oxygen transport follows the Fick's law with constant permeabilities
\begin{align}
    j_{O_2}^{y} = \lambda_{y,O_2} \left( [O_2]_{\text{ECS}} - [O_2]_{y} \right), \quad &y \in \{n, a\}.
\end{align}
In the neuronal and astrocytic compartments, the dynamics of specific metabolites are modeled; however, oxygen is exchanged between the two compartments through transmembrane transport.
In that metabolic model, only oxygen and energy-related metabolites follow a specific time-dependent dynamics, defined by ODEs, whereas for the other ones present in the model \cite{CALVETTI2018}, we assume them to be constant at their resting values. We set this simplifying assumption that reduces the complexity of the model while still capturing essential aspects of energy metabolism.
The stoichiometric matrix $S$ is introduced, with entries:
\begin{align*}
    S_{X,\text{R}}= \# \text{\ number of molecules of $X$ produced and/or depleted by reaction R,}
\end{align*}
\noindent and the differential equation of mass balance for a generic metabolite $X$ in the cellular compartment $y$ is 
\begin{align}
    \eta_y \frac{d[X]_y}{dt} = j_{X}^y + \sum_R s_{X,\text{R}} \psi_\text{R}^y,
    \quad &y \in \{n, a\},
\end{align}
\noindent where $X\in[ O_2,\ \text{ATP}, \ \text{ADP},\ \text{NADH},\ \text{NAD}^+]$ and $\psi^y_\text{R}$ are the time-dependent expressions of reaction rates. The following notation is used for the phosphorylation and redox states of the two cells 
\begin{align}
    p_n = \frac{[\text{ATP}]_n}{[\text{ADP}]_n},  
    \quad p_a = \frac{[\text{ATP}]_a}{[\text{ADP}]_a}, \quad 
    r_n = \frac{[\text{NADH}]_n}{[\text{NAD}^+]_n},  
    \quad r_a = \frac{[\text{NADH}]_a}{[\text{NAD}^+]_a}.
\end{align}
The dephosphorylation fluxes play a crucial role in the coupling mechanism: they are the energetic link between the metabolic and the electrophysiology model. The metabolism model can be summarized collecting all the concentrations in the following vector:
\begin{equation}
    \mathbf{c}(t) =
    \begin{bmatrix}
    [O_2]_b(t), & [O_2]_\text{ECS}(t), & \mathbf{c}_n(t), & \mathbf{c}_a(t)
    \end{bmatrix}^{\!\top}
    \in \mathbb{R}^{12},
    \label{vect_metabolism}
\end{equation}
where in particular $\mathbf{c}_y(t)$ represents the concentration of specific metabolites in the $y$-compartment
\begin{equation}
    \mathbf{c}_y(t) =
    \begin{bmatrix}
    [O_2]_y(t),\  [\text{ATP}]_y, \ [\text{ADP}]_y, \ [\text{NADH}]_y,\  [\text{NAD}^+]_y
    \end{bmatrix}^{\!\top}
    \in \mathbb{R}^{5}, \quad y \in \{n, a\}.
    \label{vect_an}
\end{equation}
The diagonal matrix $\mathrm{M}$ is introduced as
\begin{equation}
    \text{M} = \text{diag}(\eta_b, \eta_{\text{ECS}}, \eta_n \mathds{1}_5, \eta_a \mathds{1}_5),
\end{equation}
where $\mathds{1}_5$ is the identity matrix of dimension $\mathrm{M} \in \mathbb{R}^{5\times 5}$. The dynamics of the metabolism model can be written in the compact form:
\begin{equation}
    \text{M} \frac{d\mathbf{c}}{dt} + \boldsymbol{l}(\mathrm{u},\mathbf{c}, q_a, \psi_{\text{ATPase},n}, \psi_{\text{ATPase},a}, q) = \boldsymbol{0},
    \label{eq:metabolism}
\end{equation}
where arterial concentration $q_a$ and blood flow $q=q(t)$ are the system inputs, while through dephosphorylation fluxes energy is exchanged in order to allow proper communication between metabolism and electrophysiology. All parameters related to the metabolic model are listed in Table \ref{table:metabolism_params}.
\begin{table}[h]
    \centering
    \resizebox{0.85\textwidth}{!}{%
    \begin{tabular}{|cc|ccc|ccc|}
    \hline
    \multicolumn{2}{|c|}{\text{Volume fractions}} &
    \multicolumn{3}{c|}{\text{Blood flow parameters}} &
    \multicolumn{3}{c|}{\text{Oxygen parameters}}  \\
    \hline
    \text{Parameter} & \text{Value} &
    \text{Parameter} & \text{Value} & \text{Units} &
    \text{Parameter} & \text{Value} & \text{Units}  \\
    \hline
    $\eta_n$ & 0.4 & Hct & 0.45 &  & $q_a$ & 9.14 & [mM]  \\
    $\eta_a$ & 0.3 & Hb & 5.18 &  & $\lambda_{b,O_2}$ & 0.04 & $[\text{mM}]^{1-k}[s]^{-1}$ \\
    $\eta_{\text{ECS}}$ & 0.3 & $\text{K}_H$ & $36.4 \cdot 10^{-3}$ & [mM] & $\lambda_{n,O_2}$ & 0.94 & [s]$^{-1}$ \\
    $\eta_b$ & 0.04 & $q$ & 0.40 & [mL][min]$^{-1}$ & $\lambda_{a,O_2}$ & 0.68 & [s]$^{-1}$ \\
     &       & $F$ & $2/3$ &  & & &  \\
    \hline
    \end{tabular}
    }
    \caption{Physiological parameters of metabolic model: compartments volume fractions, blood flow parameters, arterial oxygen concentration, and Fick’s law parameters. Data taken from \cite{CALVETTI2018}.}
    \label{table:metabolism_params}
\end{table}

\subsection{Conductance-based ionic model and coupling conditions}
The electrophysiological model presented is a modified version of the classical Barreto–Cressman model \cite{cressman2009_BC1}, including astrocytic potassium clearance and potassium diffusion, as clarified in \cite{CALVETTI2018}, and also considering calcium dynamics. In order to enable the coupling with the metabolic model, the electrophysiology model reads:
\begin{subequations}
\label{eq:BC}
\begin{empheq}[left=\empheqlbrace]{alignat=4}
       \frac{d\mathrm{[Ca^{2+}]_i}}{dt} &=  -\frac{\mathrm{[Ca^{2+}]_i}}{80} - \mathrm{G_{Ca}}\frac{0.002(\mathrm{u}-E\mathrm{_{Ca}})}{1 + \exp\left(-\frac{25 + \mathrm{u}}{2.5} \right)}  \label{eq:BC1}\\
         \eta_\text{ECS}\tau \frac{d\mathrm{[K^+]_o}}{dt} &=  \gamma \eta_n I_\text{K} - 2 \psi_\text{ATPase,n}  - 2 \psi_\text{ATPase,a}  - \varepsilon_\text{ECS}(\mathrm{[K^+]_o} - \text{K}_\text{bath}), \label{eq:BC2}
       \\
        \eta_n\tau \frac{d\mathrm{[Na^+]_i}}{dt} &= - \gamma \eta_n I_\text{Na} - 3 \psi_\text{ATPase,n}, \label{eq:BC3}
       \\
        \frac{d\mathrm{w}}{dt} &= \varphi \left( \alpha_\mathrm{w}(\text{u})(1 - \mathrm{w}) - \beta_\mathrm{w}(\text{u})\mathrm{w} \right), \quad \mathrm{w} \in \{h, n, m\}, \label{eq:BC4}
\end{empheq}
\end{subequations}
where $\alpha \ \text{and}\ \beta$ are the voltage-dependent expressions of the gating variables, and the ionic currents, which are a modified version of the ones proposed in Hodgkin-Huxley model \cite{hodgkin1990_HH}, are defined as in Equations \eqref{eq:ionic_BC} and \eqref{eq:ionic_BC_2} .
\begin{align}
\label{eq:ionic_BC}
    I_{\text{Na}} &= g_{\text{Na}} m^3 h (\mathrm{u} - E_{\text{Na}}) + g_{\text{Na},\text{leak}}(t) (\mathrm{u} - E_{\text{Na}}),\\
    I_{\text{K}} &= g_{\text{K}} n^4 (\mathrm{u} - E_{\text{K}}) + g_{\text{K},\text{leak}}(t) (\mathrm{u} - E_{\text{K}}), \quad
    I_{\text{Cl}} = g_{\text{Cl}} (\mathrm{u} - E_{\text{Cl}}).
\label{eq:ionic_BC_2}
\end{align}
In Equation \eqref{eq:BC} the coefficient $\tau$ is a conversion factor from seconds to milliseconds and $\varepsilon_{\text{ECS}}=\eta_{\text{ECS}}\varepsilon$, with $\varepsilon$ defined as the potassium concentration clearance rate. The reversal potentials $E_{\text{Na}},\ E_{\text{K}},\ E_{\text{Cl}}$ are computed using the Nernst equation:
\begin{equation}
E_{\text{Ca}} = 120\,\text{mV}, \quad
E_{\text{Na}} = 26.64 \log\left( \frac{[\text{Na}^+]_\text{o}}{[\text{Na}^+]_i} \right), \quad
E_{\text{K}} = 26.64 \log\left( \frac{[\text{K}^+]_\text{o}}{[\text{K}^+]_i} \right), \quad
E_{\text{Cl}} = 26.64 \log\left( \frac{[\text{Cl}^-]_i}{[\text{Cl}^-]_\text{o}} \right).
\end{equation}
The leak conductances of sodium and potassium are time-dependent, allowing for the description of the metabolic response to different levels of synaptic activity, and vice versa. Rather than describe the full biochemical cycle of glutamate involved in synaptic transmission, its effect is described by a transient increase in sodium and potassium leak conductances related to the activation function:
\begin{align}
\left( g_{\text{Na},\text{leak}}(t), g_{\text{K},\text{leak}}(t) \right) = 
(1 + \xi(t)) \left( g_{\text{Na},\text{leak}}^0, g_{\text{K},\text{leak}}^0 \right),
\end{align}
\noindent where $g^0_\text{Na,leak}$ and $g^0_\text{K,leak}$ are the constant resting values of the leaks, and $\xi(t)$ is an activation function that models the effect of glutamate. 
In order to describe the connection between the electrophysiology and energetic dynamics, the ion mass currents of the Barreto-Cressman model are modified and defined as
\begin{equation}
    \begin{aligned}
        I_\mathrm{pump} &= \frac{p_n}{\mu_\text{pump} + p_n} 
        \left( \frac{\rho}{1 + \exp\left(\frac{25 - [\text{Na}^+]_i}{3}\right)} \right) 
        \left( \frac{1}{1 + \exp\left(5.5 - [\text{K}^+]_o\right)} \right), \\
        I_\mathrm{glia} &= \frac{p_a}{\mu_\text{glia} + p_a} 
        \left( \frac{G_\text{glia}}{1 + \exp\left(\frac{18 - [\text{K}^+]_o}{2.5}\right)} \right), \quad
        I_\mathrm{diff} = \varepsilon_{\text{ECS}} \left( [\text{K}^+]_o - \text{K}_{\text{bath}} \right),
    \end{aligned}
\end{equation}
where $\mu_{\text{pump}}=0.1$ and $\mu_{\text{glia}}=0.1$ are the affinity constants. These expressions describe the feedback mechanism from the metabolism to the electrophysiology. The sodium current associated with the neurotransmitter activity is
\begin{equation}
    \begin{aligned}
    & I_\text{Na,act} = \begin{cases}
        I_\text{Na,leak} - I^0_\text{Na,leak} \quad & \xi > 0 \\
        0 \quad & \xi = 0 ,
    \end{cases} 
    \end{aligned}
\end{equation}
where $I^0_{\text{Na}^+,\text{leak}}$ denotes the average sodium leak current during baseline activity and the leak current of sodium reads as 
\begin{equation}
    I_\text{Na,leak}(t) = g_\text{Na,leak}(t)(\mathrm{u}-E_{\text{Na}}).
\end{equation}
The following mathematical expressions describe the ATP dephosphorylation in neuron and astrocyte, considering the energetic cost of glutamate–glutamine recycling
\begin{equation}
    \begin{aligned}
    & \psi_\text{ATPase,n} = H_1 + s(\eta_n I_\text{pump,Na} + 0.33\frac{\gamma}{\sigma} I_\text{Na,act}) \\
    & \psi_\text{ATPase,a} = H_2 + s(\frac{\eta_\text{ECS}}{2} I_\text{glia} + 2.33\frac{\gamma}{\sigma} I_\text{Na,act}), \\
    \end{aligned}
\end{equation}
where $H_1$ and $H_2$ indicate the energy required to perform routine metabolic task and $s>0$ is a suitable calibration parameter. The parameter $\gamma$ is the conversion factor from electric currents to mass fluxes and $\sigma$ accounts for the energetic cost of cycling glutamate-glutamine during synaptic activity \cite{attwell2001energy, CALVETTI2018}.

In the context of ischemic conditions, the volumetric variation of cells is closely related to the pressure of ionic concentrations \cite{dijkstra2016biophysical}, particularly oxygen. The relation introduced here simplifies this relationship and establishes a direct connection between cell volume and oxygen concentration. The equation used to represent this variation is:
\begin{equation}
    \label{eq:volume_variation}
    V(t) = V_0 \left( 1 + \delta \left( 1 - \frac{[O_2](t)}{[O_2]_\text{0}} \right) \right),
\end{equation}
where $V_0$ is the initial cell volume, $\delta$ represents the percentage of blood reduction in the tissue, reflecting the extent of ischemia, $[O_2](t)$ is the oxygen concentration, and $[O_2]_0$ is the baseline oxygen concentration. The equation \eqref{eq:volume_variation} suggests that cell volume increases as oxygen concentration decreases, due to ischemia, leading to edema \cite{Hellas2021NeuronalSwelling}.
Acute ischemia and the resulting hypoxia disrupt neuronal membranes and metabolism, impairing the Na$^+$/K$^+$-ATPase and causing intracellular Na$^+$ and water accumulation. After a stroke, astrocytes become reactive and may contribute to the formation of glial scars. In particular, astrocytes exhibit a reduced ability to regulate extracellular potassium and neurotransmitters, such as glutamate, resulting in abnormal ionic and synaptic activity. These alterations, combined with changes in the extracellular space, contribute to synaptic dysfunction and the occurrence of pathological high-frequency electrical activity or even epileptic activity \cite{cell_volume_seizures_ullah, chen2022poststroke}. 
Therefore, astrocytic malfunction following ischemia is recognized as a key mechanism underlying high-frequency activity and post-stroke epilepsy.
In order to introduce the biological response delay in the model for the volume of the cells at time $t$, we take into account the oxygen concentration in neurons and astrocytes at time $t-\tau$, $[O_{2}]_c(t-\tau)$, with $c \in \{n,a\}$ denoting the neuron and astrocyte compartments, respectively. The ECS compartment volume is modified accordingly.  

\noindent The coupling between electrophysiology and metabolic model describes the energetic needs of the cells, encoded via the energetic cost of the ion pump action, supplied by ATPases.
Collecting all the electrophysiology variables into the vector 
\begin{equation}
    \mathbf{m}(t) = 
    \begin{bmatrix}
    \ [\text{Na}^+]_i(t), \ [\text{K}^+]_o(t), \ [\text{Ca}^{2+}]_i(t), \ m(t), \ n(t), \ h(t) \
    \end{bmatrix}
    \in \mathbb{R}^6
    \label{eq:vect_electro}
\end{equation}
we can write the model \eqref{eq:BC} as a non-linear system 
\begin{equation}
    \frac{d\mathbf{m}}{dt} + \boldsymbol{g}(\mathrm{u},\mathbf{m}, p_n, p_a, \xi) = \boldsymbol{0}.
    \label{eq:electro}
\end{equation}
Finally, by combining the models described by Equations \eqref{eq:metabolism} and \eqref{eq:electro}, we obtain a coupled model that captures both electrophysiological and metabolic dynamics.

\subsection{The monodomain model}
\label{sec:monodomain}
For the spatial propagation of neuronal activity, the cellular coupled model is embedded into the monodomain formulation \cite{schreiner2022simulating, leimer2024highorder}, which describes the evolution of the transmembrane potential across a tissue domain. The resulting problem couples a PDE for the potential with ODE systems describing ionic concentrations, gating variables, and metabolic states, given by Equations \eqref{eq:electro} and \eqref{eq:metabolism}.
\par 
Given an open, bounded domain $\Omega \in \mathbb{R}^d$, $(d=2,3)$, and a final time $T>0$, we introduce the transmembrane potential $\mathrm{u} = \mathrm{u}(\boldsymbol{x},t)$ with $\mathrm{u}: \Omega \times [0,T] \rightarrow \mathbb{R}$, the vector $\mathbf{m} = \mathbf{m}(\boldsymbol{x},t)$  with $\mathbf{m}: \Omega \times [0,T] \rightarrow \mathbb{R}^n, n=6,$  containing the ion concentrations and gating variables of the ionic model, and finally the vector $\mathbf{c} = \mathbf{c}(\boldsymbol{x},t)$  with $\mathbf{c}: \Omega \times [0,T] \rightarrow \mathbb{R}^m, m=12,$ containing oxygen concentrations and variables related to blood and energy dynamics. The coupled problem reads as follows: 
\par
For any time $ t \in (0,T]$, find $\mathrm{u}=\mathrm{u}(\boldsymbol{x},t)$, $\mathbf{m}=\mathbf{m}(\boldsymbol{x},t)$ and $\mathbf{c}=\mathbf{c}(\boldsymbol{x},t)$ such that:
\begin{subequations}
\label{eq:monodomain}

\begin{empheq}[left=\empheqlbrace]{alignat=4}
    \label{eq:monodomain:pde}
    \chi_m C_m  \frac{\partial \mathrm{u}}{\partial t} - \nabla \cdot (\mathbf{\Sigma} \nabla  \mathrm{u}) +  \chi_m f(\mathrm{u},\mathbf{m},\mathbf{c}) & = 0, & & \quad \mathrm{in} \; \Omega \times (0,T],
    \\
    \label{eq:monodomain:ionic}
    \frac{\partial \boldsymbol{\mathbf{m}}}{\partial t} + \boldsymbol{g}(\mathrm{u},\mathbf{m},p_n,p_a,\xi) & = \boldsymbol{0}, & & \quad \mathrm{in} \; \Omega \times (0,T],    
    \\
    \label{eq:monodomain:metab}
    \frac{\partial \boldsymbol{\mathbf{c}}}{\partial t} + \boldsymbol{l}(\mathrm{u},\mathbf{c},q_a,\psi_\text{ATPase,n},\psi_\text{ATPase,a},q) & = \boldsymbol{0}, & & \quad \mathrm{in} \; \Omega \times (0,T],  
    \\
    \label{eq:monodomain:bcs}
    \mathbf{\Sigma} \nabla \mathrm{u} \cdot \boldsymbol{n} & = 0,  & & \quad\mathrm{on}\;  \partial \Omega  \times (0,T],    
    \\
    \label{eq:monodomain:ics}
    \mathrm{u}(0) = \mathrm{u}^0, \; \mathbf{m}(0)  = \mathbf{m}^0, \; \mathbf{c}(0) & = \mathbf{c}^0,  & & \quad \mathrm{in}\; \Omega.
\end{empheq}

\end{subequations}

In Equation~\eqref{eq:monodomain:pde}, $\boldsymbol{\Sigma}$ is the conductivity tensor, $\chi_m$ is the membrane capacitance per unit area, $C_m$ is the membrane capacitance, and $f=f(\mathrm{u},\mathbf{m},\mathbf{c})$ represents the ionic forces. In Equation~\eqref{eq:monodomain:ionic}, $\boldsymbol{g}=\boldsymbol{g}(\mathrm{u},\mathbf{m},\mathbf{c})$ represents the evolution of the ion concentrations and in Equation~\eqref{eq:monodomain:metab} $\boldsymbol{l}=\boldsymbol{l}(\mathrm{u},\mathbf{m},\mathbf{c})$ represents the evolution of concentration related to blood and oxygen dynamics. We impose homogeneous Neumann boundary conditions in Equation~\eqref{eq:monodomain:bcs}, where $\boldsymbol{n}$ represents the normal to the boundary $\partial \Omega$. Finally, we enforce the initial conditions $\mathrm{u}^0$, $\mathbf{m}^0$ and $\mathbf{c}^0$ in Equation~\eqref{eq:monodomain:ics}.

The propagation of electrical signals in brain tissue depends on its electrical conductivity, which can be influenced by micro structural variations. 
A well-established relationship links the apparent diffusion coefficient $\boldsymbol{D}^*$ of molecules in the extracellular space to the geometric property known as tortuosity, defined as $\lambda = \lambda(\boldsymbol{x})$ \cite{nicholson1998diffusion_tortuosity}. 
This quantity depends on the type of brain tissue being modeled. In the present work, we focus exclusively on ischemic regions within the gray matter and the tortuosity is assumed to be constant within each such region. The relation reads as 
\begin{equation}
    \boldsymbol{D}^* = \frac{1}{\lambda^2}\boldsymbol{D},
\label{eq:diffusion_tortuosity}
\end{equation}
where $\boldsymbol{D}$ is the free diffusion coefficient in an unobstructed medium, and $\lambda$ quantifies the increase in path length caused by the complex geometry of the ECS and cellular volume variation.
This relation, widely validated for the diffusion of metabolites and ions, shows that an increase in tortuosity reduces the apparent diffusion \cite{nicholson1998diffusion_tortuosity}. We recall that the conductivity tensor $\boldsymbol{\Sigma}$ can be described via an isotropic part, that is exploited for the modelization of grey matter, and an anisotropic part where we take into account axonal directions, fundamental for the representation of white matter tissue:
\begin{equation}
    \boldsymbol{\Sigma} = \sigma_\text{ext} \mathds{1}  + \sigma_\text{axn} \boldsymbol{n} \otimes \boldsymbol{n}.   
\end{equation}
By analogy, the conductivity tensor can be scaled with the same geometric factor $\lambda^2$. The conductivity in the ischemic region is therefore defined as:
\begin{equation}
    \boldsymbol{\Sigma}^* = \frac{1}{\lambda^2}\,\boldsymbol{\Sigma} = \frac{1}{\lambda^2}\sigma_\text{ext} \mathds{1}  + \frac{1}{\lambda^2}\sigma_\text{axn} \boldsymbol{n} \otimes \boldsymbol{n},  
    \label{eq:conductivity_lambda}
\end{equation}
which preserves the same principal directions defining anisotropy. Under pathological ischemic condition, the value of the parameter $\lambda$ increases to $\lambda = 2.2$ (see \cite{nicholson1998diffusion_tortuosity}). 

%% file: PolyDG.tex
\label{sec:polyDG}
To derive the weak formulation of the problem, we introduce the Sobolev space $V=H^1(\Omega)$, and we employ a standard definition of the scalar product in $L^2(\Omega)$, denoted by $(\cdot,\cdot)_{\Omega}$. The induced norm is denoted by $\lVert \cdot \rVert$. 
We assume that the forcing terms, physical parameters, and initial conditions are sufficiently regular to have a well-defined formulation, $i.e.$: $f(\mathrm{u},\mathbf{m},\mathbf{c}) \in L^2(0,T;L^2(\Omega))$, $I^\mathrm{ext}(\mathrm{u},\mathbf{m},\mathbf{c}) \in L^2(0,T;L^2(\Omega))$, $\boldsymbol{g}(\mathrm{u},\mathbf{m},\mathbf{c}) \in [L^2(0,T;L^2(\Omega))]^n$, $\boldsymbol{l}(\mathrm{u},\mathbf{m},\mathbf{c}) \in [L^2(0,T;L^2(\Omega))]^m$, $\chi_m$ and $C_m \in L_+^\infty(\Omega)$ where $L^{\infty}_+(\Omega):=\{v \in L^\infty(\Omega): v \ge 0 \text{ a.e. in } \Omega\}$, $\mathrm{u}^0\in L^2(\Omega)$,  $\mathbf{m}^0\in [L^2(\Omega)]^n$ and $\mathbf{c}^0\in [L^2(\Omega)]^m$.
The weak formulation of the problem \eqref{eq:monodomain} reads: $ \forall \: t \in (0,T] $ find $\mathrm{u}(t) \in V ,\mathbf{m}(t) \in \mathbf{V},\mathbf{c}(t) \in \mathbf{V}$ such that:
\begin{equation}
\left\{
\begin{aligned}
    \left(\chi_m C_m\frac{\partial \mathrm{u}(t)}{\partial t},v\right)_{\Omega}
    + (\boldsymbol{\Sigma} \nabla \mathrm{u},\nabla v)_{\Omega}
    + \left(\chi_m f\left(\mathrm{u}(t),\mathbf{m}(t), \mathbf{c}(t)\right),v\right)_\Omega
    &= (I^\text{ext},v)_\Omega
    && \forall\, v \in V, \\
    \left(\frac{\partial \mathbf{m}(t)}{\partial t},\boldsymbol{w}\right)_{\Omega}
    + \left(\boldsymbol{g}(\mathrm{u}(t),\mathbf{m}(t), \mathbf{c}(t)),\boldsymbol{w}\right)_{\Omega}
    &= \boldsymbol{0}
    && \forall\, \boldsymbol{w} \in [V]^n, \\
    \left(\frac{\partial \mathbf{c}(t)}{\partial t},\boldsymbol{z}\right)_{\Omega}
    + \left(\boldsymbol{l}(\mathrm{u}(t),\mathbf{m}(t), \mathbf{c}(t)),\boldsymbol{z}\right)_{\Omega}
    &= \boldsymbol{0}
    && \forall\, \boldsymbol{z} \in [V]^m, \\
    \mathrm{u}(0)=\mathrm{u}^0,\quad
    \mathbf{m}(0)=\mathbf{m}^0,\quad
    \mathbf{c}(0)&=\mathbf{c}^0
    && \text{in } \Omega.
\end{aligned}
\right.
\label{eq:weak_general}
\end{equation}

We now present the PolyDG semi-discrete formulation of the problem described in Equation \eqref{eq:weak_general}. Let $\mathcal{T}_h$ represent a polytopal mesh partition of the domain $\Omega$, consisting of disjoint elements $K$. For each element $K$, we define its diameter as $h_K$ and set $h = \max_{K \in \mathcal{F}_h} h_K < 1$. The interfaces are defined as the intersections of the $(d-1)$-dimensional facets of neighboring elements. We denote by $\mathcal{F}_h^I$ the union of all interior faces contained within $\Omega$ and by $\mathcal{F}_h^N$ those lying on the boundary $\partial \Omega$. In the following, we assume that the underlying grid is polytopic regular in the sense of \cite{cangiani_hp-version_2014,cangiani_hp-version_2017}.
\par
We define $\mathbb{P}^{p_K}(K)$ as the space of polynomials of degree $p_K\geq1$ over the element $K$ and the discontinuous finite element space as: 
\begin{equation*}
V_h^{\mathrm{DG}} = \{v_h \in L^2(\Omega) : v_h|_{K} \in \mathbb{P}^{p_K}(K) \quad \forall \: K \in \mathcal{T}_h\},
\end{equation*}
Let $F \in \mathcal{F}_h^I$ be the face shared by the elements $K^{\pm}$, and let $\boldsymbol{n}^{\pm}$ denote the normal unit vectors pointing outward to $K^\pm$, respectively. For a regular enough scalar-valued function $v$ and a vector-valued function $\boldsymbol{q}$, the trace operators are defined as follows \cite{arnold_unified_2002}:
\begin{equation*}
    \begin{aligned}
      &\averagel v \averager = \frac{1}{2} (v^+ + v^-) , \quad & \jumpl v  \jumpr& = v^+ \boldsymbol{n}^+ + v^- \boldsymbol{n}^-, \quad &\text{on }F \in \mathcal{F}_h^I,&\\
      &\averagel \boldsymbol{q} \averager =  \frac{1}{2} (\boldsymbol{q}^+ + \boldsymbol{q}^-) , \quad & \jumpl \boldsymbol{q} \jumpr & = \boldsymbol{q}^+\cdot \boldsymbol{n}^+ + \boldsymbol{q}^-\cdot \boldsymbol{n}^-,  \quad &\text{on }F \in \mathcal{F}_h^{I},&\\
    \end{aligned}
\end{equation*}
where the superscripts $\pm$ indicate the traces of these functions on $F$ taken in the interiors of $K^{\pm}$, respectively. The definition of the penalization parameter reads as follows $\eta : \mathcal{F}_h^I \cup \mathcal{F}_h^D \rightarrow \mathbb{R}_+$:
\begin{equation}
    \eta = \eta(\boldsymbol{p},h,\boldsymbol{\Sigma}) = \eta_0 
    \begin{cases}
        \{\boldsymbol{\Sigma}_K\}_A \dfrac{\{p^2_K\}_A}{\{h_K\}_H} & \text{on } F \in \mathcal{F}_h^I, \\
        \boldsymbol{\Sigma}_K \dfrac{p^2_K}{h_K} & \text{on } F \in \mathcal{F}_h^D, \\
    \end{cases}
    \label{eq:eta}
\end{equation}
which depends explicitly on both the local degrees and the mesh size. The PolyDG framework is extended in order to exploit adaptivity with respect to the polynomial degree, to reduce the overall computational costs and the total number of degrees of freedom of the system. For further details the exploited $p$-adaptive algorithm is described in \cite{leimer2025p}, in which the polynomial degree vector $\boldsymbol{p} = \{p_K\}_{K \in \mathcal{T}_h}$ is dynamically updated. The penalization parameter $\eta$ is modified accordingly, while $\eta_0$ is selected large enough to ensure stability of the method. 
In Equation \eqref{eq:eta} we set $\boldsymbol{\Sigma}_K = \|\sqrt{\boldsymbol{\Sigma}}|_K\|^2_{L^2(K)}$ and we consider both the harmonic average operator $\{\cdot\}_H$, and the arithmetic average operator $\{\cdot\}_A$ on $F \in \mathcal{F}_h^I$.
This setting enable us to introduce the following bilinear form $\mathcal{A}(\cdot,\cdot): V_h^{\mathrm{DG}}\times V_h^{\mathrm{DG}} \rightarrow \mathbb{R}$:
\begin{equation}
    \begin{aligned}
       \mathcal{A}(u,v) = \int_{\Omega} \mathbf{\Sigma}\nabla_h u \cdot  \nabla_h v \;dx+ \sum_{F \in \mathcal{F}_h^I} \int_F (\eta  \jumpl u \jumpr \cdot  \jumpl v  \jumpr - \averagel \boldsymbol{\Sigma} \nabla u \averager \cdot  \jumpl v\jumpr   - \jumpl u \jumpr \cdot \averagel \boldsymbol{\Sigma} \nabla v\averager) d\sigma \quad \forall \: u,v \in V^{\text{DG}}_h ,
    \end{aligned}
    \label{eq:coer}
\end{equation}
where $\nabla_h$ is the element-wise gradient. The semi-discrete formulation of problem in Equation \eqref{eq:monodomain} reads:
\par
For any $t \in (0,T]$, find $(\mathrm{u}_h(t),\textbf{m}_h(t), \textbf{c}_h(t)) \in V^{\mathrm{DG}}_h \times \mathbf{V}^{\mathrm{DG}}_h \times \mathbf{V}^{\mathrm{DG}}_h$ such that:
\begin{equation}
\left\{
\begin{aligned}
\left(\chi_m C_m  \dfrac{\partial \mathrm{u}_h(t)}{\partial t}, v_h \right)_\Omega 
+ \mathcal{A}(\mathrm{u}_h(t),v_h)
+ \left(\chi_m f\!\left(\mathrm{u}_h(t),\mathbf{m}_h(t),\mathbf{c}_h(t)\right),v_h\right)_\Omega
&= (I^{\text{ext}}_h, v_h)_\Omega 
&& \forall v_h \in V_h^{DG}, \\[0.4em]
\left( \dfrac{\partial \mathbf{m}_h(t)}{\partial t}, \boldsymbol{w}_h \right)_\Omega 
+ \left(\boldsymbol{g}\!\left(\mathrm{u}_h(t),\mathbf{m}_h(t), \mathbf{c}_h(t)\right),\boldsymbol{w}_h\right)_{\Omega}
&= 0
&& \forall \boldsymbol{w}_h \in [V_h^{DG}]^n, \\[0.4em]
\left( \dfrac{\partial \mathbf{c}_h(t)}{\partial t}, \boldsymbol{z}_h \right)_\Omega 
+ \left(\boldsymbol{l}\!\left(\mathrm{u}_h(t),\mathbf{m}_h(t), \mathbf{c}_h(t)\right),\boldsymbol{z}_h\right)_{\Omega}
&= 0
&& \forall \boldsymbol{z}_h \in [V_h^{DG}]^m, \\[0.4em]
\mathrm{u}_h(0) = \mathrm{u}_h^0,\quad
\mathbf{m}_h(0) = \mathbf{m}_h^0,\quad
\mathbf{c}_h(0) &= \mathbf{c}_h^0
&& \text{in } \Omega.
\end{aligned}
\right.
\label{eq:mono_semi_discrete}
\end{equation}

We denote the dimension of the discrete space as $N_h(\boldsymbol{p})$, to make explicit its dependence on the vector of local polynomial distribution. Let $N_h(\boldsymbol{p})$ be the dimension of $V_h^\mathrm{DG}$ and let $(\varphi_j)^{N_h(\boldsymbol{p})}_{j=0}$ be a suitable basis for $V_{h}^\mathrm{DG}$, then $\mathrm{u}_h(t) = \sum_{j=0}^{N_h(\boldsymbol{p})} \text{U}_j(t)\varphi_j$, $\text{m}_l(t) = \sum_{j=0}^{N_h(\boldsymbol{p})} \text{M}^l_j(t)\varphi_j$ for all $l = 1,...,n$, and $\text{c}_q(t) = \sum_{j=0}^{N_h(\boldsymbol{p})} \text{C}^q_j(t)\varphi_j$ for all $q=1,...,m$. We denote $\mathbf{U} \in \mathbb{R}^{N_h(\boldsymbol{p})}$, $\mathbf{M}_l \in \mathbb{R}^{N_h(\boldsymbol{p})}$  for all $l=1,...,n$ and $\mathbf{C}_q \in \mathbb{R}^{N_h(\boldsymbol{p})}$ for all $q=1,...,m$. Moreover we define $\mathbf{M} = [\mathbf{M}_1,...,\mathbf{M}_n]^\top$ and $\mathbf{C} = [\mathbf{C}_1,...,\mathbf{C}_m]^\top$. We define the matrices:
\begin{equation}
    \begin{aligned}
      [\mathbf{M}_\text{prj}]_{ij} &= (\varphi_i,\varphi_j)_{\Omega}, \; &\text{(Mass matrix),}& \quad  i,j = 1,...,N_h(\boldsymbol{p}) \\
      [\mathbf{F}]_{j} &= (I^\text{ext},\varphi_j)_{\Omega},  \; &\text{(Forcing term),}& \quad j = 1,...,N_h(\boldsymbol{p})\\
      [\mathbf{I}(\mathrm{u},\mathbf{m},\mathbf{c})]_{j} &= (f(\mathrm{u},\mathbf{m},\mathbf{c}),\varphi_j)_{\Omega},  \; &\text{(Non-linear ionic forcing term),}& \quad  j = 1,...,N_h(\boldsymbol{p})\\
[\mathbf{G}_l(\mathrm{u},\mathbf{m},\mathbf{c})]_{j} &= (\boldsymbol{g}_l(\mathrm{u},\mathbf{m},\mathbf{c}),\boldsymbol{\varphi}_j)_{\Omega},  \; &\text{(Electrophysiological model),}& \quad  j = 1,...,N_h(\boldsymbol{p}), \;l=1,...,n  \\ 
[\mathbf{L}_q(\mathrm{u},\mathbf{m},\mathbf{c})]_{j} &= (\boldsymbol{l}_q(\mathrm{u},\mathbf{m},\mathbf{c}),\boldsymbol{\varphi}_j)_{\Omega},  \; &\text{(Metabolic model),}& \quad  j = 1,...,N_h(\boldsymbol{p}), \;q=1,...,m  \\ 
      [\mathbf{A}]_{ij} &= \mathcal{A}(\varphi_i,\varphi_j) &\text{(Stiffness matrix),}& \quad  i,j = 1,...,N_h(\boldsymbol{p}).
      \label{eq::matrixFull}
    \end{aligned}
\end{equation}
We are now ready to introduce the fully-discrete formulation.
We partition the interval \([0, T]\) into \(N\) sub-intervals \((t^{(k)}, t^{(k+1)}]\), each of length \(\Delta t\), such that \(t^{(k)} = k\Delta t\) for \(k = 0, \dots, N-1\). For temporal discretization, we adopt the second-order Crank-Nicolson scheme for the linear part, with the ionic term discretized using a second-order explicit extrapolation. Given the initial conditions \(\mathbf{U}_0\), \(\mathbf{M}_0\), and \(\mathbf{C}_0\) the discrete scheme is: find $\mathbf{U}^{(k+1)} \simeq \mathbf{U}(t^{(k+1)})$, $\mathbf{M}^{(k+1)} \simeq \mathbf{M}(t^{(k+1)})$ and $\mathbf{C}^{(k+1)} \simeq \mathbf{C}(t^{(k+1)})$ for $k=0,...,N-1$, such that
\begin{equation*}
\label{eq:Full_discrete_complete}
    \begin{aligned}
      \left(\chi_m C_m \mathbf{M}_\text{prj} + \frac{\Delta t}{2}\mathbf{A} \right) \mathbf{U}^{(k+1)}  &=  \,  \left(\chi_m C_m \mathbf{M}_\text{prj} - \frac{\Delta t}{2}\mathbf{A} \right) \mathbf{U}^{(k)}  - \frac{\mathrm{\chi_m} \Delta t}{2} (3\mathbf{I}^{(k)}-\mathbf{I}^{(k-1)}) + \frac{\Delta t}{2}( \mathbf{F}^{(k+1)}+\mathbf{F}^{(k)}) , \\
      \mathbf{M}^{(k+1)} &=  \, \mathbf{M}^{(k)} - \Delta t \mathbf{G}^{(k)},\\
     \mathbf{C}^{(k+1)} &=  \, \mathbf{C}^{(k)} - \Delta t \mathbf{L}^{(k)},\\
      (\mathbf{U}^{(0)},\mathbf{M}^{(0)}, \mathbf{C}^{(0)}) &= \, (\mathbf{U}_0,\mathbf{M}_0,\mathbf{C}_0).
    \end{aligned}
\end{equation*}

%% file: Numerical_result0D.tex
We present a set of numerical tests aimed at assessing the numerical performance of the method in approximating pathophysiological scenarios of brain function. We aim to investigate the interaction between metabolism and electrophysiology, as well as how localized ischemic regions influence the electrophysiology and the propagation of electrical signals in the surrounding tissue, providing novel insights into their coupled behaviour.

\subsection{Neuronal coupled model}
\label{sec:0D}
In this section, we investigate the effects of oxygen supply reduction on neuronal and astrocytic dynamics, 
highlighting how different levels of ischemia shape the electro-metabolic response of the system. 
In addition, we perform a sensitivity analysis with respect to the potassium clearance rate $\varepsilon$, 
which regulates the reactivity of the extracellular space to ionic perturbations, in order to assess how changes in this parameter influence the transition from physiological to pathological activity.

\subsubsection{Ischemic influence on pathological action potential evolution}
We first focus on a scenario in which we analyze how drops in blood flow can affect physiological neuronal evolution related to the transmembrane potential. 
From a clinical viewpoint, this drop affects the metabolic activity, triggering a rise in the potassium concentration in the ECS, increasing the excitability of the neuron. The piecewise smooth blood reduction function is described in Equation \eqref{eq:At_indicator}, taken from \cite{CALVETTI2018}.


\begin{equation}
\begin{aligned}
A(t) ={}&
\mathds{1}_{[t_0,\, t_1)}(t)
+ \left(1 - \delta \frac{t - t_1}{r_1}\right)\mathds{1}_{[t_1,\, t_1 + r_1)}(t)
+ (1 - \delta)\mathds{1}_{[t_1 + r_1,\, t_2)}(t) \\
&+ \left(1 - \delta \left(1 - \frac{t - t_2}{r_2}\right)\right)\mathds{1}_{[t_2,\, t_2 + r_2)}(t)
+ \mathds{1}_{[t_2 + r_2,\, T)}(t).
\end{aligned}
\label{eq:At_indicator}
\end{equation}

\begin{figure}[h]
    \centering
\begin{subfigure}[b]{0.3\textwidth}
    \centering
\includegraphics[width=\textwidth]{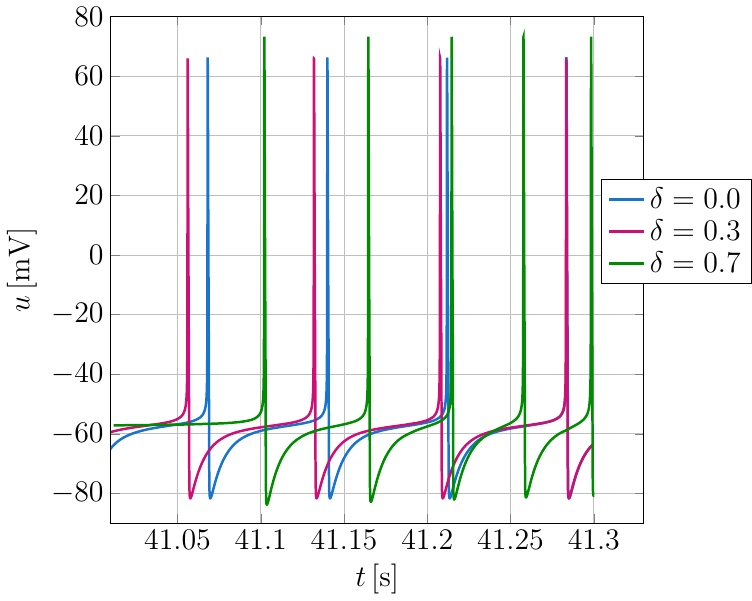}%
\end{subfigure}
\begin{subfigure}[b]{0.3\textwidth}
    \centering
\includegraphics[width=\textwidth]{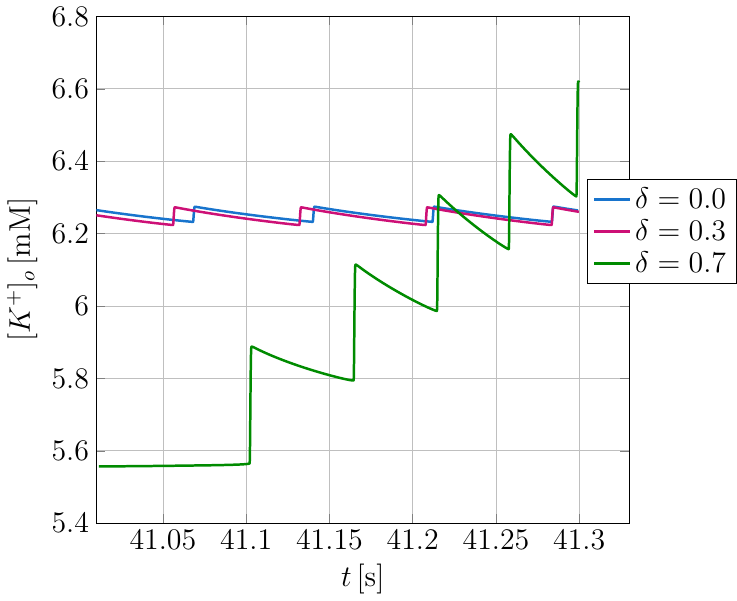}%
\end{subfigure}
\begin{subfigure}[b]{0.3\textwidth}
    \centering
\includegraphics[width=\textwidth]{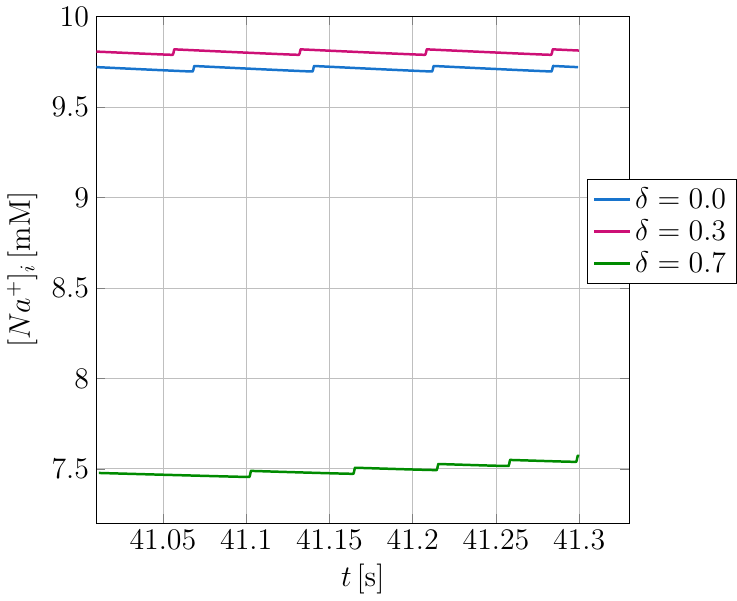}%
\end{subfigure}
\begin{subfigure}[b]{0.3\textwidth}
    \centering
\includegraphics[width=\textwidth]{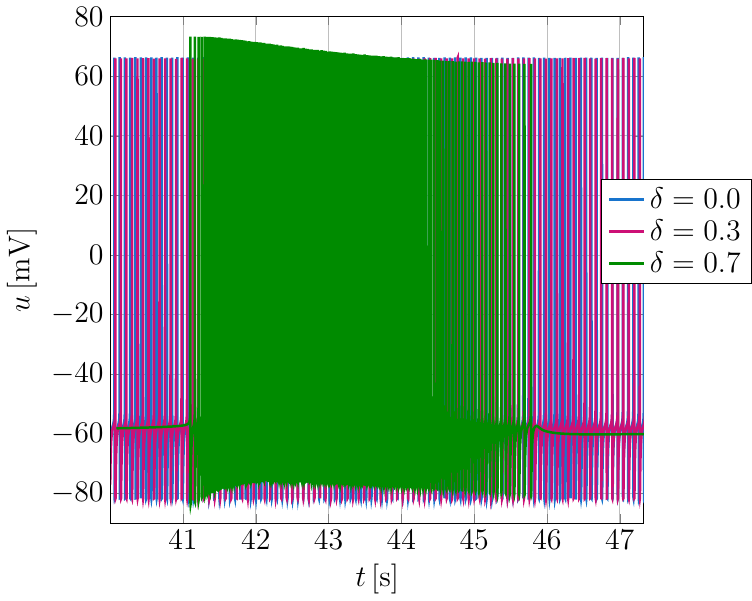}%
\end{subfigure}
\begin{subfigure}[b]{0.3\textwidth}
    \centering
\includegraphics[width=\textwidth]{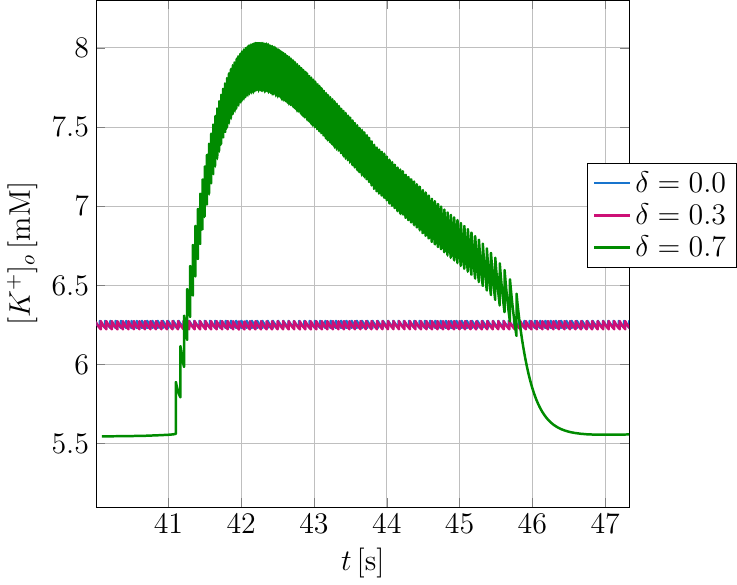}%
\end{subfigure}
\begin{subfigure}[b]{0.3\textwidth}
    \centering
\includegraphics[width=\textwidth]{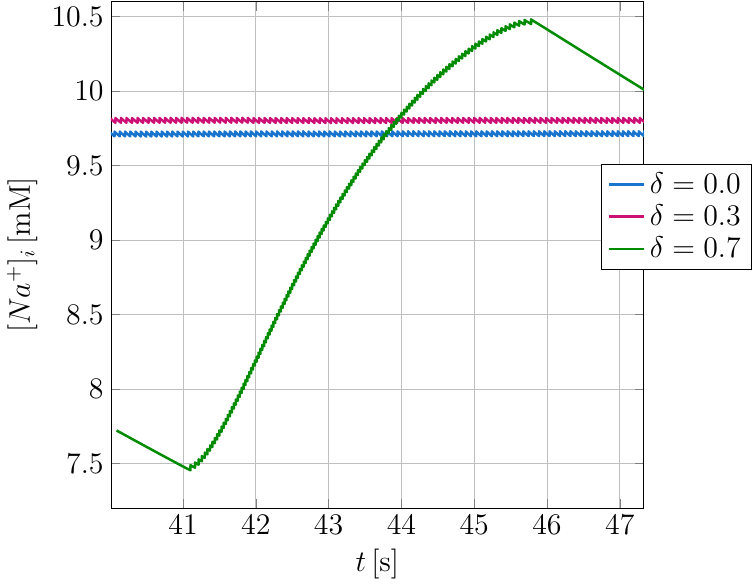}%
\end{subfigure}
\caption{Sensitivity to ischemic severity \(\delta\) for temporal zoom \(t\in[41.0,\,41.35]\) s (first row) and extended window \(t\in[40,\,47]\)s (second row).
Left: transmembrane potential \(\mathrm{u}\). Middle: extracellular potassium \(\text{[K$^+$]}_o\). Right: intracellular sodium \(\text{[Na$^+$]}_i\).}
    \label{fig:simulation_ischemia_short}
\end{figure}

\begin{figure}[H]
    \centering
\begin{subfigure}[b]{0.24\textwidth}
    \centering
\includegraphics[width=\textwidth]{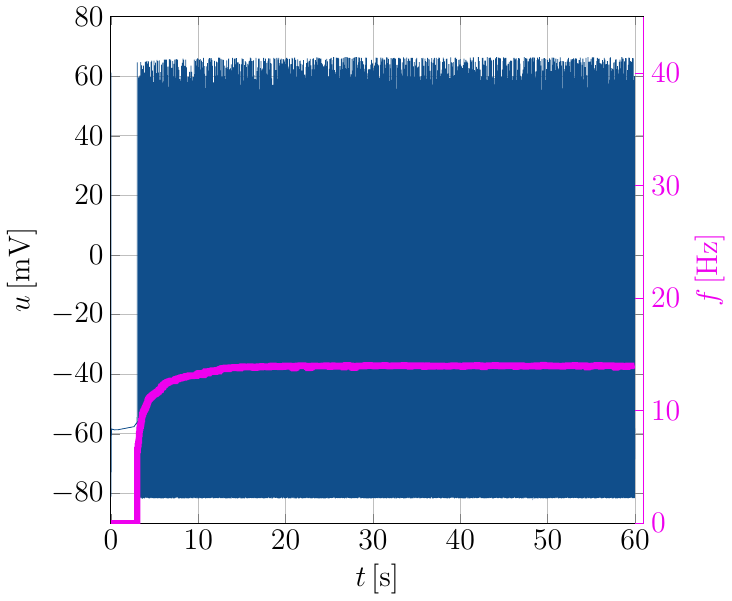}%
\end{subfigure}
\begin{subfigure}[b]{0.24\textwidth}
    \centering
\includegraphics[width=\textwidth]{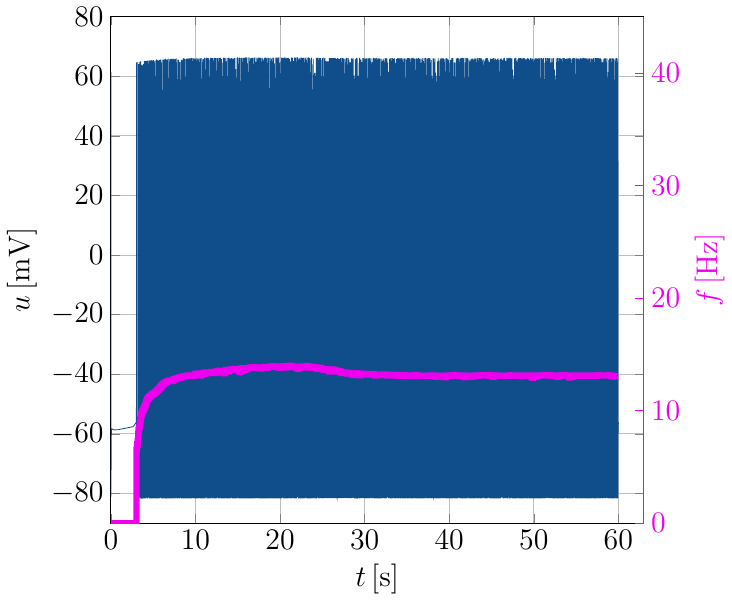}%
\end{subfigure}
\begin{subfigure}[b]{0.24\textwidth}
    \centering
\includegraphics[width=\textwidth]{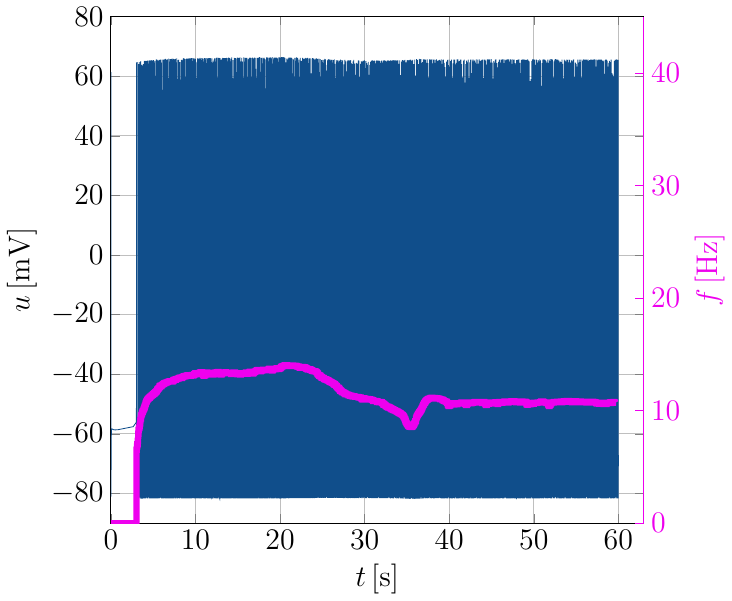}%
\end{subfigure}
\begin{subfigure}[b]{0.24\textwidth}
    \centering
\includegraphics[width=\textwidth]{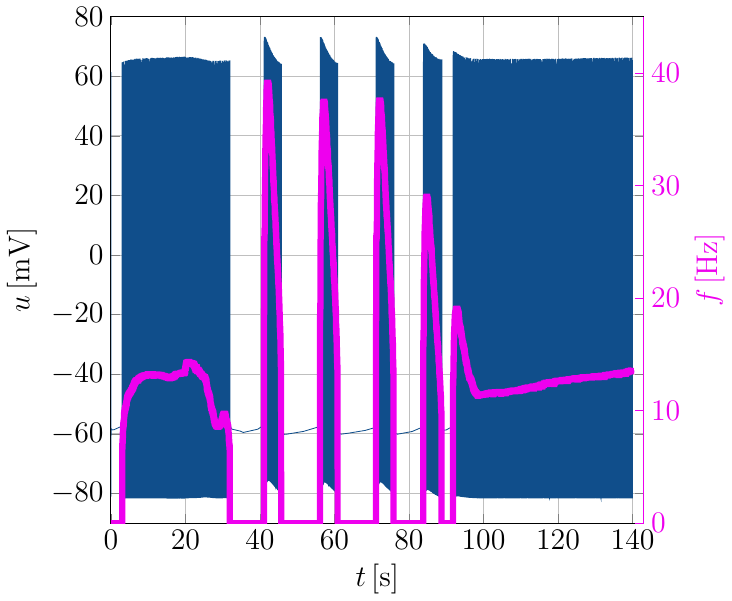}%
\end{subfigure}
\begin{subfigure}[b]{0.24\textwidth}
    \centering
\includegraphics[width=1.05\textwidth]{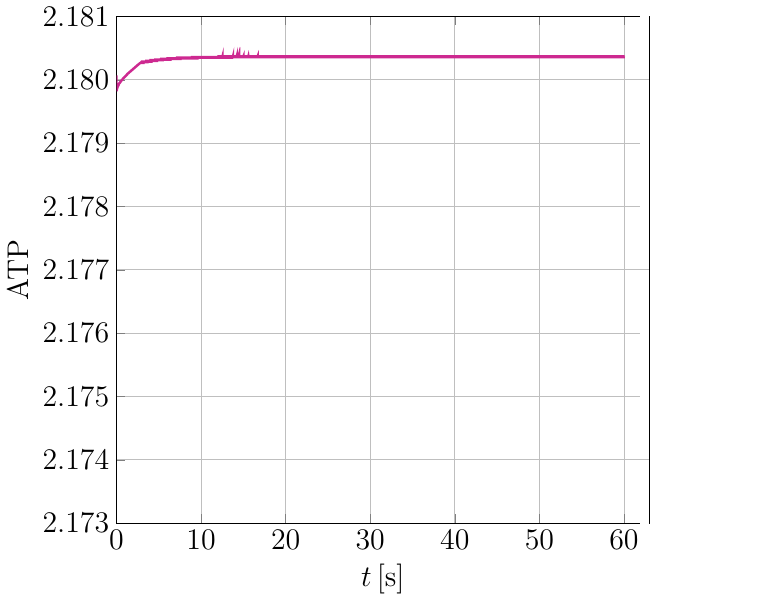}%
\end{subfigure}
\begin{subfigure}[b]{0.24\textwidth}
    \centering
\includegraphics[width=1.05\textwidth]{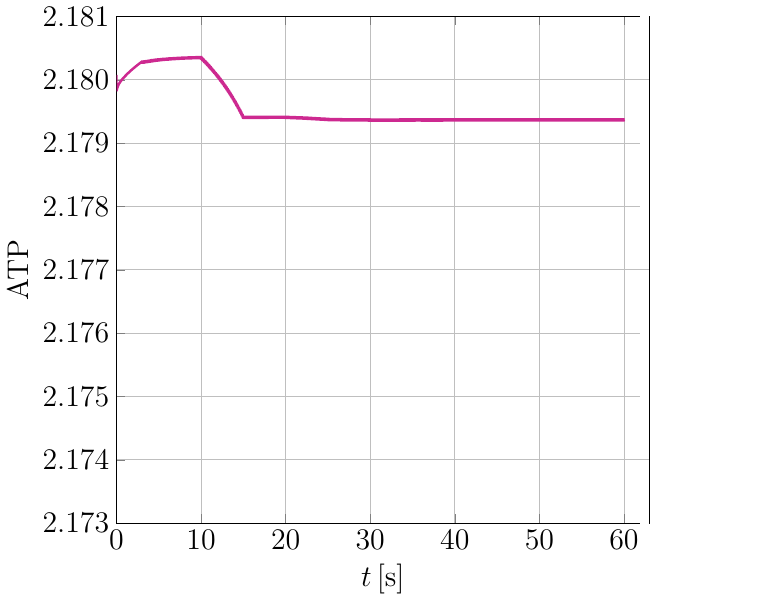}%
\end{subfigure}
\begin{subfigure}[b]{0.24\textwidth}
    \centering
\includegraphics[width=1.05\textwidth]{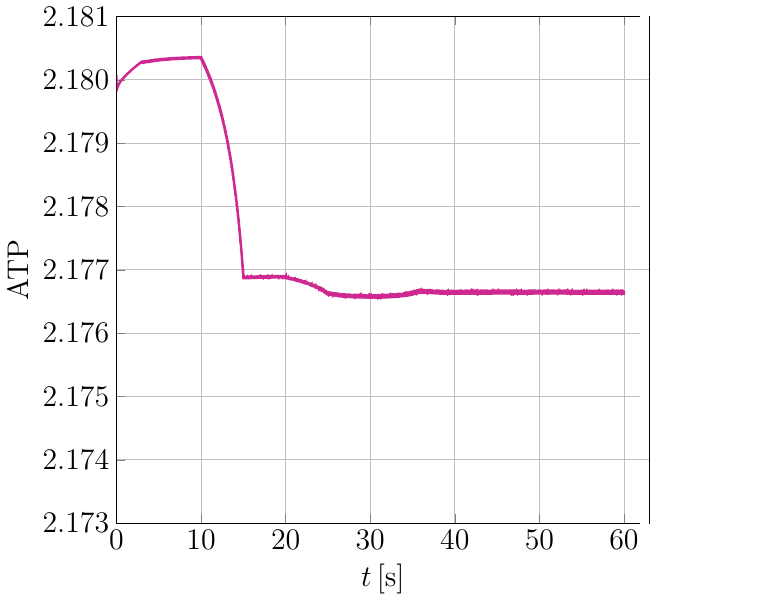}%
\end{subfigure}
\begin{subfigure}[b]{0.24\textwidth}
    \centering
\includegraphics[width=1.05\textwidth]{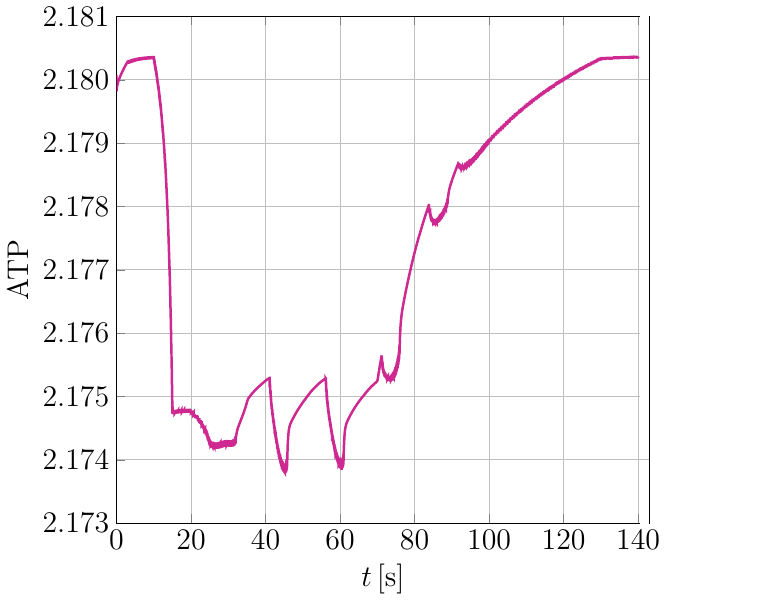}%
\end{subfigure}
\begin{subfigure}[b]{0.24\textwidth}
    \centering
\includegraphics[width=0.97\textwidth]{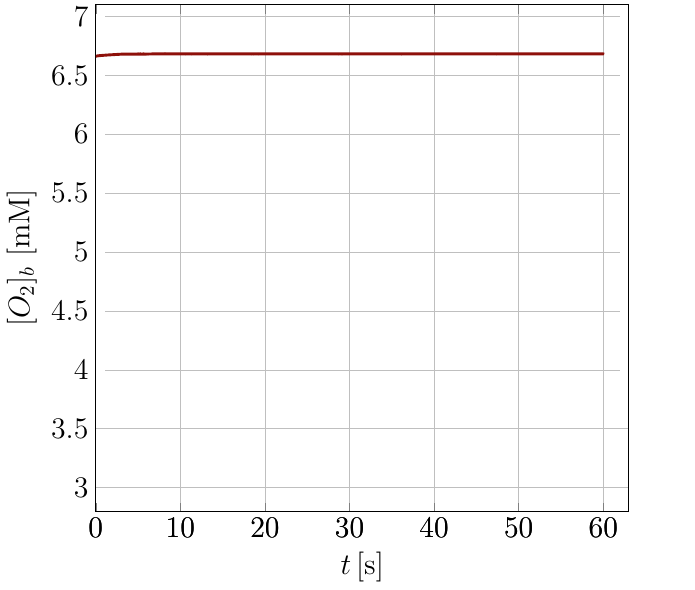}%
\end{subfigure}
\begin{subfigure}[b]{0.24\textwidth}
    \centering
\includegraphics[width=0.97\textwidth]{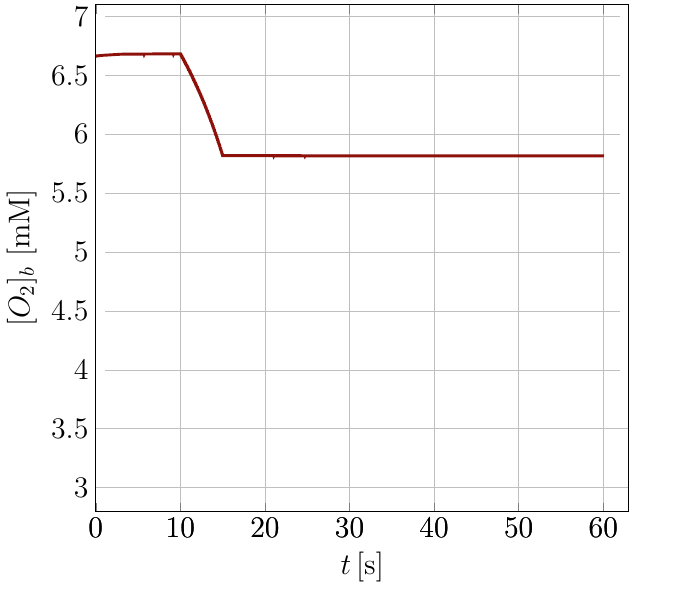}%
\end{subfigure}
\begin{subfigure}[b]{0.24\textwidth}
    \centering
\includegraphics[width=0.97\textwidth]{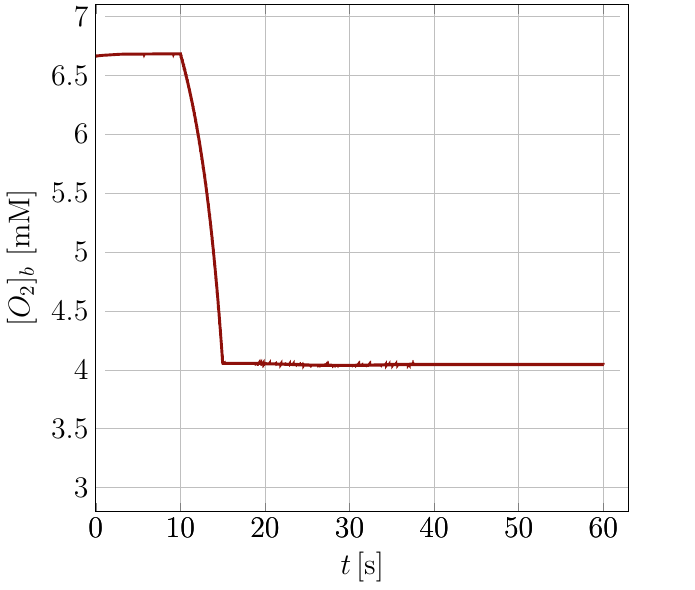}%
\end{subfigure}
\begin{subfigure}[b]{0.24\textwidth}
    \centering
\includegraphics[width=0.97\textwidth]{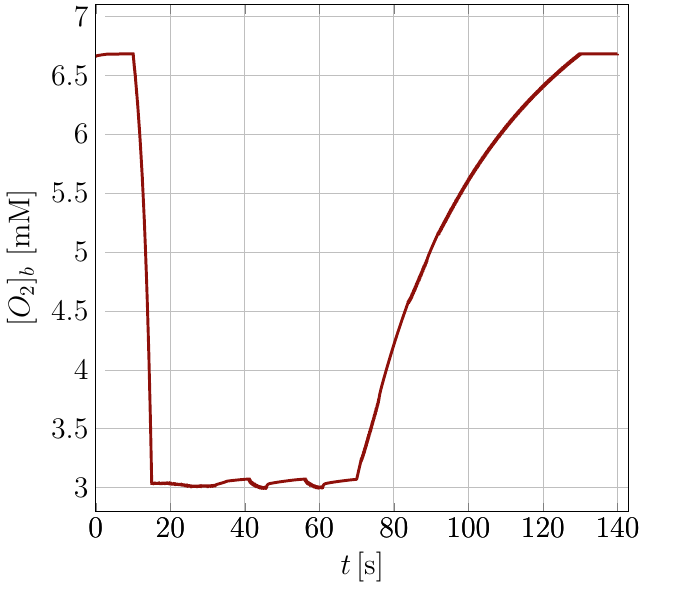}%
\end{subfigure}
\begin{subfigure}[b]{0.24\textwidth}
    \centering
\includegraphics[width=1.03\textwidth]{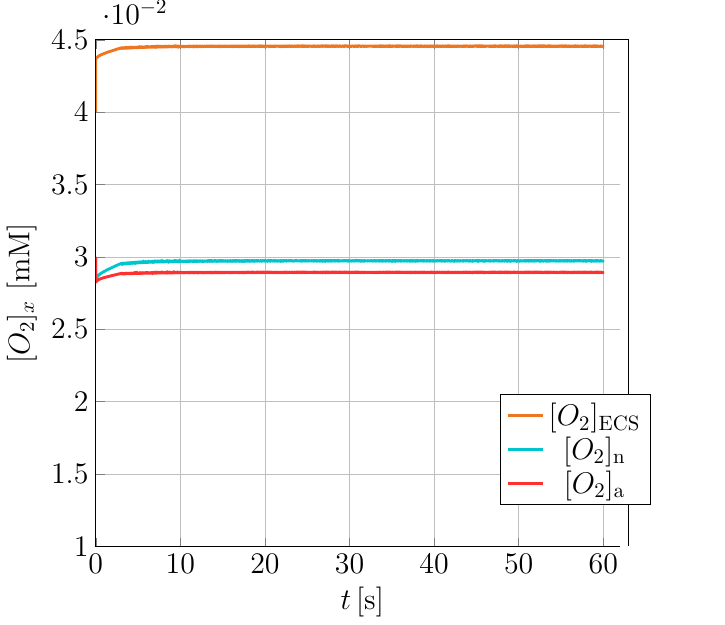}%
\end{subfigure}
\begin{subfigure}[b]{0.24\textwidth}
    \centering
\includegraphics[width=1.03\textwidth]{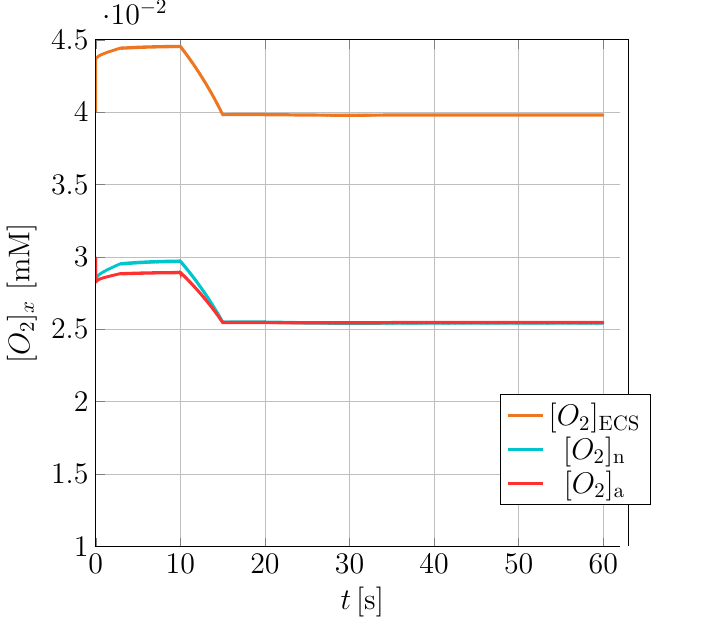}%
\end{subfigure}
\begin{subfigure}[b]{0.24\textwidth}
    \centering
\includegraphics[width=1.03\textwidth]{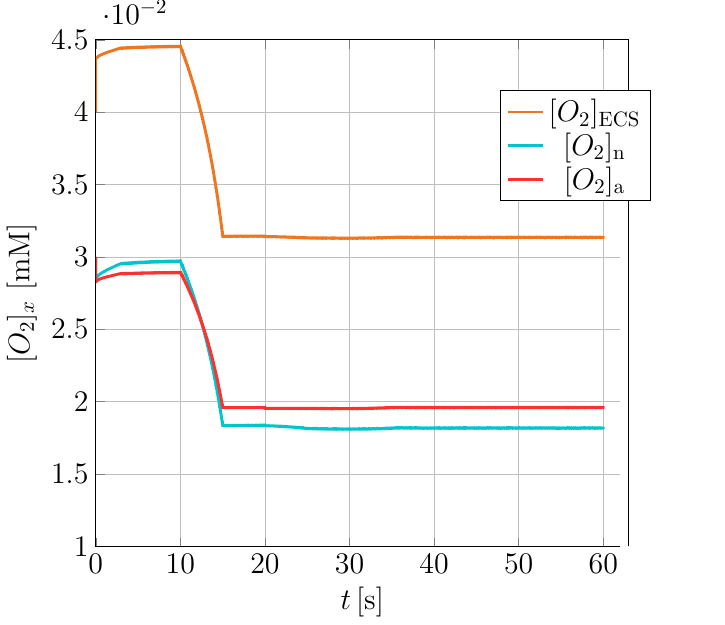}%
\end{subfigure}
\begin{subfigure}[b]{0.24\textwidth}
    \centering
\includegraphics[width=1.03\textwidth]{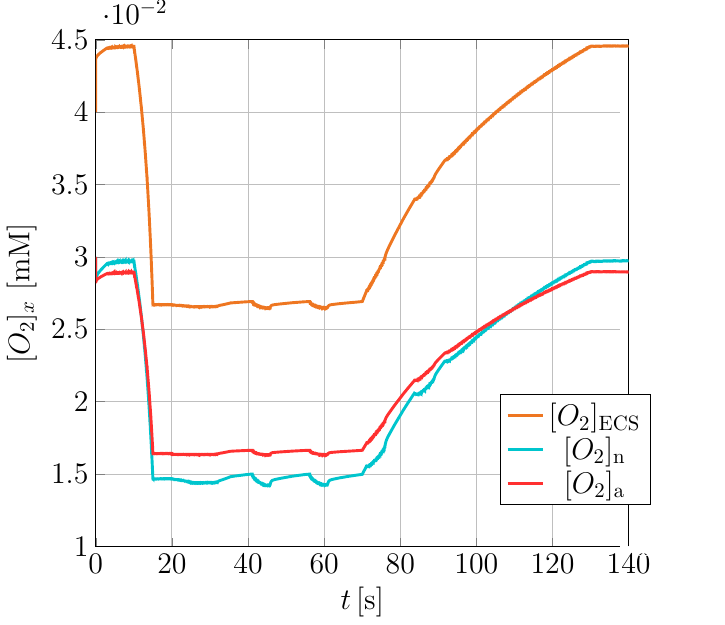}%
\end{subfigure}
\caption{Coupled electro–metabolic response at increasing ischemic severity \(\delta\).
Columns correspond to \(\delta=0\) (first column), \(\delta=0.3\) (second column, subclinical), \(\delta=0.6\) (third column, subclinical), and \(\delta=0.7\) (right, severe). First row: transmembrane potential \(\mathrm{u}\) (blue, left axis) and instantaneous firing frequency \(f\) (right axis). Second row: ATP concentration. Third row: blood oxygen concentration \([{\rm O}_2]_b\).
Fourth row: oxygen concentration in extracellular space (ECS), neuron, and astrocyte.}
    \label{fig:simulation_ischemia_20}
\end{figure}

We analyze the effect of different percentage reductions in blood volume compared to the baseline value. Based on the blood depression trend, we model a $30\%$, $60\%$ and $70\%$ decrease in blood concentration, reproducing both a subclinical ($30\%-60\%$) and severe ($70\%$) ischemia.
The parameters for the blood flow trend are listed in Table \ref{tab:ischemia_params_30}. 
Since there is no significant variation in the electro-metabolic cerebral system after cellular swelling for subclinical ischemia, the simulation is performed in a reduced time interval, without waiting for the end of the ischemic episode. The initial conditions for both the metabolic and electrical models are listed in Table \ref{table:initial_conditions_0D}. Figure \ref{fig:simulation_ischemia_short} shows a comparison of different scenarios in a limited time window: we analyze the sensitivity of the neuronal model with respect to different ischemic severity levels $\delta$. The first row presents a temporal zoom 
($t \in [41,\,41.35]\,\mathrm{s}$), where the simulation with severe ischemia exhibits a higher firing frequency compared to the subclinical cases. In the second row we extend the time window ($t \in [40,\,47]\,\mathrm{s}$) in order to appreciate the complete bursting behavior of the severe case. Figure \ref{fig:simulation_ischemia_20} compares a subclinical ischemia (second and third columns: 30\%-60\% cerebral blood flow reduction with onset over 10--15\,s) with a severe ischemia (right column: 70\% reduction with prolonged depression), each column displaying the transmembrane potential and ionic concentrations, ATP/ADP dynamics, and oxygen levels in blood, extracellular space (ECS), neuron, and astrocyte. 
The subclinical case represents a transient episode of neurological dysfunction caused by the focal brain ischemia without tissue injury. A mild reduction in cerebral blood flow may not cause permanent neuronal damage due to the activation of compensatory mechanisms, such as cerebral autoregulation \cite{cipolla2009cerebral} and energy coupling, that can keep the system close to physiological condition. These processes help the brain to maintain an appropriate supply of oxygen and nutrients, even under conditions of reduced blood flow. The membrane potential shows only a slight, transient frequency modulation at onset, preserving a physiological pattern.

\begin{table}[h]
    \centering
    \footnotesize
    \setlength{\tabcolsep}{6pt}
    \renewcommand{\arraystretch}{1.15}
    \begin{tabular}{|l|c|c|c|c|c|c|}
        \hline
        &$ \delta$ & $t_1$ & $r_1$ & $t_2$ & $r_2$ & $T$ \\
        \hline
        Subclinical ischemia &0.3-0.6 & 10\,s & 5\,s & 90\,s & 80\,s & 60\,s \\
        Severe ischemia& 0.7 & 10\,s & 5\,s &70\,s & 60\,s & 150\,s \\
        \hline
    \end{tabular}
    \caption{Parameters for subclinical and severe ischemia simulations: blood percentage drop, ramping times, and final observation time.}
    \label{tab:ischemia_params_30}
\end{table}

\begin{table}[h]
\centering
\resizebox{\textwidth}{!}{%
\begin{tabular}{|lcc|lcc|lcc|lcc|}
\hline
\text{Variable} & \text{Value} & \text{Units} &
\text{Variable} & \text{Value} & \text{Units}&
\text{Variable} & \text{Value} & \text{Units}&
\text{Variable} & \text{Value} & \text{Units}\\
\hline
$\mathrm{u}^0$ & $-56.199$ & [mV] & $\mathrm{[Ca]}_i^0$ & $0$ & [mM] & $[O_{2}]^0_{a}$ & $0.03$ &[mM] & $\text{NADH}^0_n$ & $1.2e-3$ &[mM] \\
$\mathrm{[K]}_o^0$ & $6.2773$ & [mM] & $\mathrm{[Na]}_i^0$ & $11.56$ & [mM] & $\text{ATP}^0_n$ & $2.18$ & [mM] & $\text{NAD}^0_n$ & $0.03$ & [mM] \\
$m^0$ & $0.0936$ & & $[O_{2}]^0_\text{ECS}$ & $0.04$ &[mM] & $\text{ADP}^0_n$ & $6.3\cdot 10^{-3}$ & [mM] & $\text{NADH}^0_a$ & $1.2\cdot 10^{-3}$ &[mM] \\
$n^0$ & $0.1558$ && $[O_{2}]^0_{n}$ & $0.03$ & [mM] & $\text{ATP}^0_a$ & $2.17$ & [mM] & $\text{NAD}^0_a$ & $0.03$ & [mM] \\
$h^0$ & $0.9002$ && $[O_{2}]^0_{b}$ & $6.67$ & [mM] & $\text{ADP}^0_a$ & $0.03$ & [mM] & $ $ & $ $ & \\
\hline
\end{tabular}
}
\caption{Severe ischemia, initial conditions for the variables of the coupled ionic model. }
\label{table:initial_conditions_0D}
\end{table}

In the case of subclinical ischemia, corresponding to the reduction in blood flow that occurs within the 10–15\,s interval, the variation in frequency is minimal and shows a temporal shift due to the physiological delay of the cellular response. This variation, nevertheless, remains within the physiological range. Consequently, there is a reduction in oxygen concentration in all compartments, along with a small variation in the resting value of ATP concentration, which remains close to baseline. Compartmental $[\mathrm{O}_2]$ decreases but it is sufficient to sustain Na$^+$/K$^+$-ATPase activity and astrocytic K$^+$ uptake, preventing pathological accumulation of $[\mathrm{K}^+]_\text{o}$.  
Such a small and transient ischemic scenario is supported by clinical evidence, indicating that it causes no permanent brain injury, as compensatory mechanisms preserve tissue viability.
In the severe case, reduced oxygen delivery undergoes oxidative phosphorylation and thus ATP availability, lowering the phosphorylation ratios $p_n=[\mathrm{ATP}]_n/[\mathrm{ADP}]_n$ and $p_a=[\mathrm{ATP}]_a/[\mathrm{ADP}]_a$ that scale the pump and glial clearance currents, allowing  $[\mathrm{K}^+]_\text{o}$ to accumulate. The Nernst potentials shift, and the membrane potential enters irregular dynamics with high frequency behaviour (40 Hz) \cite{Jiruska2011_high_freq}, while ATP decays and ADP rises, leading to an energetic crisis and $[{O}_2]$ collapses across all compartments. Once neurons and astrocytes react to the reduced blood flow, cellular swelling occurs, accompanied by irregular firing patterns. The frequency of bursts indicates the presence of high neuronal activity \cite{Jiruska2011_high_freq}, which can also represent the onset of epileptiform scenarios, which may appear in subtle forms without overt convulsions. During these bursts, oxygen levels in all compartments drop, reflecting the increased metabolic demand needed to sustain the higher firing rates. The effect of oxygen deficit depends on both the severity and the duration of the ischemic episode.

\begin{figure}[h]
    \centering
    \begin{subfigure}[b]{0.24\textwidth}
    \centering
\includegraphics[width=1\textwidth]{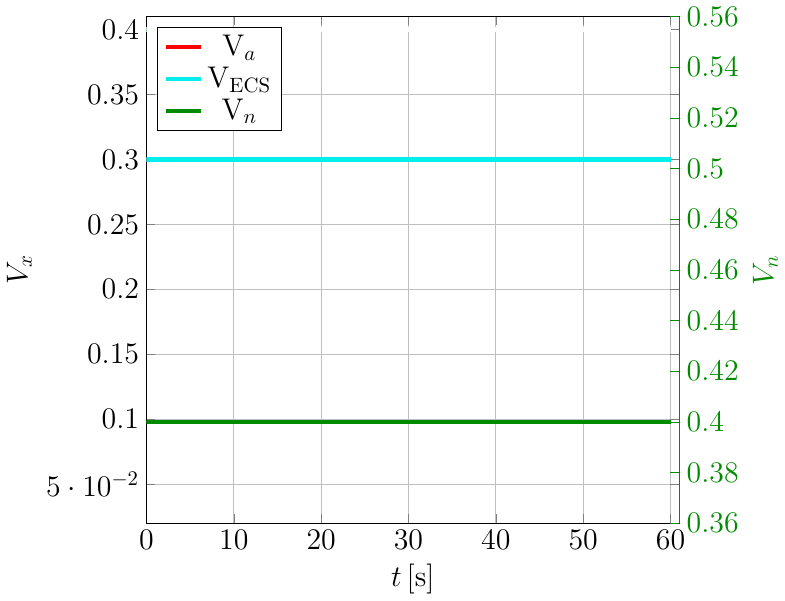}%
\caption{$\delta = 0.0$}
    \end{subfigure}
        \begin{subfigure}[b]{0.24\textwidth}
    \centering
\includegraphics[width=1\textwidth]{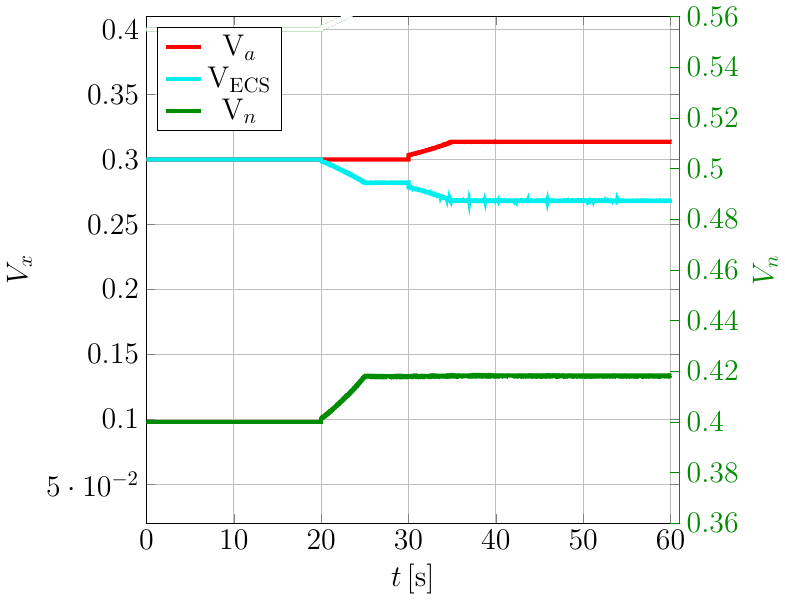}%
\caption{$\delta = 0.3$}
\end{subfigure}
\begin{subfigure}[b]{0.24\textwidth}
\centering
\includegraphics[width=1\textwidth]{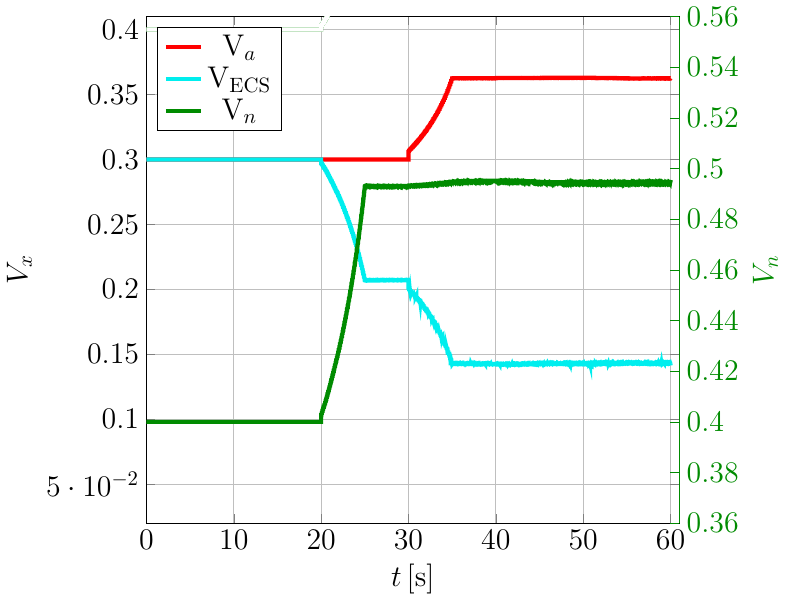}
\caption{$\delta = 0.6$}
\end{subfigure}
\begin{subfigure}[b]{0.24\textwidth}
\centering
\includegraphics[width=1\textwidth]{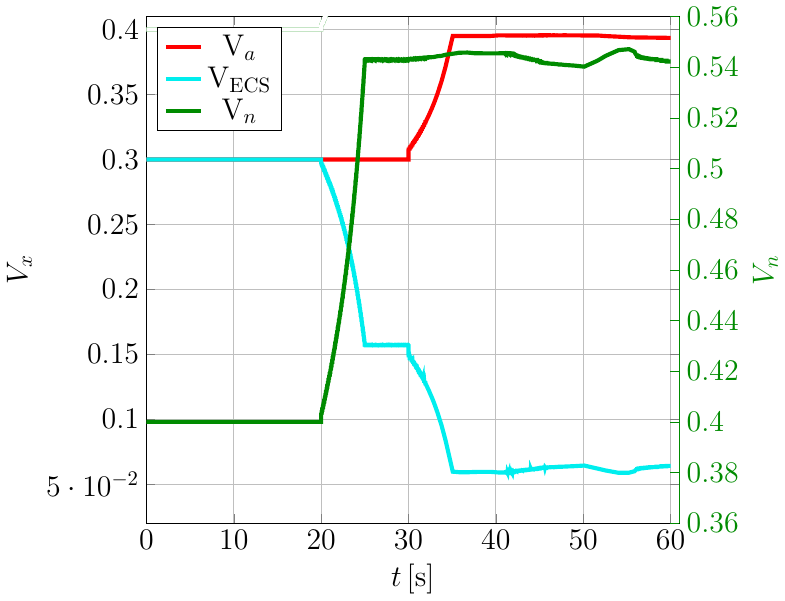}
\caption{$\delta = 0.7$}
\end{subfigure}
\caption{Volume evolution of neurons, astrocytes and extracellular space with respect to different values of $\delta$.}
    \label{fig:volume_comparison}
\end{figure}

Referring to the volumetric formulation \eqref{eq:volume_variation}, we further illustrate the  temporal evolution of cell volume variations under ischemic conditions. These results provide a quantitative view of how oxygen depletion and ionic imbalance drive swelling dynamics within neuronal and astrocytic compartments. 
As illustrated in Figure~\ref{fig:volume_comparison}, increasing values of $\delta$ lead to a progressive swelling of both neuronal and astrocytic compartments, accompanied by a marked reduction in the extracellular space volume. This behavior reflects the intracellular accumulation of ions and water under ischemic conditions, resulting in a decreased extracellular volume available for ionic diffusion and, consequently, a disruption of ionic homeostasis.  
Furthermore, a delay in the physiological response of neurons and astrocytes can be observed, reflecting their slower adaptation to ischemic event.
Higher values of $\delta$ therefore intensify cellular edema and reduce the tissue’s ability to maintain physiological equilibrium, contributing to increased neuronal excitability and ischemic damage. The evolution of the volumes for the different compartments are illustrated in Figure \ref{fig:volume_comparison}.

\subsubsection{Sensitivity analysis on potassium clearance rate}
\begin{figure}[h]
    \centering
    \begin{subfigure}[b]{0.23\textwidth}
    \centering
\hspace{-1.7em}\includegraphics[width=1\textwidth]{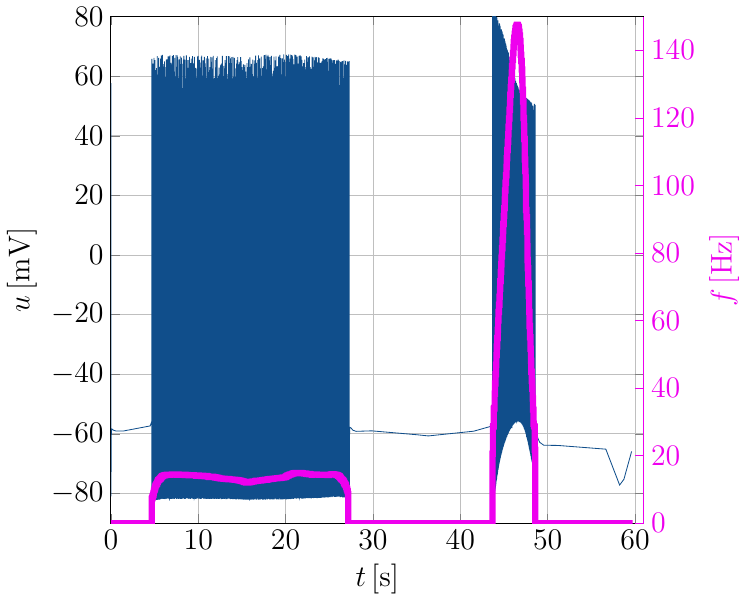}%
    \end{subfigure}
        \begin{subfigure}[b]{0.23\textwidth}
    \centering
\hspace{-1.4em}\includegraphics[width=1\textwidth]{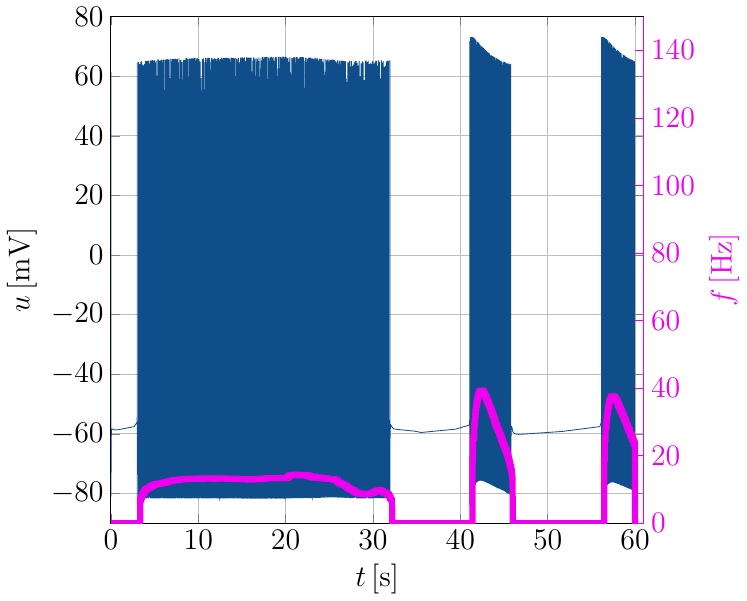}%
\end{subfigure}
    \begin{subfigure}[b]{0.23\textwidth}
\centering
\hspace{-0.6em}\includegraphics[width=1\textwidth]{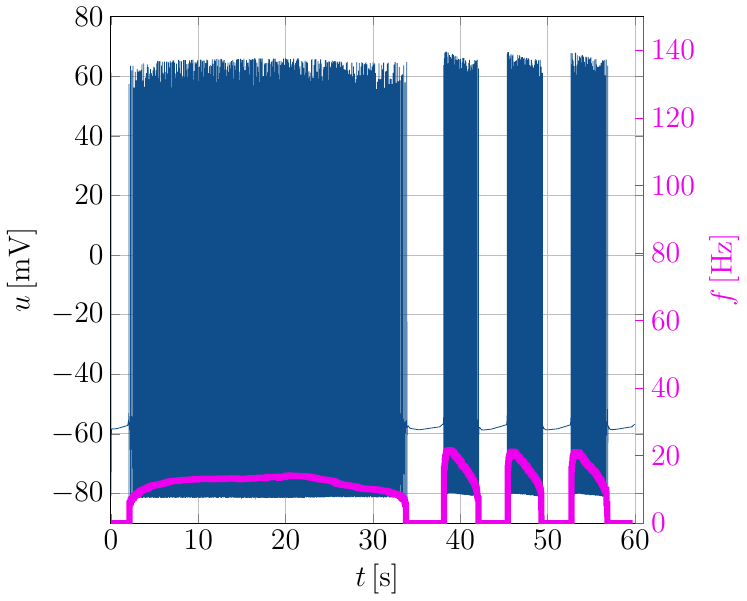}
\end{subfigure}
    \begin{subfigure}[b]{0.23\textwidth}
    \centering
\hspace{-0.8em}\includegraphics[width=1\textwidth]{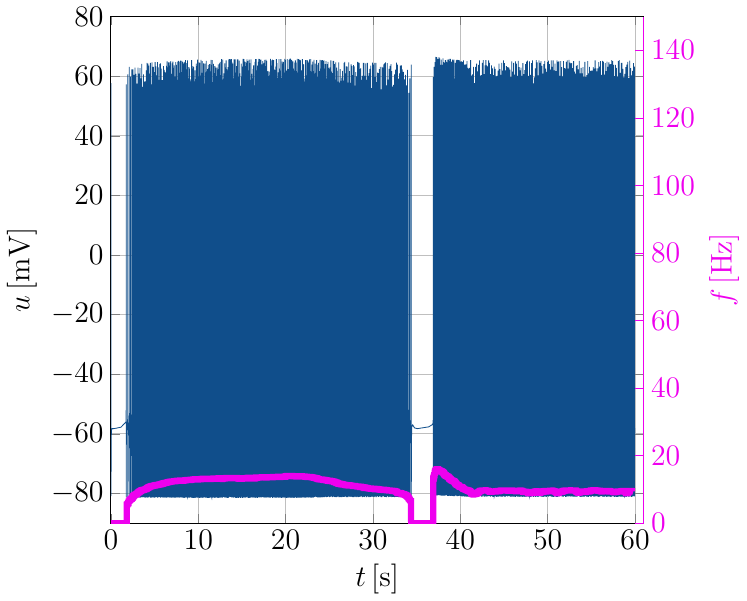}%
\end{subfigure}
    \begin{subfigure}[b]{0.23\textwidth}
    \centering
\hspace{-0.9em}\includegraphics[width=0.9\textwidth]{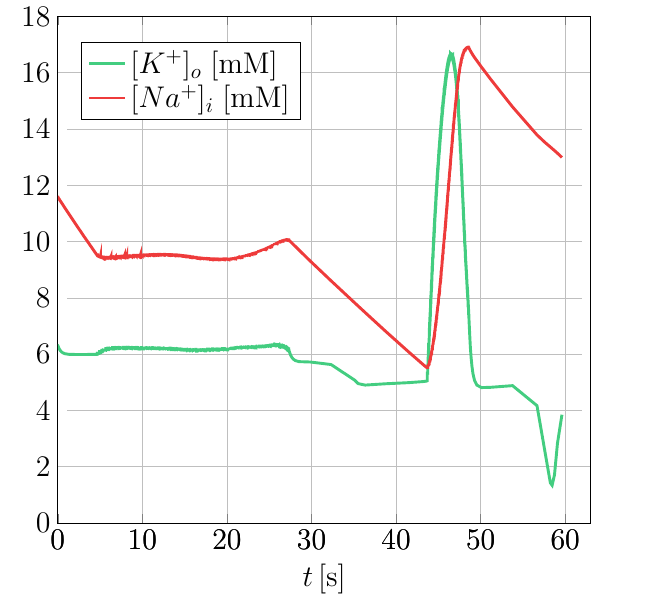}%
\caption{}\label{fig:eps2.7}
\end{subfigure}
\begin{subfigure}[b]{0.23\textwidth}
    \centering
\hspace{-0.7em}\includegraphics[width=0.9\textwidth]{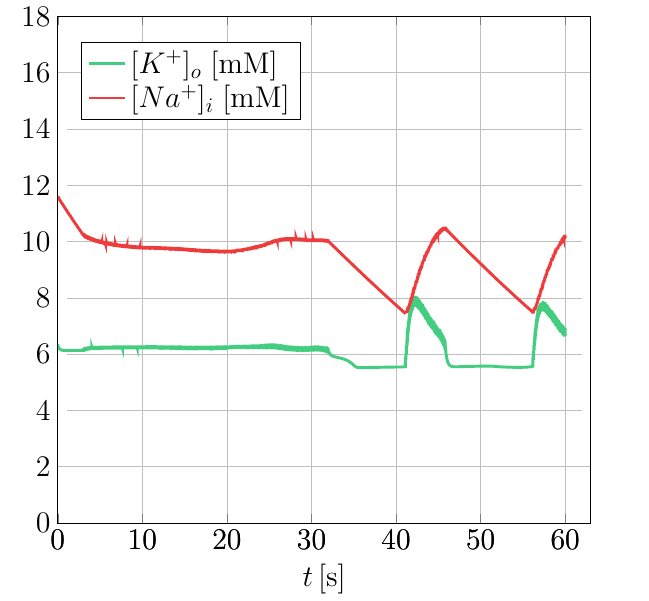}%
\caption{}\label{fig:eps5}
\end{subfigure}
    \begin{subfigure}[b]{0.23\textwidth}
    \centering
\hspace{-0.3em}\includegraphics[width=0.9\textwidth]{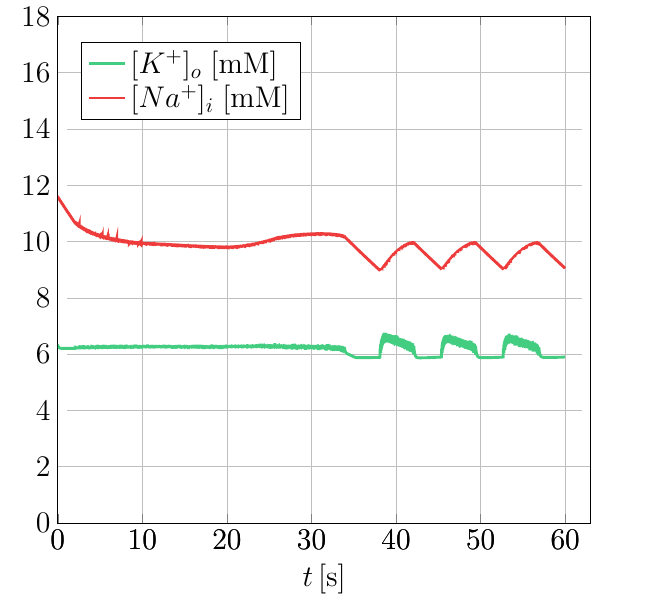}%
\caption{}\label{fig:eps9}
\end{subfigure}
    \begin{subfigure}[b]{0.23\textwidth}
    \centering
\hspace{4.9em}\includegraphics[width=0.9\textwidth]{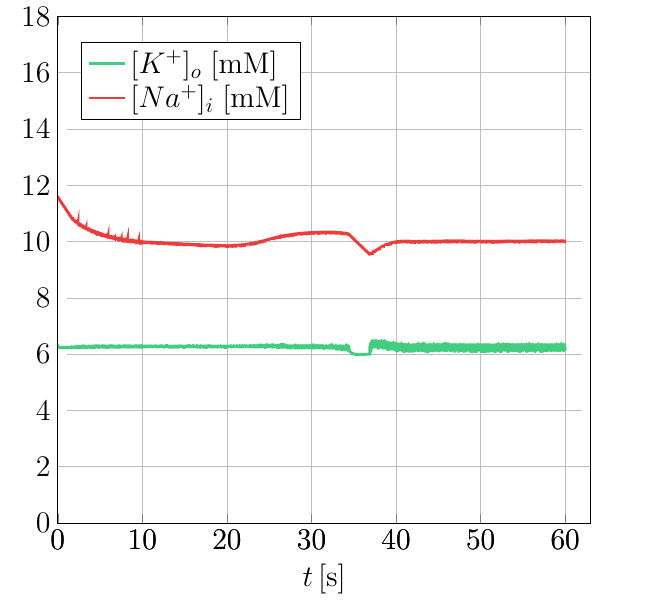}
\caption{}\label{fig:eps13}
    \end{subfigure}
\caption{Sensitivity analysis with respect to potassium clearance rate ($\varepsilon$) in severe ischemic condition. (a) $\varepsilon=2.7$.  (b)  $\varepsilon=5$. (c)  $\varepsilon=9.33$. (d) $\varepsilon=13$.}
    \label{fig:eps_comparison}
\end{figure}

Regarding the potassium clearance rate $\varepsilon$, higher values represent almost physiological conditions with slower reactivity to ionic perturbations. In contrast, lower values represent more pathological scenarios in which neuronal discharges rapidly accumulate $[\text{K}^+]_{\mathrm{o}}$ and trigger epileptic activity. 
A sensitivity analysis with respect to $\varepsilon$ is reported in Figure~\ref{fig:eps_comparison}, which shows both the transmembrane potential with the corresponding firing frequency (top row) and the ionic concentrations (bottom row). Lower values of $\varepsilon$ indicate a more pathological condition, as the frequency of reactive epileptic bursts increases and the ionic concentrations reach non-physiological regimes. 

For comparison across scenarios, simulations are run for a final time of $T=60\,\mathrm{s}$, which is sufficient to capture all relevant differences. 
In the highly pathological regime ($\varepsilon=2.7$, Figure \ref{fig:eps2.7}), the reduced diffusion leads to the redistribution of ions. As a consequence, neuronal discharges accumulate extracellular potassium, producing a runaway increase of $[\text{K}^+]_{{o}}$ above $16\,\mathrm{mM}$ and an abnormal depolarization of the transmembrane potential. The cell enters a hyper-excitable state, characterized by epileptic-like bursts with frequency exceeding $140\,\mathrm{Hz}$. 
Exploiting $\varepsilon=5$ and $\varepsilon=9.33$ as in Figures \ref{fig:eps5} and \ref{fig:eps9}, respectively, it is clear that the system reacts to ischemic stress with recurrent bursts of activity. The extracellular potassium concentration oscillates around $10$--$12\,\mathrm{mM}$, while the intracellular sodium concentration remains elevated but stable. This scenario captures a realistic pathological condition, where ischemia leads to periodic high-frequency discharges sustained by ionic dysregulation. In Figure \ref{fig:eps13}, the potassium diffusion is particularly high ($\varepsilon = 13$), and this leads to a rapid compensation of ionic perturbations. In this case, $[\text{K}^+]_{{o}}$ and $[\text{Na}^+]_{i}$ remain close to stable values, and the membrane potential shows only a brief transient. Although this scenario is not physiologically realistic, it demonstrates that high $\varepsilon$ prevents pathological bursting by dampening ionic feedback.
Overall, Figure \ref{fig:eps_comparison} highlights that the system dynamics strongly depend on the potassium diffusion rate $\varepsilon$: decreasing $\varepsilon$ from $13$ to $2.7$ shifts the behavior from stable ionic homeostasis to pathological regimes characterized by extracellular potassium accumulation above $16\,\mathrm{mM}$ and firing frequencies exceeding $140\,\mathrm{Hz}$, with intermediate values ($\varepsilon = 5$--$9.33$) producing periodic high-frequency bursting sustained by $\mathrm{[K^+]}_o$ oscillations around $10$--$12\,\mathrm{mM}$.

%% file: Num2D.tex
\subsection{Ischemic induced spontaneous spiking in idealized two-dimensional
domain}
\begin{figure}[h]
    \centering
    \begin{subfigure}[b]{0.5\textwidth}
    \centering
\includegraphics[width=\textwidth]{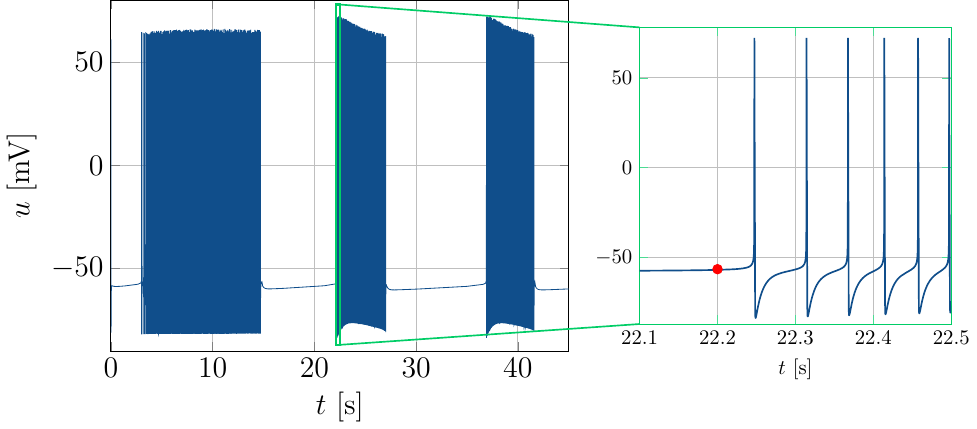}%
    \end{subfigure}
    \caption{Neuronal model simulation under severe ischemia (cellular swelling delay neglected). Zoom-in of the first pathological spikes of the transmembrane potential.}
    \label{fig:testcase1_0d_zoom}
\end{figure}

In this test case, we consider a simplified square geometry $\Omega = (0,1.5)\times(0,1.5)$ where the entire domain is characterized by grey matter tissue. The objective is to investigate the electrical activity of a small pathological ischemic subregion.
The domain is divided into two macro subregions: the pathological region, which represents grey matter tissue under pathological, severe ischemic conditions, and the surrounding region, where all variables and initial conditions are set to represent physiological tissue. We consider, in order to analyze the coupled behavior, the parameters $q_0 = 0.3q$, $\eta_{n} = 0.5409$, $\eta_{a} = 0.3938$, and $\eta_{ecs} = 0.0653$ for the ischemic subregion. The ischemic region is initialized in a transient resting-potential state, while a subregion, $\Omega_0$, models the tissue immediately prior to the onset of high-frequency activity, where initial conditions of the model are described in Table \ref{table:initial_conditions_testcase1} taken from \cite{leimer2024highorder, CALVETTI2018}.
Concerning the ischemic part of the domain, all values are taken once the decrease in blood flow has ended, in order to establish a region with blood flow deficit, right before the onset of pathological spikes, at time $T = 22.2\ \text{s}$, as shown in Figure \ref{fig:testcase1_0d_zoom}. At the selected time, all parameters and initial conditions are taken from the cellular neuronal simulation to characterize the two-dimensional ischemic tissue. 

The cellular response delay to volumetric changes is neglected; this simplification is considered acceptable, as the main focus here is on the propagation of self-induced spikes originating from a pathological ischemic region. 
Equation~\eqref{eq:conductivity_indicator_compact} represents the conductivity values associated with the grey matter under physiological and ischemic conditions:

\begin{equation}
\sigma_{n,l}(\mathbf{x}) =
2.735\,\mathds{1}_{\Omega_{\mathrm{GM}}}(\mathbf{x})
+ 0.565\,\mathds{1}_{\Omega_{\mathrm{IGM}}}(\mathbf{x}),
\label{eq:conductivity_indicator_compact}
\end{equation}

The quantities $\sigma_n$ and $\sigma_l$ denote the electrical conductivities in the normal and tangential directions, respectively, with respect to the local fiber orientation. 
In the numerical model, the computational domain $\Omega$ is partitioned into the physiological grey matter region, $\Omega_{\mathrm{GM}}$, and the ischemic grey matter region, $\Omega_{\mathrm{IGM}}$, and the corresponding conductivities are defined as piecewise constant functions over these subdomains.

\begin{figure}[h]
\centering
\includegraphics[width=\textwidth]{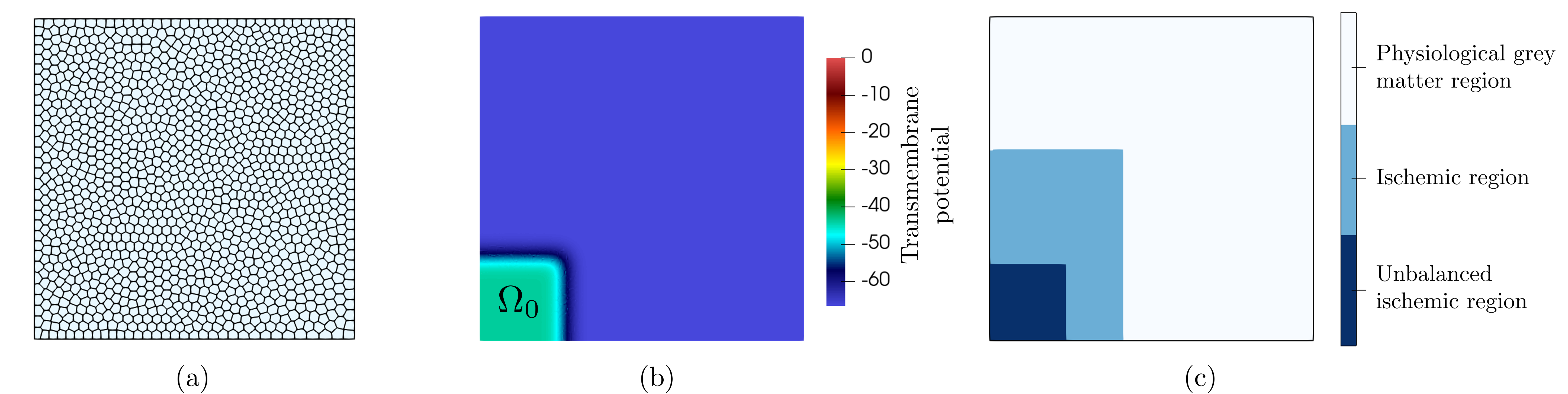}  
\caption{Initial setup: (a) Computational polytopal mesh with 1800 elements. (b) Initial condition of the transmembrane potential; $\Omega_0$: tissue immediately prior to the onset of high-frequency activity. (c) Computational ischemic region in grey matter tissue.}
\label{fig:setting_testcase1}
\end{figure}

\begin{figure}[H]
\centering
  \begin{subfigure}{\textwidth}
  \centering
  \includegraphics[width=0.4\textwidth]{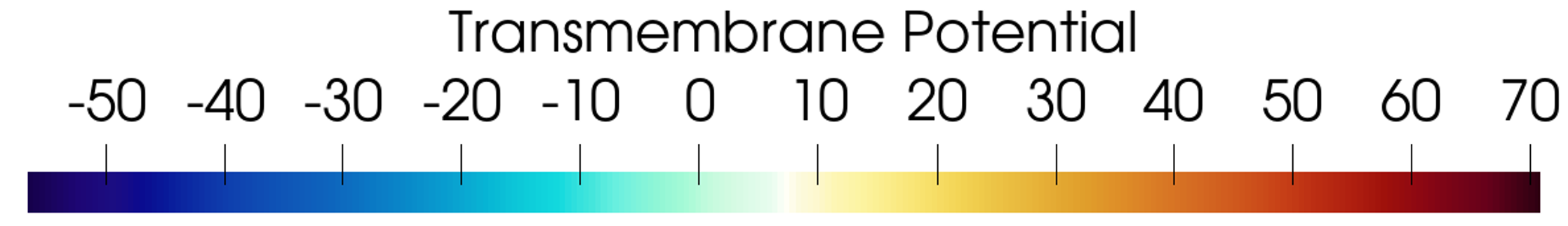}
  \end{subfigure}
  \begin{subfigure}[t]{0.24\textwidth}
    \centering
    \includegraphics[width=\textwidth]{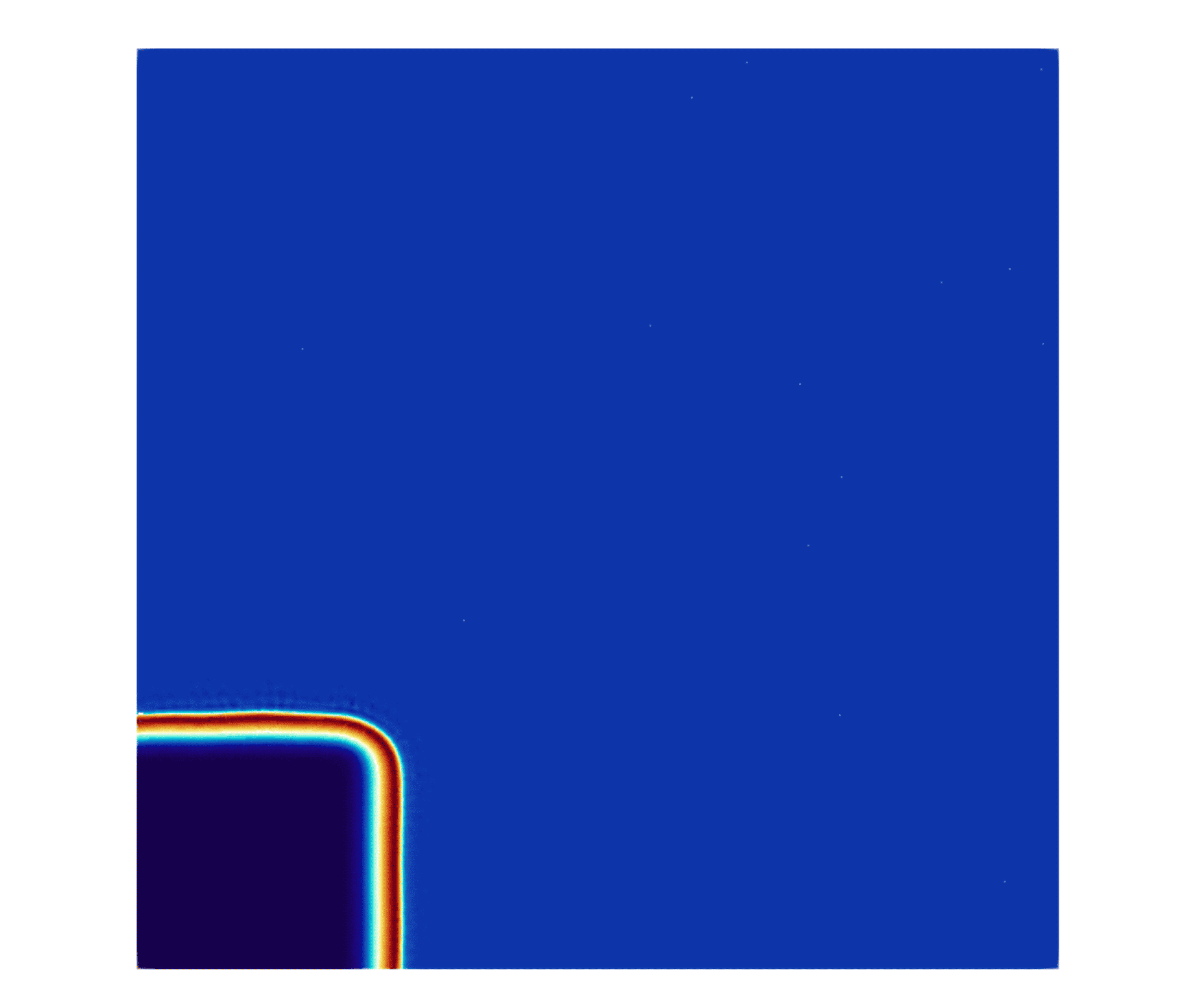}
\caption{$t = 1.8\ \mathrm{ms}$}
  \end{subfigure}%
  \begin{subfigure}[t]{0.24\textwidth}
    \centering
    \includegraphics[width=\textwidth]{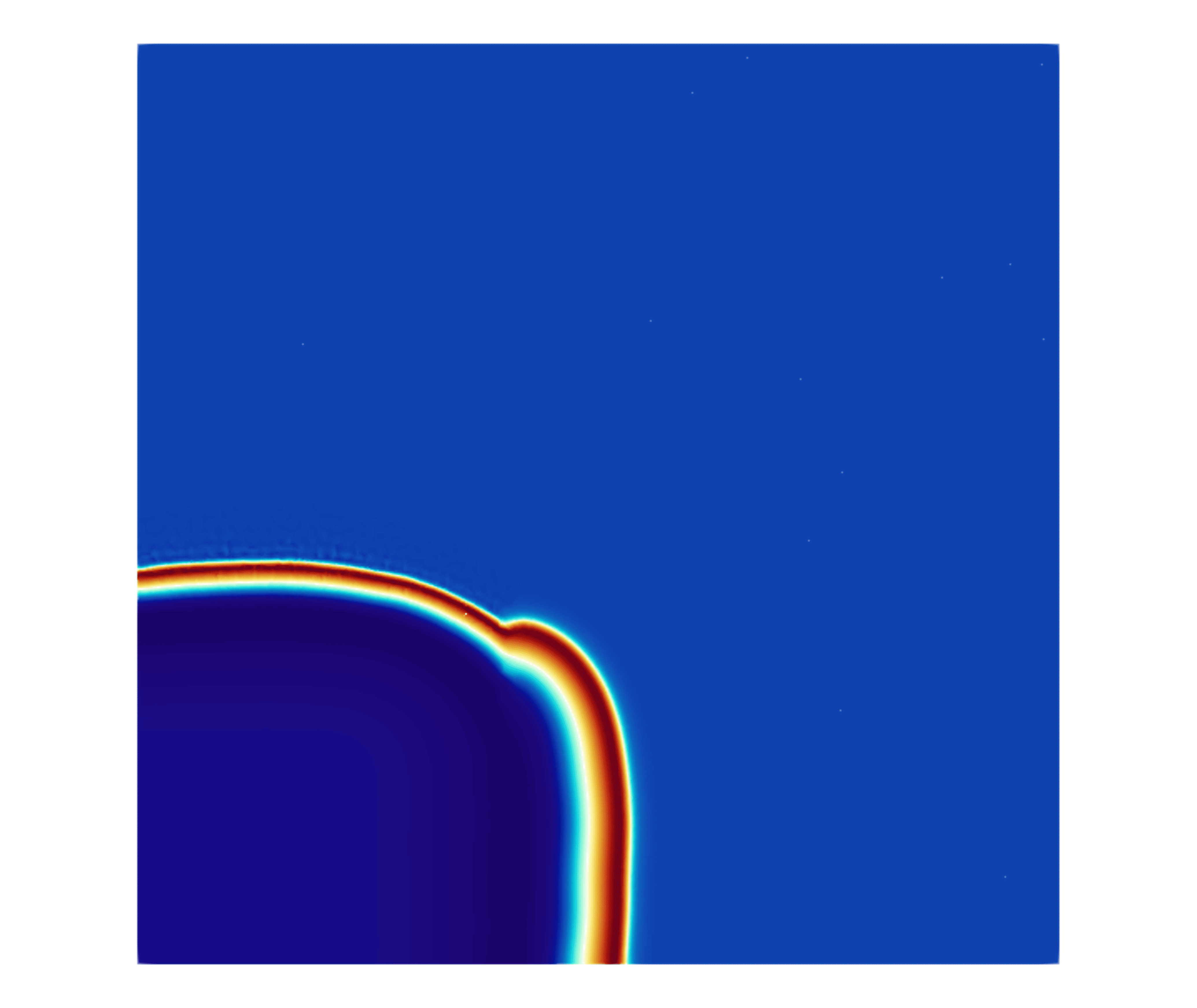}
\caption{$t = 3.15\ \mathrm{ms}$}
  \end{subfigure}
    \begin{subfigure}[t]{0.24\textwidth}
    \centering
    \includegraphics[width=\textwidth]{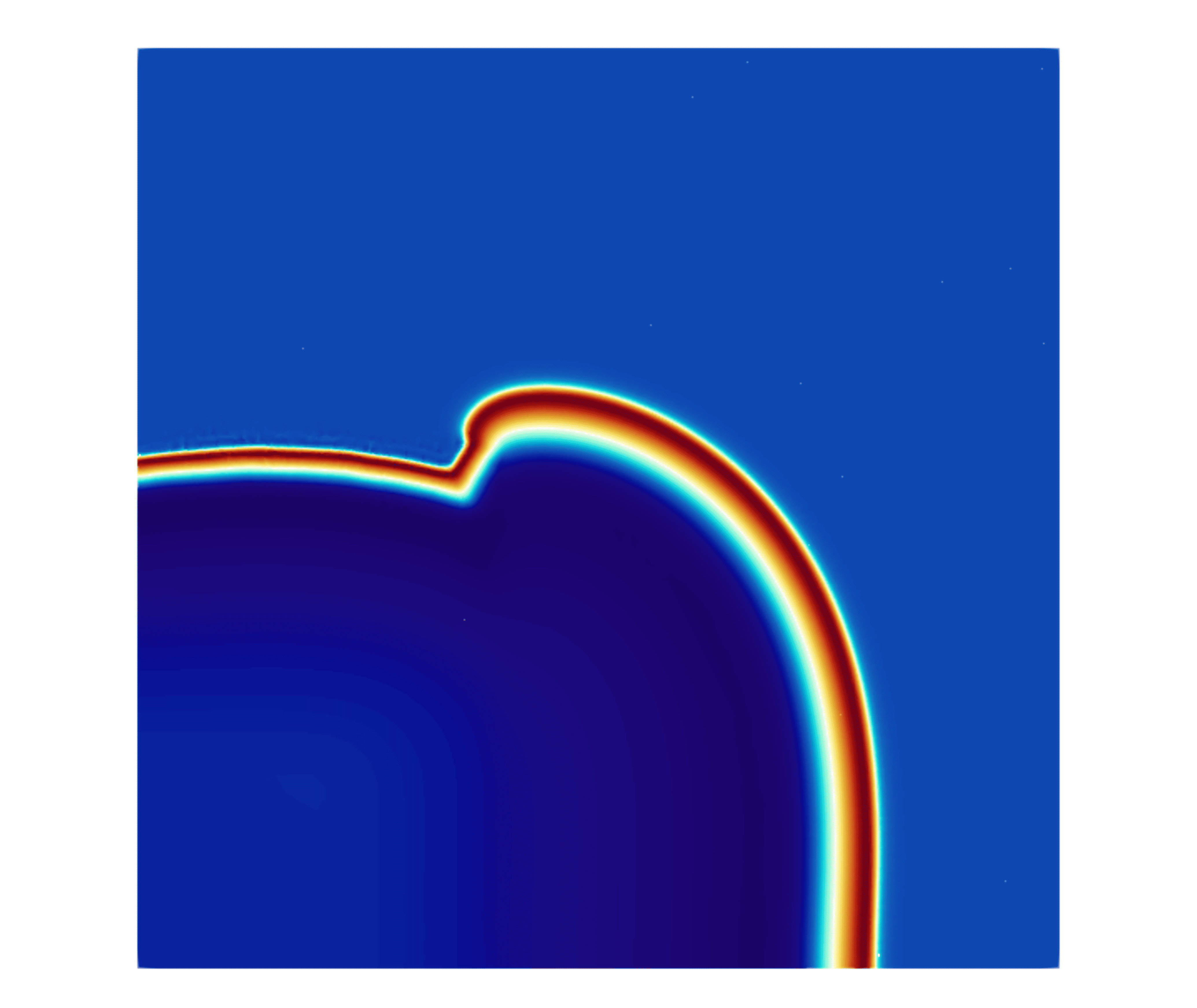}
\caption{$t = 5.25\ \mathrm{ms}$}
  \end{subfigure}
      \begin{subfigure}[t]{0.24\textwidth}
    \centering
    \setcounter{subfigure}{6}
    \includegraphics[width=1.1\textwidth]{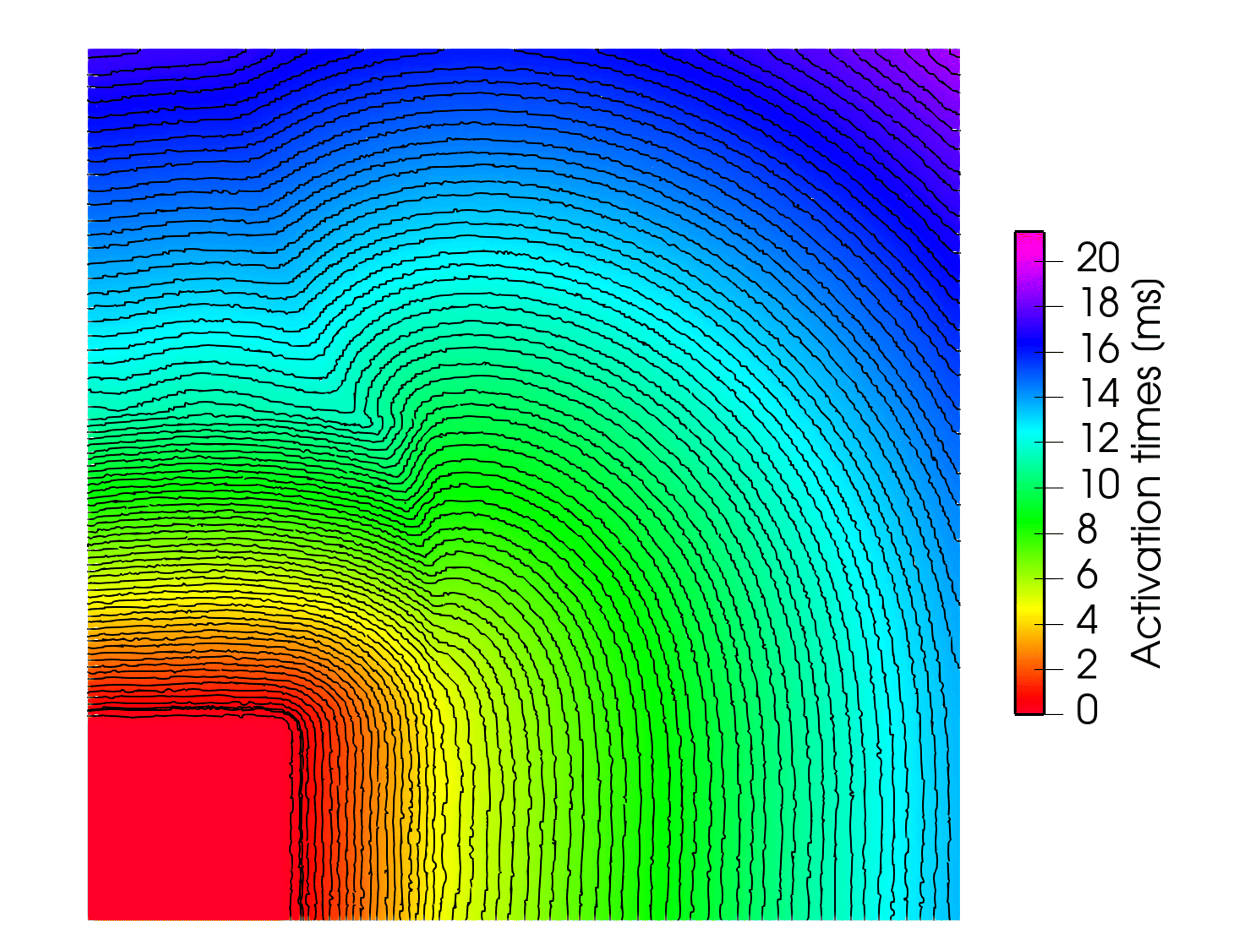}
\caption{Activation times for the first wavefront $[\mathrm{ms}]$}
  \end{subfigure}
\hfill
  \begin{subfigure}[t]{0.24\textwidth}
    \centering
    \setcounter{subfigure}{3}
    \includegraphics[width=\textwidth]{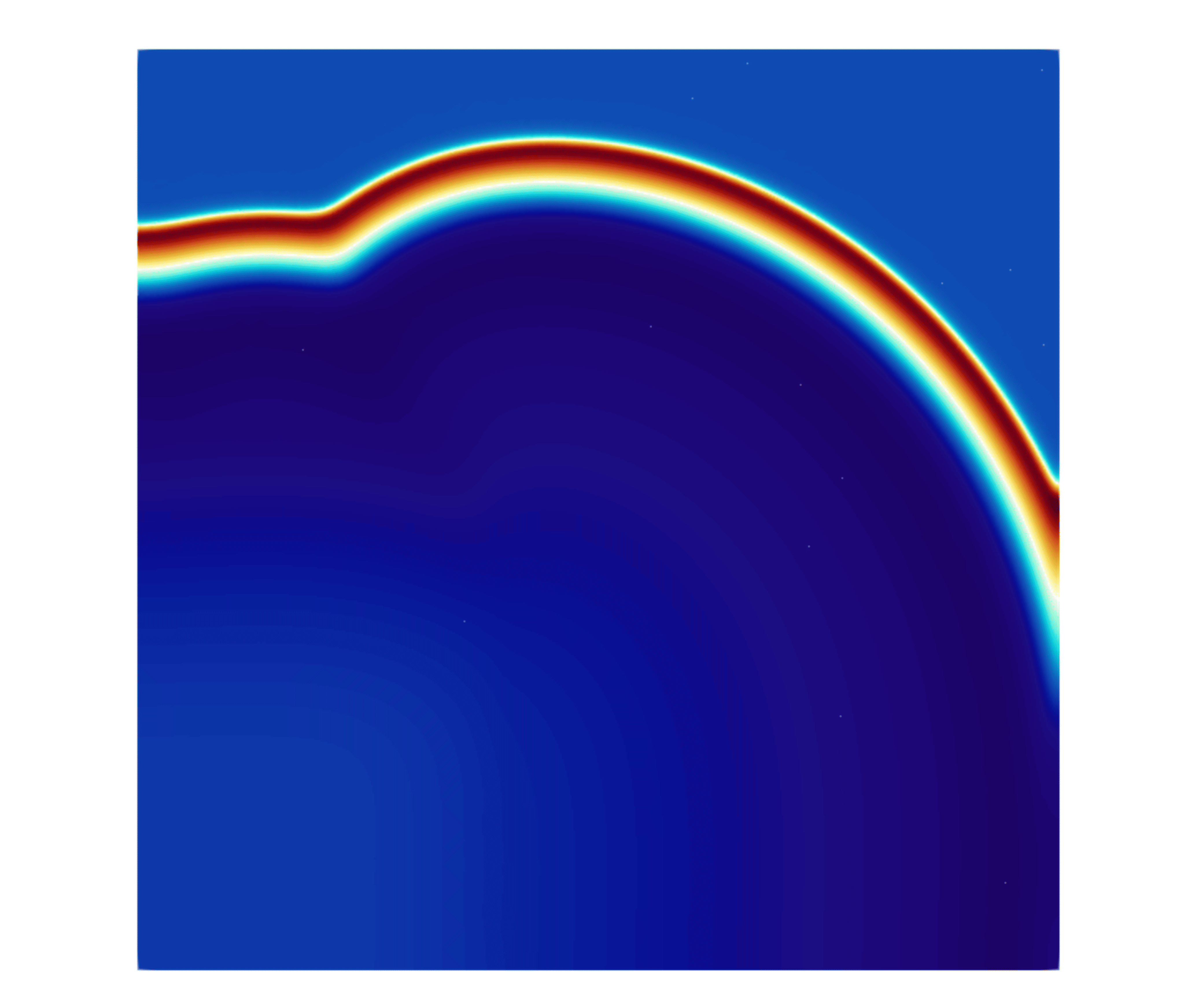}
\caption{$t = 8.85\ \mathrm{ms}$}
  \end{subfigure}%
  \begin{subfigure}[t]{0.24\textwidth}
    \centering
  \includegraphics[width=\textwidth]{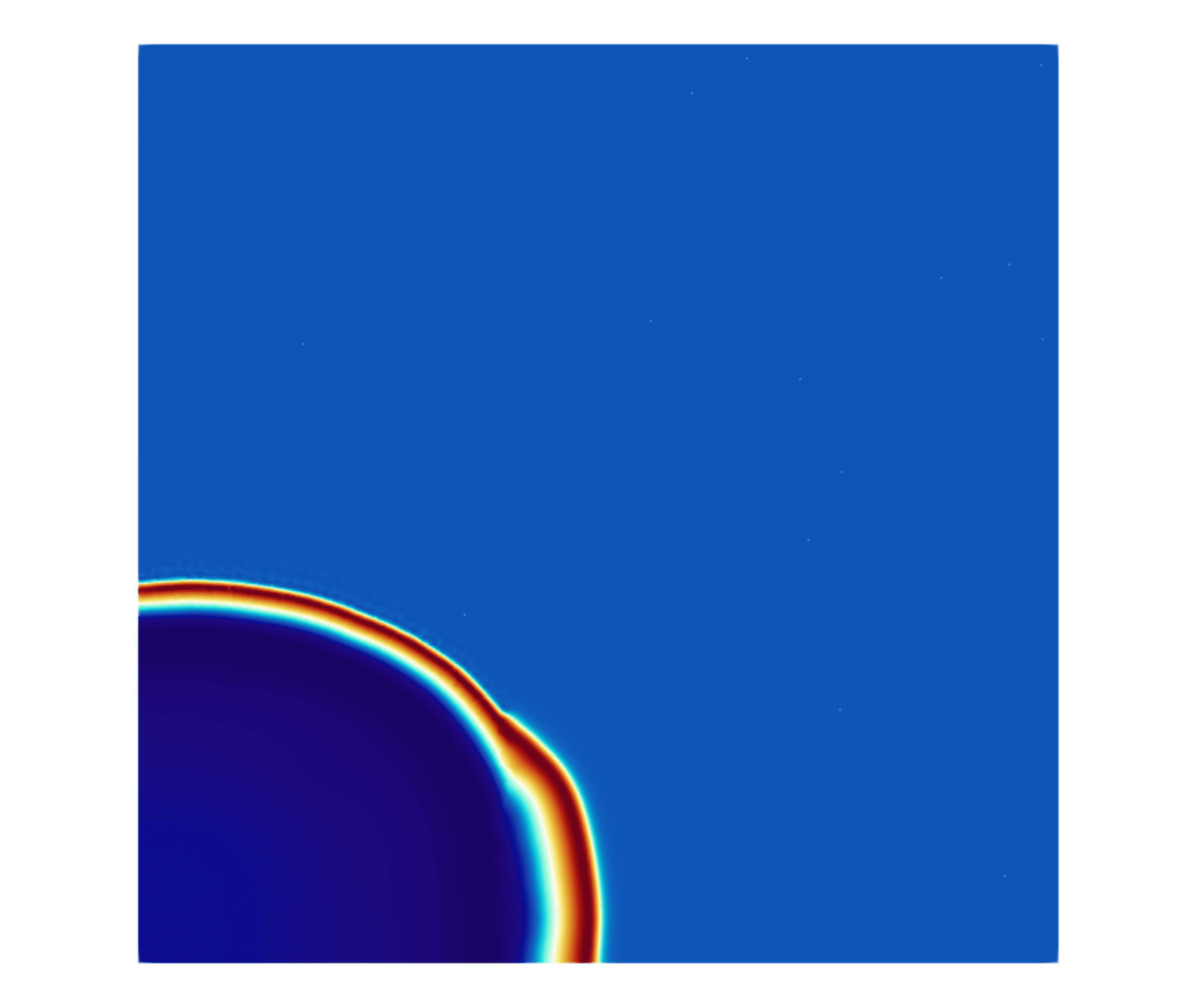}
\caption{$t = 42.75\ \mathrm{ms}$}
  \end{subfigure}
    \begin{subfigure}[t]{0.244\textwidth}
    \centering
    \includegraphics[width=\textwidth]{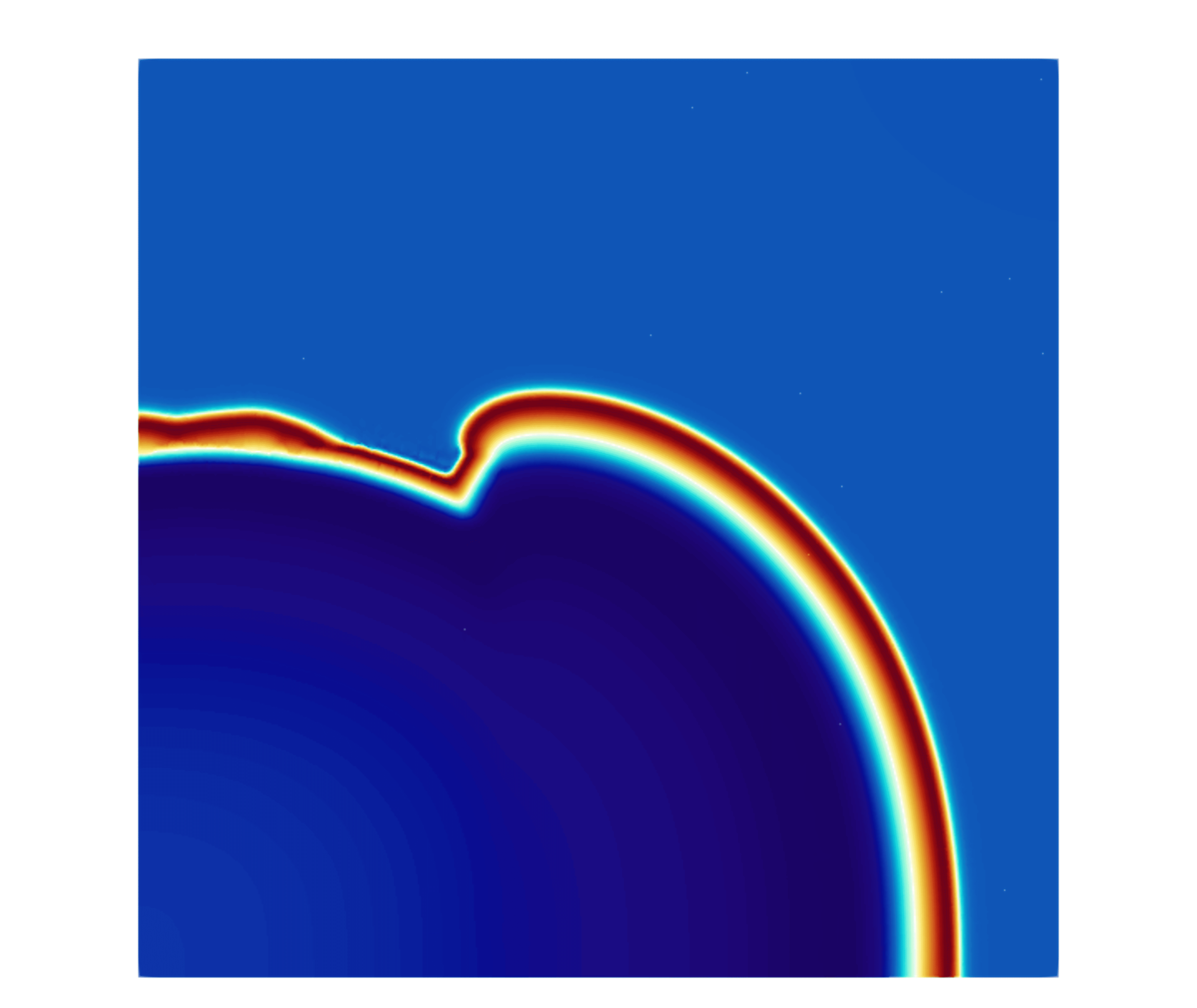}
\caption{$t = 48.25\ \mathrm{ms}$}
  \end{subfigure}
    \begin{subfigure}[t]{0.24\textwidth}
    \centering
    \includegraphics[width=1.1\textwidth]{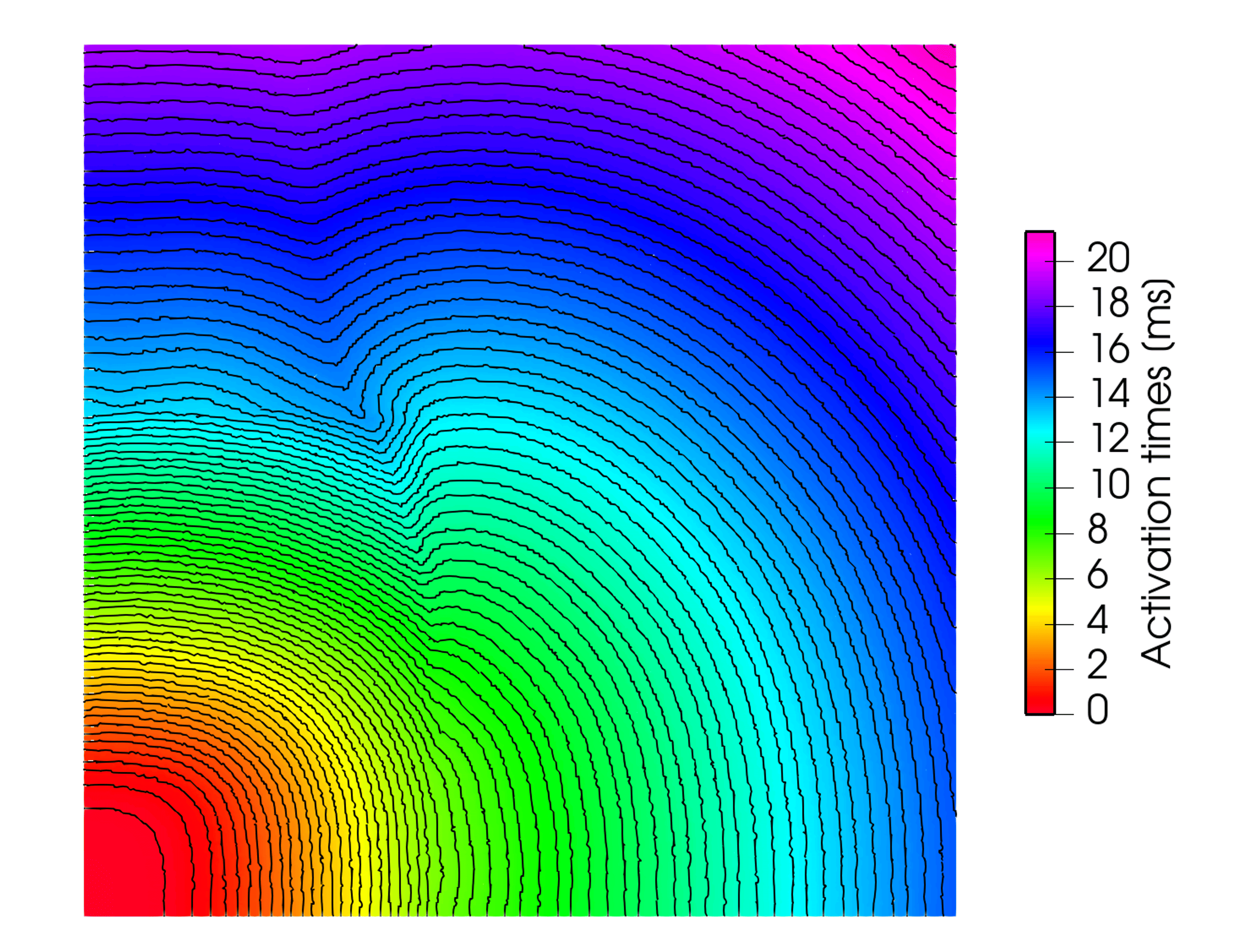}
\caption*{ (h)\, Activation times for the second wavefront $[\mathrm{ms}]$}
  \end{subfigure}
\caption{(a)-(f) Snapshots of the simulation at different time steps, illustrating the evolution of the transmembrane potential in a squared tissue of grey matter with an unstable pathological ischemic condition for first and second auto-induced wavefront. (g)-(h) Activation times for the first and second travelling wavefront in the domain.}
\label{fig:test_case1_potential}
\end{figure}

\begin{table}[h]
\centering
\begin{tabular}{|lcc|lcc|lcc|}
\hline
\text{Variable} & \text{Value} & \text{Units} &
\text{Variable} & \text{Value} &
\text{Units} &
\text{Variable} & \text{Value} & 
\text{Units} \\
\hline 
$\mathrm{u}^0$ & $-46.735$ & [mV] & $[O_{2}]^0_\text{ECS}$ & $0.269$ &[mM] &  $\text{ATP}^0_a$ & $2.17$ & [mM]  \\
$\mathrm{[Ca]}_i^0$ & $1.168 \cdot 10^{-3}$ & [mM] & $[O_{2}]^0_{a}$ & $0.0166$ & [mM] & $\text{ATP}^0_n$ & $2.18$ & [mM] \\
$\mathrm{[K]}_o^0$ & $5.6336$ & [mM] & $[O_{2}]^0_{n}$ & $0.0149$  & [mM] &  $\text{ADP}^0_a$ & $0.03$ & [mM] \\
$\mathrm{[Na]}_i^0$ & $8.2451$  & [mM] & $[O_{2}]^0_{b}$ & $3.0654$  & [mM] & $\text{NADH}^0_n$ & $1.2\cdot 10^{-3}$ & [mM] \\
$m^0$ & $0.1189$ && $\text{ADP}^0_n$ & $6.3\cdot 10^{-3}$  & [mM] &  $\text{NAD}^0_n$ & $0.03$ & [mM] \\
$n^0$ & $0.2067$ && $\text{NAD}^0_a$ & $0.03$ & [mM]& $\text{NADH}^0_a$ & $1.2\cdot 10^{-3}$ & [mM] \\
$h^0$ & $0.8338$ && &&& &&\\
\hline
\end{tabular}
\caption{Initial conditions for variables of the ionic coupled model for a simulation of severe ischemia.}
\label{table:initial_conditions_testcase1}
\end{table}
In Figure \ref{fig:setting_testcase1}, we represent the polytopal mesh exploited for the simulation, with 1800 elements ($h = 0.084$), the initial condition of the transmembrane potential, and the definition of the ischemic region in the domain. The simulation is performed exploiting a $\Delta t =2.5\cdot 10^{-3} \ \text{ms} $, and final time $T = 50\ \text{ms}$. The polynomial degree is automatically adapted in the range between 1 and 5, exploiting the $p$-adaptive algorithm \cite{leimer2025p}. The evolution of the transmembrane potential over the domain at selected snapshots is shown in Figure \ref{fig:test_case1_potential}. The main goal of this simulation is to understand how electrical activity originating from a pathological ischemic region affects the surrounding tissue.
The results represented in Figure \ref{fig:test_case1_potential} show that high-frequency spikes are generated within the ischemic area and then propagate into the healthy tissue, suggesting that even a localized pathological region can lead to large-scale effects on the whole system. In particular, we note that the activation originates from the pathological ischemic zone, where the conductivity of the model is reduced. This lower conductivity slows down the propagation of the electrical wave, where the wavefront advances more slowly within the ischemic region before invading the surrounding healthy tissue.
We observe that this evolution is consistent with the underlying biochemical mechanisms; in fact, ischemia is known to induce significant structural and metabolic alterations in neural tissue. In particular, cellular swelling, lack of oxygen, reduced ATP production, and loss of ionic homeostasis all contribute to reducing the ECS volume and increasing its tortuosity.  This results in a lower effective electrical conductivity of the tissue, which is expected to slow down the propagation of depolarization waves. Because the local currents are weaker, the propagating wavefront becomes not only slower but also thinner \cite{ShawRudy1997_slow_conduc}.

\subsection{Activation of a pathological ischemic subregion in idealized two-dimensional domain}
\label{sec:external}
\begin{figure}[h]
\centering
\includegraphics[width=\textwidth]{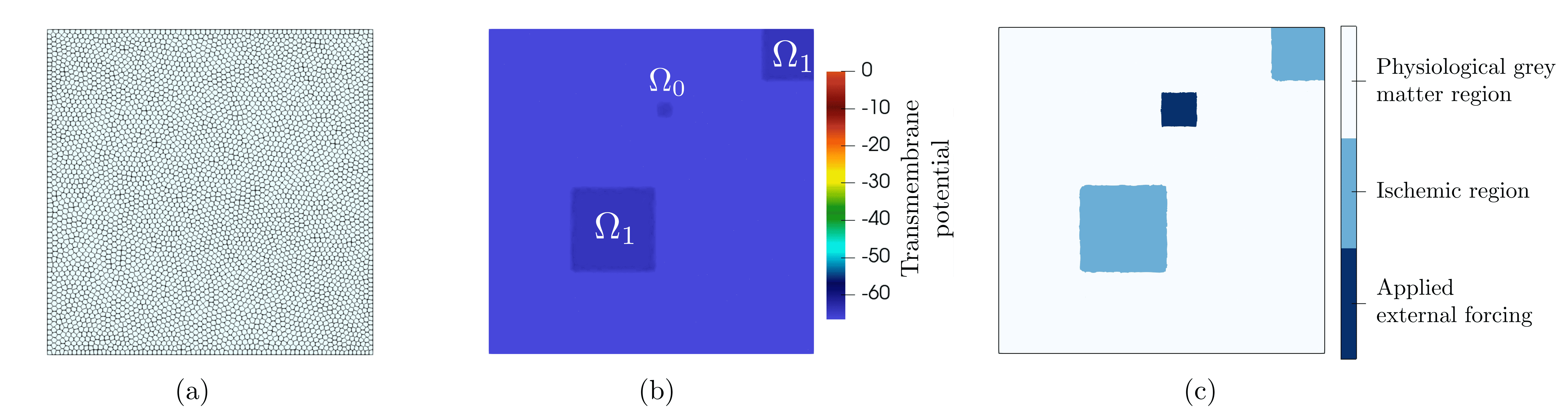}  
\caption{Initial setup: (a) Computational polytopal mesh with 2200 elements. (b) Initial condition of the transmembrane potential; $\Omega_1$: ischemic grey matter subregions; $\Omega_0$: subregion triggered by external stimulus. (c) Computational ischemic regions in grey matter tissue.}
\label{fig:setting_testcase2}
\end{figure}

In this second test case, instead of focusing on the propagation of self-induced spikes from an ischemic subregion,  we aim to investigate the ischemic tissue behavior when it is stimulated by external impulses originating from surrounding tissue. The computational geometric setting is illustrated in Figure \ref{fig:setting_testcase2}.

To avoid self-induced high-frequency activation of the ischemic region, as shown in the previous test case, we set the initial condition for the pathological subregion to a transient resting value for the transmembrane potential.
An additional external forcing is applied in a limited region of the physiological grey matter, in order to generate several spikes that subsequently propagate in the surrounding heterogeneous tissue. 
To analyze the coupled behavior, we consider the parameters $q_0 = 0.3q$, $\eta_{n} = 0.5409$, $\eta_{a} = 0.3938$, and $\eta_{ecs} = 0.0653$ for the ischemic subregions.   
In Figure \ref{fig:setting_testcase2} we represent the polytopal mesh exploited for the simulation, where $\Omega = (0, 1.7)^2$ with 2800 elements ($h = 0.064$), the initial condition of the transmembrane potential and the definition of the ischemic regions in the grey matter tissue ($\Omega_1$) and the part of the domain where the external forcing is applied ($\Omega_\text{0}$). In the part of the domain $\Omega \setminus ( \Omega_0 \cup \Omega_1)$ we set physiological grey matter tissue; the external term applied in $\Omega_0$ is $I_{\text{Act}}=40\ \mu \mathrm{A}$, in order to induce a series of sustained action potentials in the brain tissue.
The initial conditions for all the variables of the ionic model are taken from the neuronal simulation represented in Figure \ref{fig:testcase2_0d_zoom}, at time $ \hat{t} = 20 \; \mathrm{s}$. 
The numerical simulation is performed exploiting $\Delta t = 2.5\cdot 10^{-3}\ \text{ms}$ and we considered different conductivity values for the tissue: in phyisiological grey matter the conductivity is set to $\sigma_n = \sigma_l = 1.735\ \text{Sm}^{-1}$, while in the ischemic grey matter it is reduced to $\sigma_n = \sigma_l = 0.358\ \text{Sm}^{-1}$. 

\begin{figure}[H]
\centering
  \begin{subfigure}[b]{\textwidth}
    \centering
  \includegraphics[width=0.4\textwidth]{Photos/2D/scale.png}
  \end{subfigure}
  \begin{subfigure}[b]{0.20\textwidth}
    \centering
    \includegraphics[width=\textwidth]{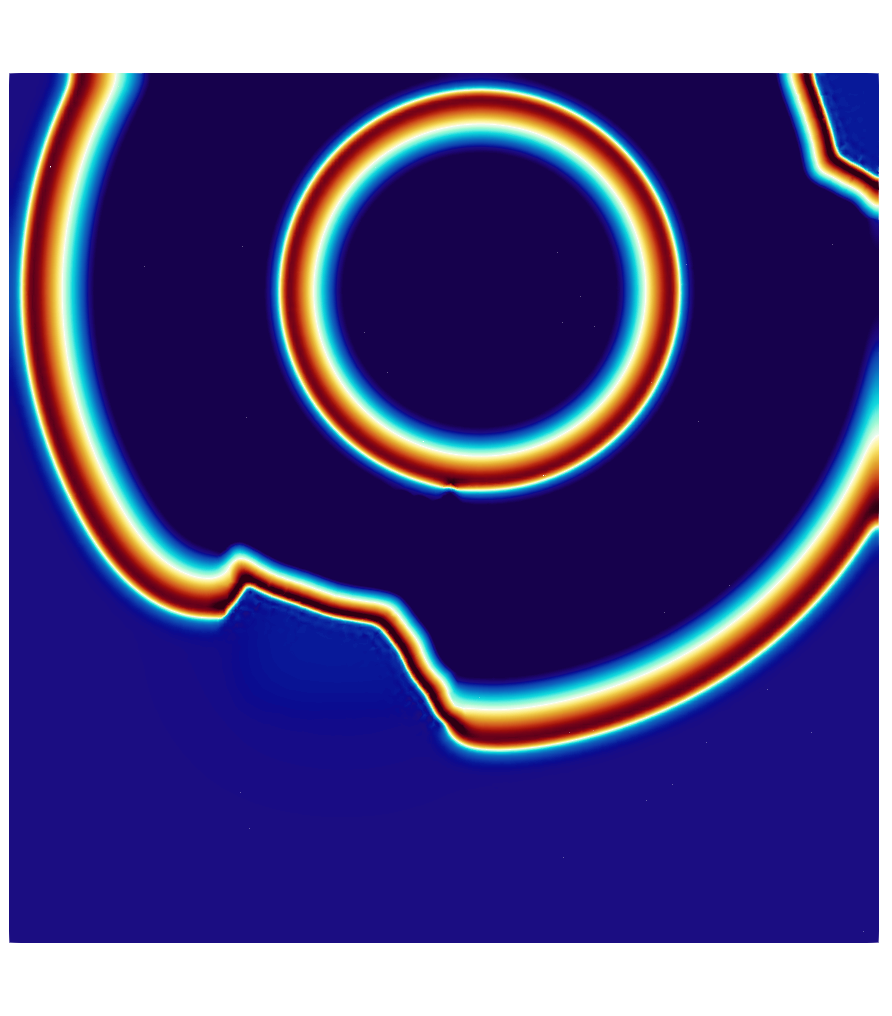}
    \caption{$t = 6.30\ \mathrm{ms}$}
  \end{subfigure}%
  \hfill
  \begin{subfigure}[b]{0.20\textwidth}
    \centering
    \includegraphics[width=\textwidth]{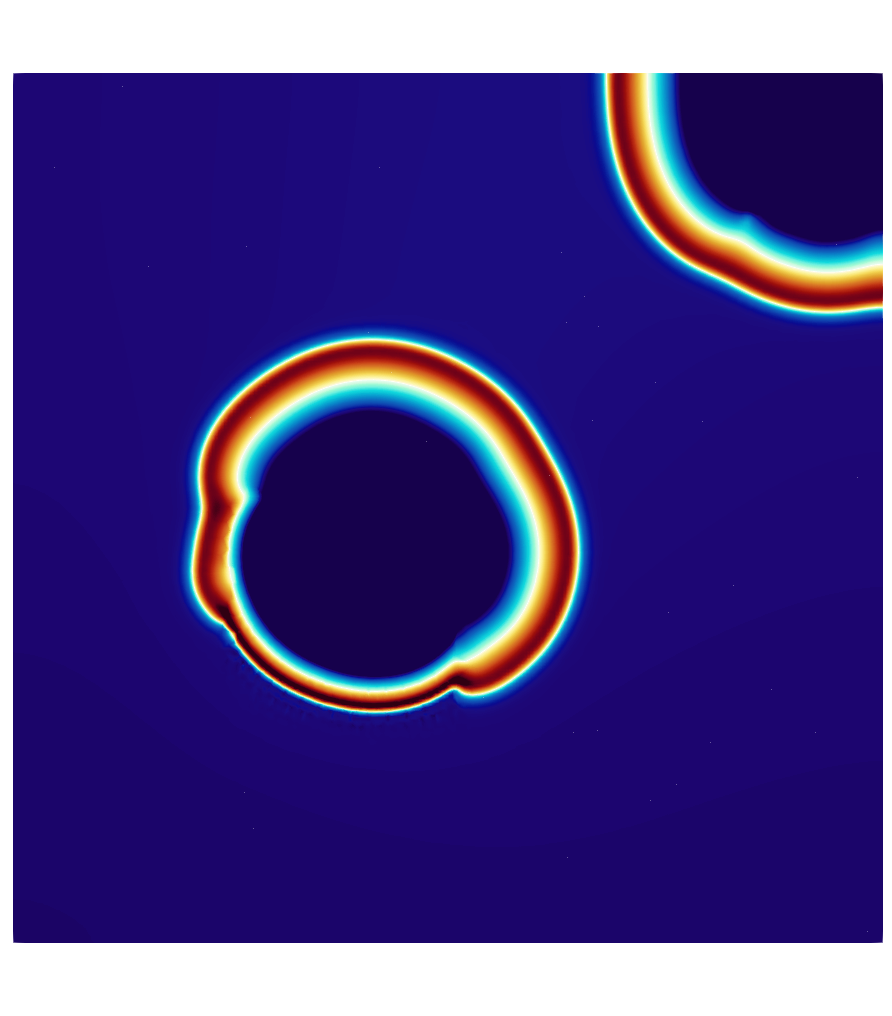}
    \caption{$t = 52.4\ \mathrm{ms}$}
  \end{subfigure}
  \hfill
  \begin{subfigure}[b]{0.20\textwidth}
    \centering
    \includegraphics[width=\textwidth]{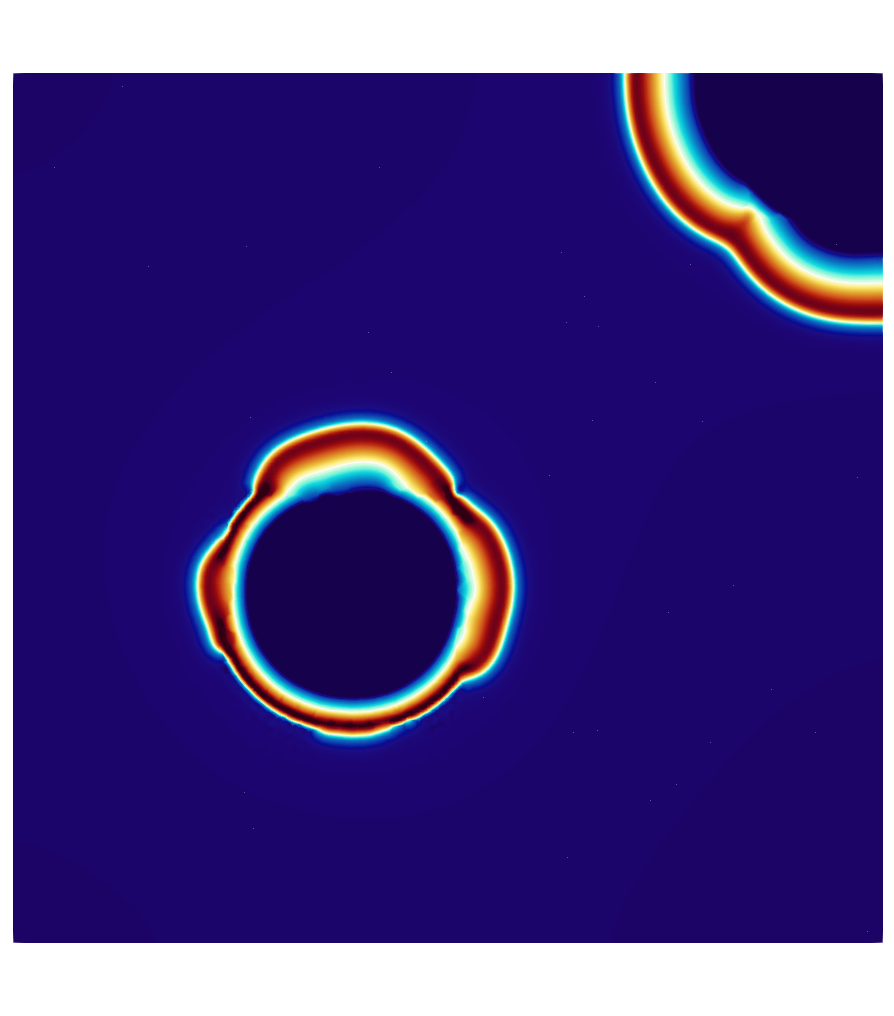}
    \caption{$t = 109.3\ \mathrm{ms}$}
  \end{subfigure}%
  \hfill
  \begin{subfigure}[b]{0.20\textwidth}
    \centering
    \includegraphics[width=\textwidth]{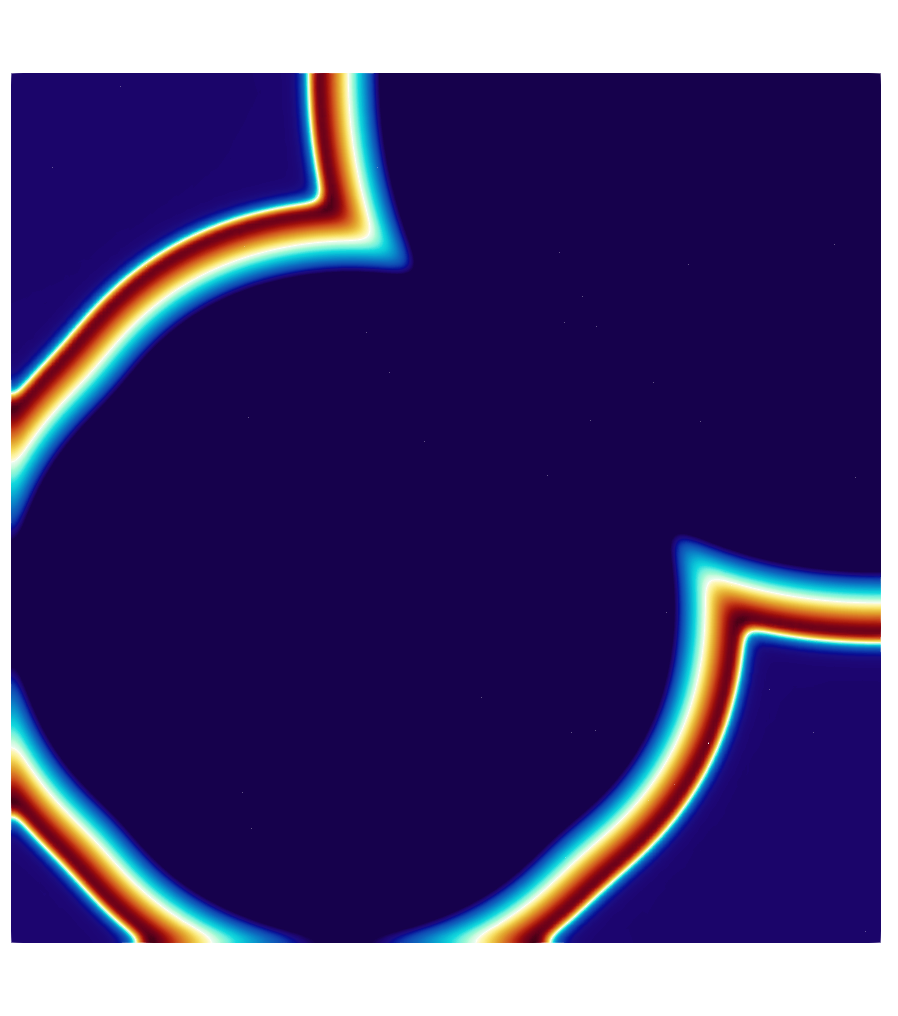}
    \caption{$t = 162.6\ \mathrm{ms}$}
  \end{subfigure}\vspace{.3em}
   \begin{subfigure}[b]{\textwidth}
    \centering
  \includegraphics[width=0.5\textwidth]{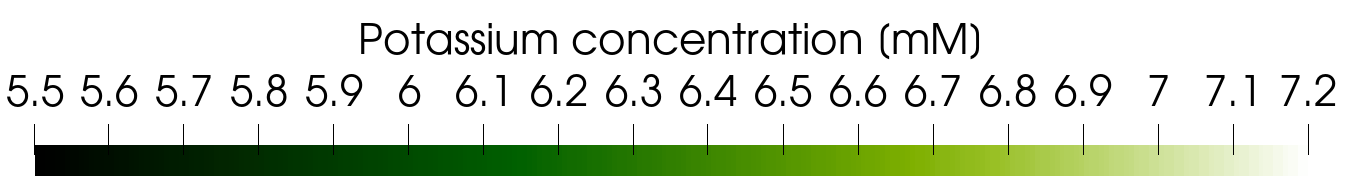}
  \end{subfigure}
  \begin{subfigure}[b]{0.20\textwidth}
    \centering
    \includegraphics[width=0.985\textwidth]{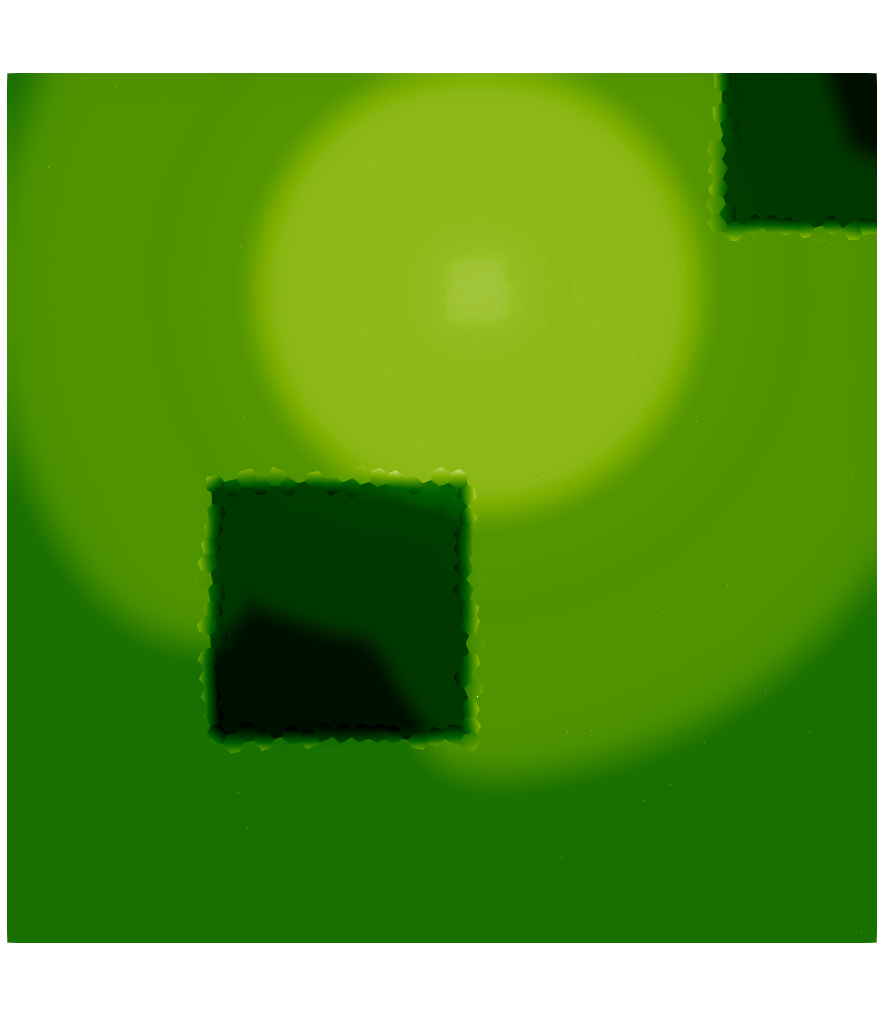}
    \caption{$[\text{K}]_\text{o} \text{ at }  t= 6.30\ \mathrm{ms}$}
  \end{subfigure}%
  \hfill
  \begin{subfigure}[b]{0.20\textwidth}
    \centering
    \includegraphics[width=\textwidth]{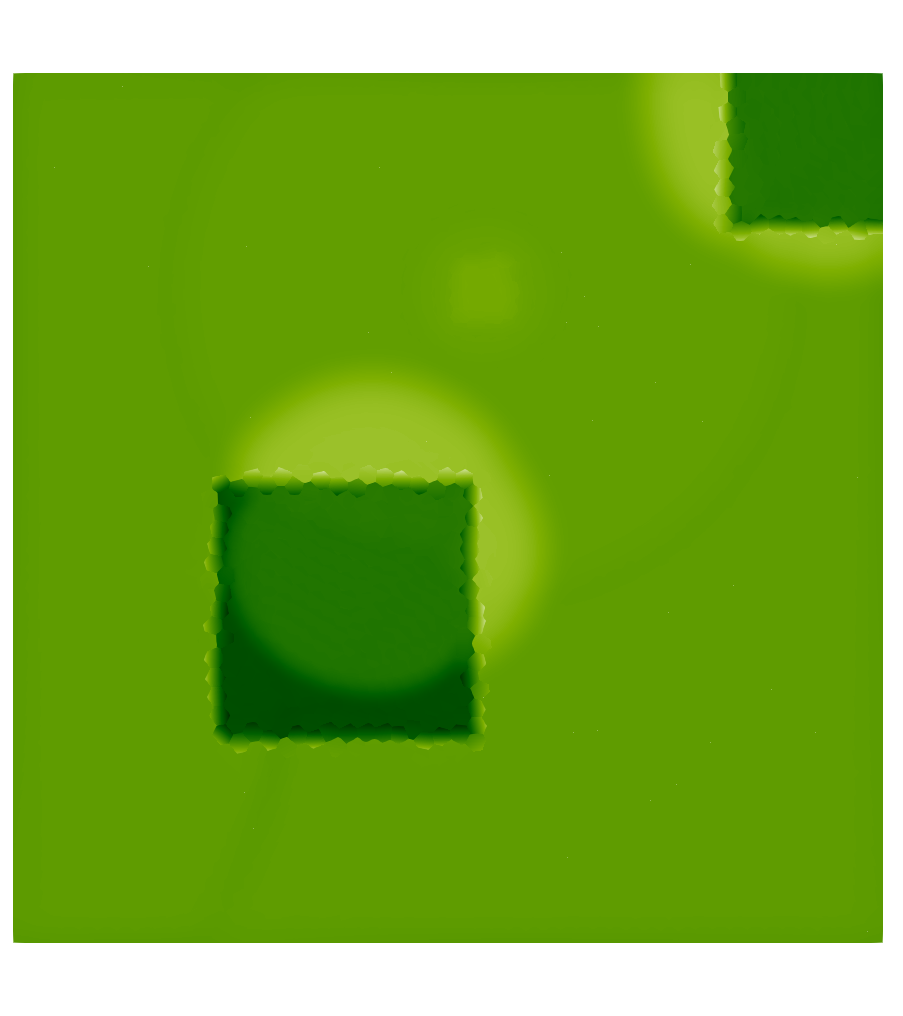}
    \caption{$[\text{K}]_\text{o} \text{ at }  t= 52.4\ \mathrm{ms}$}
  \end{subfigure}
  \hfill
  \begin{subfigure}[b]{0.20\textwidth}
    \centering
    \includegraphics[width=\textwidth]{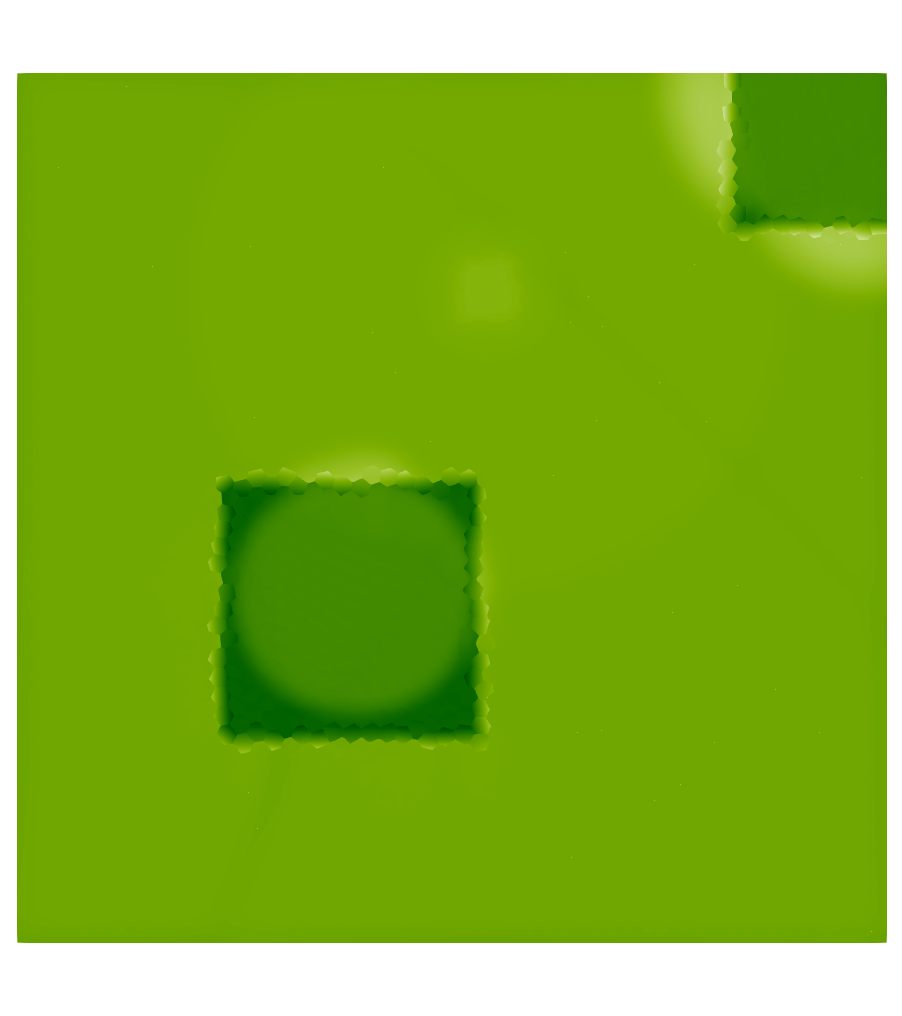}
    \caption{$[\text{K}]_\text{o} \text{ at }  t= 109.3\ \mathrm{ms}$}
  \end{subfigure}%
  \hfill
  \begin{subfigure}[b]{0.20\textwidth}
    \centering
    \includegraphics[width=\textwidth]{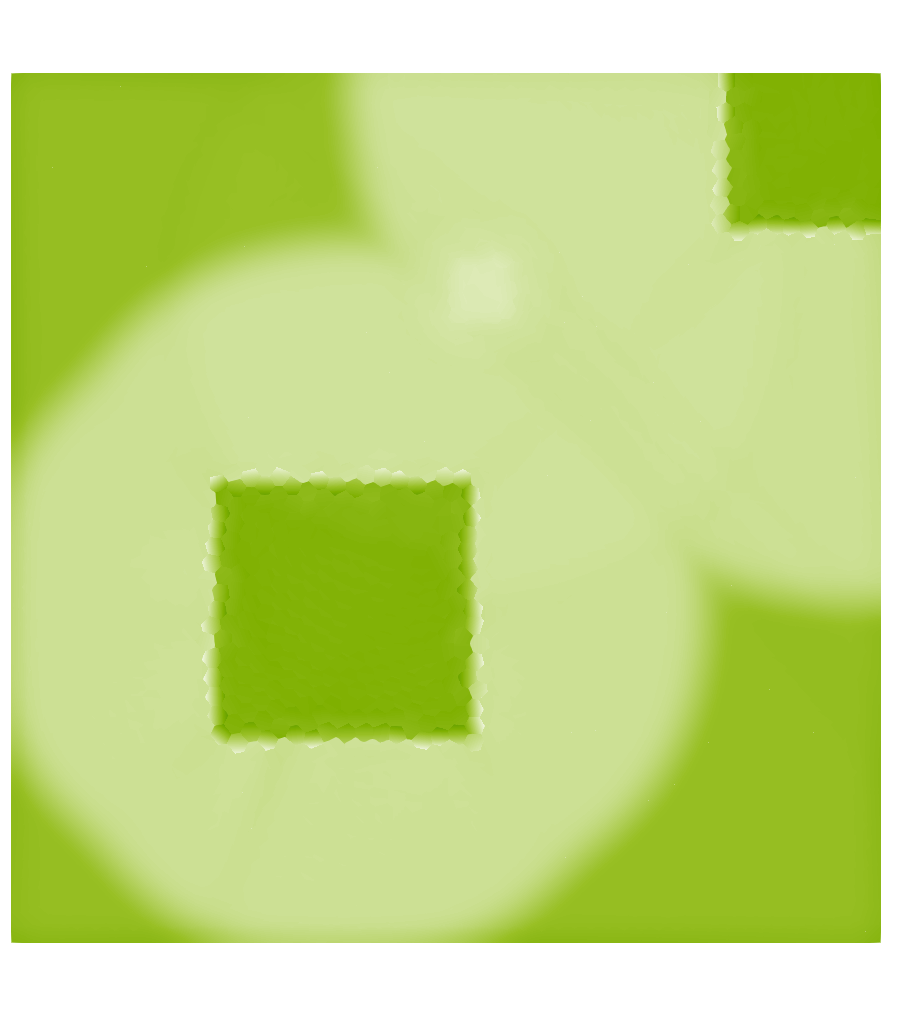}
    \caption{$[\text{K}]_\text{o} \text{ at }  t= 162.6\ \mathrm{ms}$}
  \end{subfigure}
\caption{(a)-(d) Different snapshots of the transmembrane potential evolution in grey matter tissue with externally activated part of physiological brain tissue. (e)-(h) Different snapshots of the extracellular potassium concentration in physiological and ischemic grey matter tissue.}
\label{fig:test_case2_potential}
\end{figure}

Figure~\ref{fig:test_case2_potential} illustrates the spatio-temporal evolution of the electrophysiological response in grey matter tissue containing two ischemic subregions, by reporting snapshots of the transmembrane potential and the extracellular potassium concentration \([\text{K}^+]_\text{o}\) at different time snapshots.

Figure~\ref{fig:test_case2_potential}, in panels (a)--(d), shows that the external forcing applied to the healthy gray matter initially generates multiple activation wavefronts, which propagate through the surrounding tissue. As the wavefronts expand across the computational domain, they reach the ischemic subregions, where the altered electrophysiological properties facilitate their activation into the ongoing dynamics. After being activated, the ischemic areas can sustain pathological electrical activity, effectively acting as secondary sources that further shape the spatio-temporal propagation pattern. Overall, this sequence highlights how a spatially localized stimulus can evolve into a self-sustained and more complex activity due to interactions between wavefronts and ischemic regions.
Figure~\ref{fig:test_case2_potential} reports also the extracellular potassium concentration \([\text{K}^+]_\text{o}\) (panels (e)--(h)) in the range \(5.5\text{--}7.2~\text{mM}\): the tissue is initially close to physiological levels, while during wave propagation a localized increase is observed up to values close to \(7.2~\text{mM}\), particularly in correspondence with regions undergoing repeated activation.
The overall increase is also expected, since each action potential is associated with potassium efflux, leading to a rise of \([\text{K}^+]_\text{o}\) during repeated spiking activity.
This accumulation of extracellular potassium is consistent with a progressive impairment of membrane repolarization and contributes to sustaining the firing dynamics.
Overall, the coupled evolution of the transmembrane potential and potassium concentration suggests a transition from a localized perturbation to a global abnormal regime, in which the ischemic regions become additional drivers of pathological sustained activity through their interaction with the propagating impulses.

%% file: NumericalresultBrain.tex
\subsection{Epileptic seizure simulations and interaction with ischemic region in realistic domain}
\begin{figure}[h]
\centering
\includegraphics[width=0.8\textwidth]{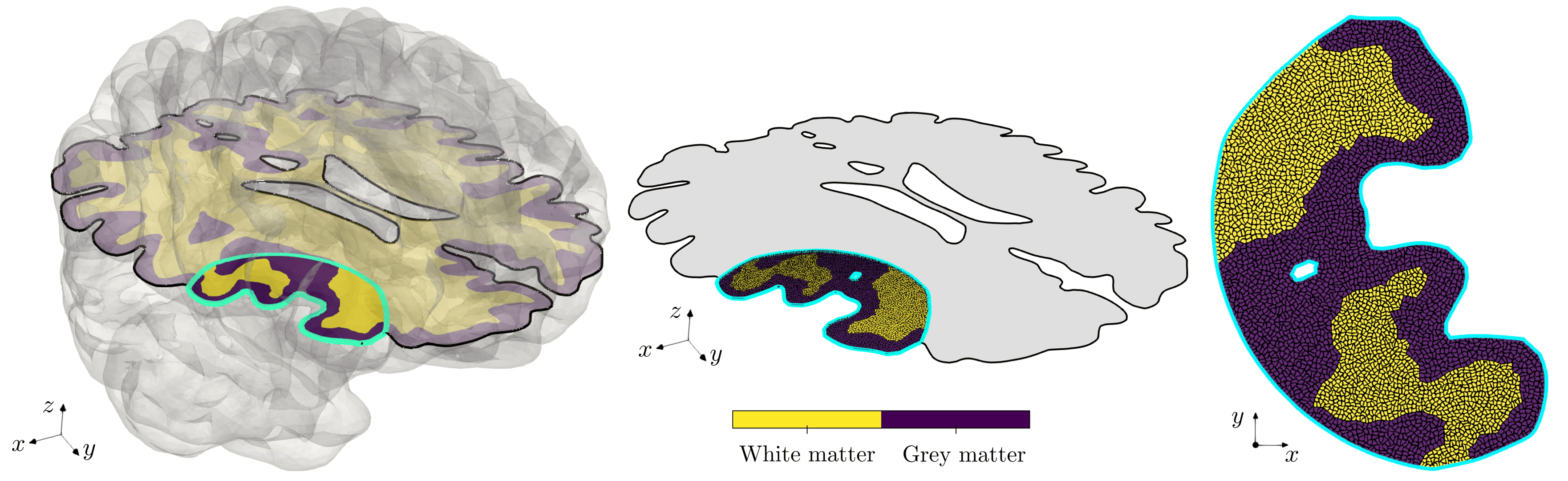}
\caption{Transversal section of the human brain. The computational domain is extracted from the temporo-occipital region of one cerebral hemisphere, including grey matter (purple) and white matter (yellow).}
\label{fig:setupbrain}
\end{figure}

\begin{figure}[h]
\centering
\includegraphics[width=\textwidth]{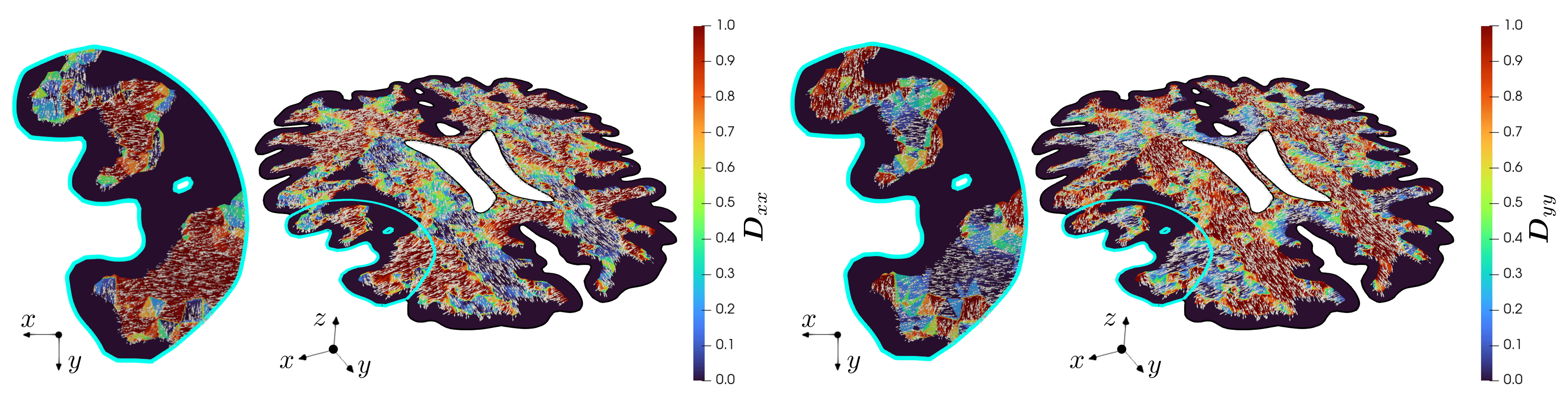}
\caption{Anisotropic diffusion tensors derived from diffusion-weighted imaging (DWI) and used to model axonal conductivity in the white matter.}
\label{fig:setupbrain1}
\end{figure}

In this section, we investigate the electro–metabolic dynamics on a realistic two-dimensional brain sub-domain, extracted from the transversal section covering the temporo-occipital region of one cerebral hemisphere. The domain is discretized using a polytopal mesh specifically generated to conform to the anatomical 
subregion; the resulting mesh consists of 4004 polygonal elements, represented in Figure~\ref{fig:setupbrain}. 
The computational domain includes both grey and white matter, reproducing their anatomical interface. Grey matter is treated as an isotropic conductor, reflecting the homogeneous distribution of neurons and unmyelinated fibers, whereas the white matter is modeled as anisotropic, with enhanced conductivity along the local axonal directions obtained from diffusion–weighted imaging (DWI), where component $\boldsymbol{D}_{xx}$ and $\boldsymbol{D}_{yy}$ are illustrated in Figure~\ref{fig:setupbrain1}. This distinction enables a more realistic representation of the heterogeneous electrophysiological behavior in cortical and subcortical regions.

\begin{figure}[h]
\centering
\includegraphics[width=\textwidth]{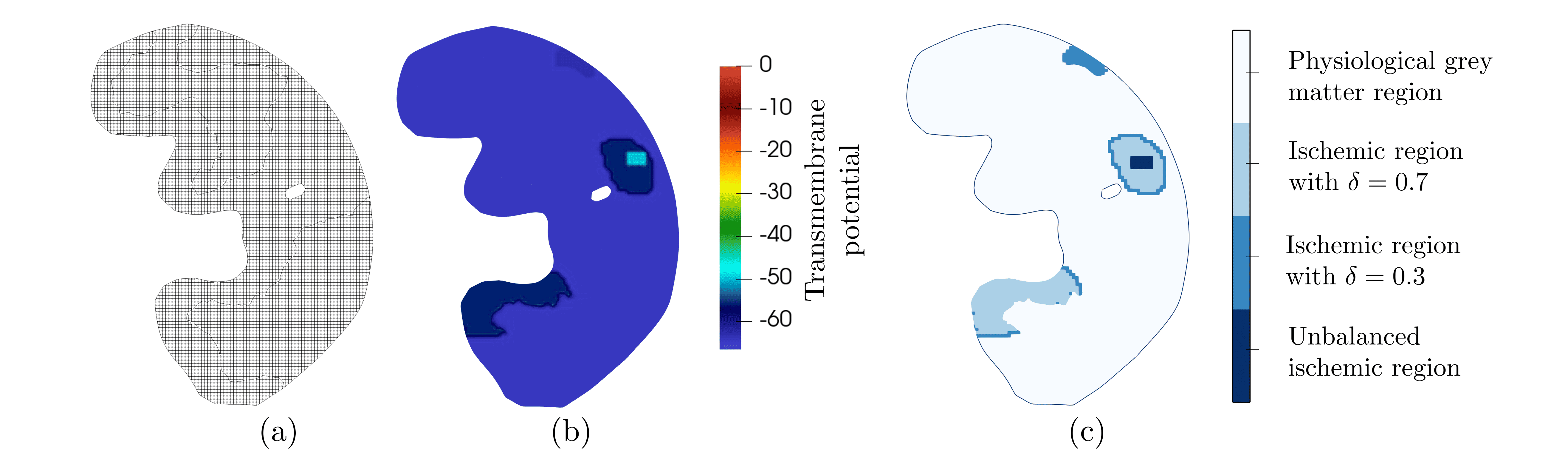}
\caption{Initial setup: (a) Computational mesh of the realistic brain section. (b) Initial condition of the transmembrane potential. (c) Distribution of different physiological and pathological regions, including ischemic subregions of different severities; one of them is externally activated (dark blue region). The computational domain accounts for anisotropic axonal directions influencing the propagation dynamics.}
\label{fig:setupbraininitial}
\end{figure}

\begin{figure}[H]
\centering
\includegraphics[width=\textwidth]{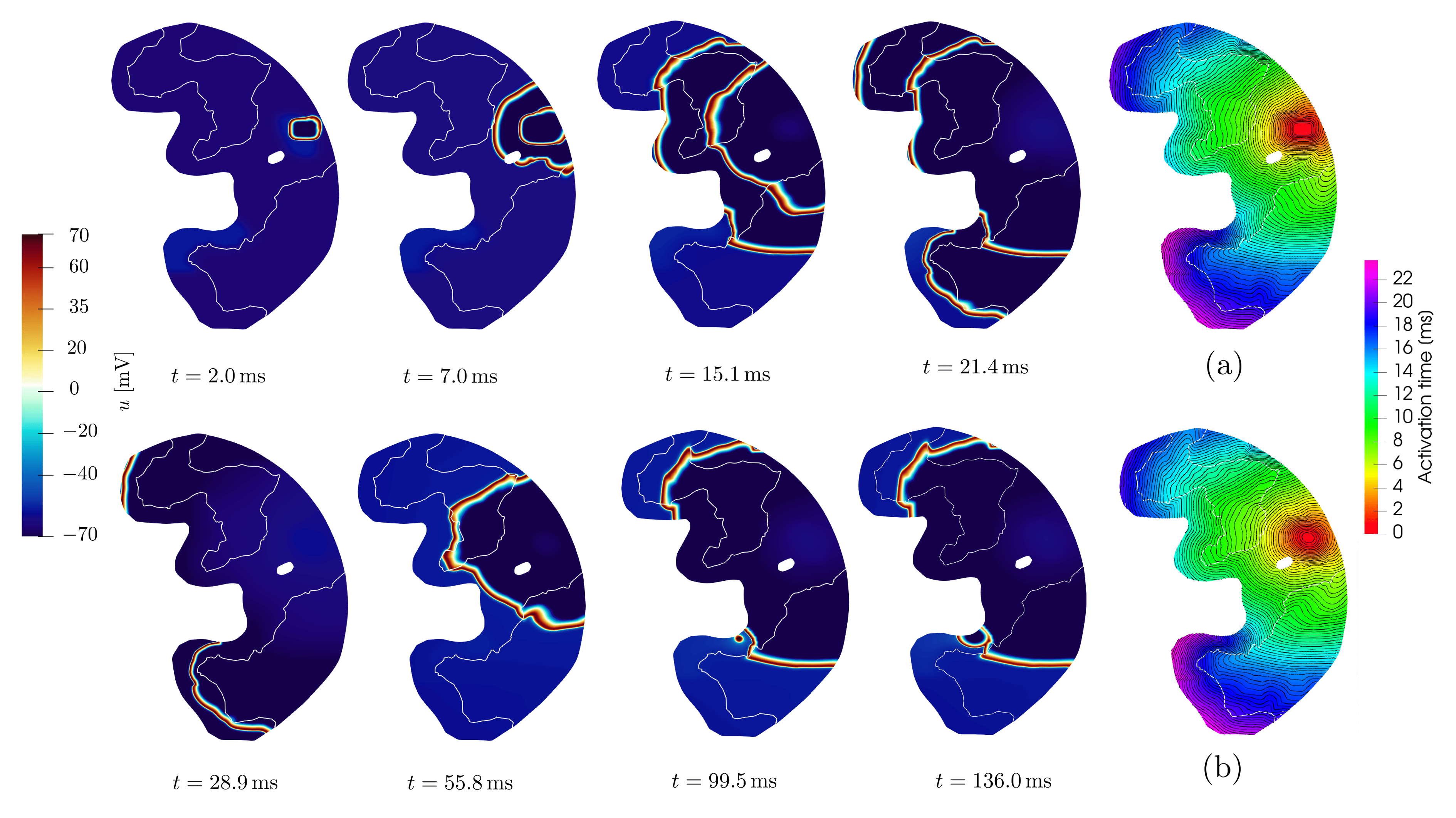}
\caption{Different snapshots of the transmembrane potential evolution in a realistic brain section, showing multiple ischemic regions of varying severity (from subclinical to severe). One of the ischemic zones is externally activated, generating two propagating wavefronts that interacts with surrounding tissue. The computational domain includes both white and grey matter regions, and anisotropic axonal directions that influence the propagation dynamics of the membrane potential over time. (a) Activation times for the first wavefront. (b) Activation times for the second wavefront}
\label{fig:brain_potential}
\end{figure}

\noindent
Following this setup, the computational model is initialized as shown in Figure~\ref{fig:setupbraininitial}, which illustrates the mesh, the initial condition of the transmembrane potential, and the distribution of the physiological and pathological regions.
The domain includes multiple ischemic regions of different severities ($\delta = 0.3$ and $\delta = 0.7$), prescribed a priori to represent subclinical and severe ischemic conditions. These regions are located within the grey matter, consistently with physiological evidence indicating that ischemic damage predominantly affects grey matter tissue. The initial condition is constructed following the configuration shown in Figure \ref{fig:setupbraininitial} and the temporal evolution of the transmembrane potential is presented in Figure~\ref{fig:brain_potential}.

In this test case, we aim to investigate the ischemic tissue behavior when it is stimulated by an external impulse, introduced through an additional forcing term in the coupled model in Equation \eqref{eq:monodomain}. The additional external stimulus is applied for a limited time interval, and it is defined as a constant value $I_\text{Act}=10\ \mu \text{A}$. The external forcing term is applied at early times within a localized ischemic area characterized by a severe blood flow reduction ($\delta$ = 0.7) (see Figure \ref{fig:setupbraininitial}). 
As these wavefronts evolve, they traverse regions of grey and white matter, encountering different conductivity and excitability properties.

\begin{figure}[h]
    \centering
    \begin{subfigure}[b]{.295\textwidth}
    \centering
\includegraphics[width=\textwidth]{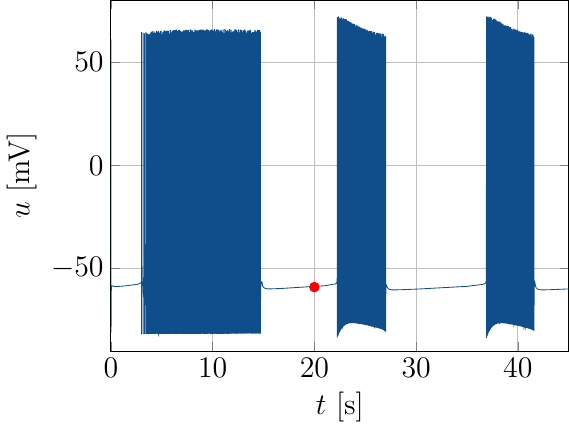}%
\caption{\label{fig:sec2dtest2}}
    \end{subfigure}\hspace{4em}
    \begin{subfigure}[b]{.5\textwidth}
    \centering
\includegraphics[width=\textwidth]{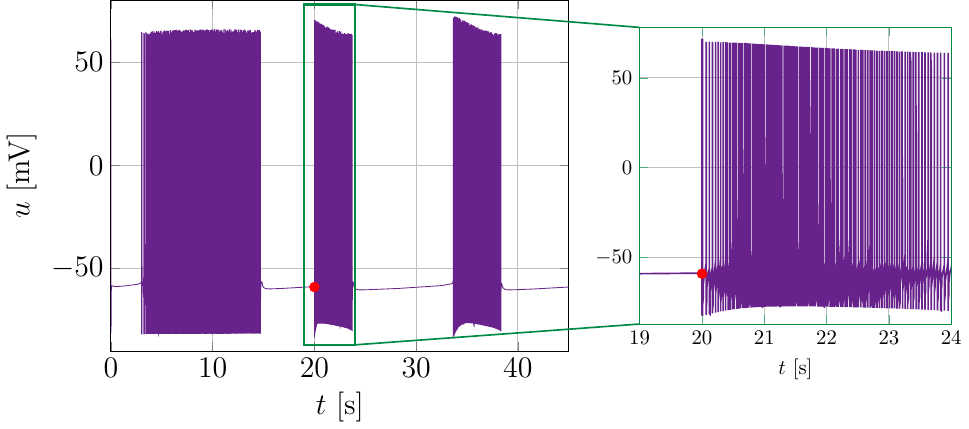}%
\caption{\label{fig:sec2dtest2_zoom}}
\end{subfigure}
\caption{Severe ischemia with external stimulus (the physiological delays are neglected). (a) Pathological evolution with $\varepsilon = 5$ and an external stimulus applied at time $T=20 \;\mathrm{s}$ (red point) for only few time steps. (b) An entire firing episode generated by the external stimulus.}
    \label{fig:testcase2_0d_zoom}
\end{figure}

We consider, in order to analyze the coupled behavior, the parameters $q_0 = 0.3q$, $\eta_{n} = 0.5409$, $\eta_{a} = 0.3938$, and $\eta_{ecs} = 0.0653$ for the severe ischemic subregions. The additional forcing term is both temporal and spatially limited: from the modeling point of view, the forcing is added to the currents of the monodomain model 
\begin{equation}
f(\mathrm{u}, \mathbf{m}, \mathbf{c}) = I_{\text{ion}}(\mathrm{u}, \mathbf{m}, \mathbf{c}) = I_{Na}(u, \mathbf{m}, \mathbf{c}) + I_K(\mathrm{u},\mathbf{m}, \mathbf{c}) + I_{Cl}(\mathrm{u}, \mathbf{m}, \mathbf{c}) + I_\text{Act}\mathbbm{1}_{\{t_1 \le t \le t_2 \} \{ \boldsymbol{x} \in \Omega_\text{1}\}} ,
\label{eq:force_act}
\end{equation}
where $\mathbbm{1}_{\{t_1 \le t \le t_2 \} \{ \boldsymbol{x} \in \Omega_\text{1}\}} $ represents the indicator function on the elements of the mesh corresponding to a small region of the ischemic grey matter, that can be defined as $\Omega_\text{1}$. For the simulation, we define $t_1 = 0 \;\mathrm{ms}$ and $t_2 = 1.2 \;\mathrm{ms}$, so that as soon as the first 2 spikes are triggered, the external force is turned off. In the part of the domain in which we do not impose ischemic condition, we set physiological tissue. The initial conditions for all the variables of the ionic model are taken from the neuronal simulation represented in Figure \ref{fig:testcase2_0d_zoom}, at time $ \hat{t} = 20 \; \mathrm{s}$. 

In this scenario, we can observe that the time-limited additional stimulus generates two action potentials. However, before stabilizing again, the transmembrane potential undergoes a complete bursting behavior with high-frequency dynamics (see Figure \ref{fig:testcase2_0d_zoom}). This time-limited external impulse triggers a pathological behavior when applied in an ischemic region, which proves to be particularly sensitive to external stimuli. In white matter, the directional preferences associated with axonal fiber orientation are captured by the anisotropic conductivity tensor $\boldsymbol{\Sigma}$, which enhances propagation along the dominant fiber direction. In contrast, in grey matter we observe a spatially uniform propagation. Specifically, the conductivity values are listed in Table \ref{table:condtestbrain}, and for the white matter region, we exploit the axonal directions represented in Figure \ref{fig:setupbrain1}. Within the ischemic tissue the conductivity is further reduced, according to Equation \eqref{eq:diffusion_tortuosity}.

\begin{table}[h]
    \centering
    \begin{tabular}{|c|c|c|}
        \hline
        \text{Tissue type} & $\sigma_n$ [Sm$^{-1}$] & $\sigma_l$ [Sm$^{-1}$] \\
        \hline
        Grey matter   &  1.3354 & 1.3354  \\
        Ischemic grey matter & 0.2712 & 0.2712 \\
        \hline
    \end{tabular}
\caption{Conductivity values for grey matter in physiological and pathological conditions. 
    $\sigma_n$ and $\sigma_l$ denote the conductivities in the normal and tangential directions, respectively.}
    \label{table:condtestbrain}
\end{table}

The evolution of the transmembrane potential is reported in Figure \ref{fig:brain_potential}. We can observe that the ischemic region is particularly sensitive with respect to external stimuli; once the ischemic region is triggered and two initial action potentials are generated, it remains active and continues to generate pathological impulses. The first auto-induced wave starts at approximately $t = 45 \; \mathrm{ms}$ without any external forcing term; this pathological evolution is in fact auto-induced by the coupled model and generated from the ischemic region.

In subsequent frames ($t > 90~\mathrm{ms}$), new pathological fronts originate from a different severely ischemic core and invade the adjacent healthy tissue, effectively perturbing the normal propagation dynamics and synchronizing the whole domain into a pathological activation pattern.
In Figure \ref{fig:brain_potential}a and \ref{fig:brain_potential}b are shown the activation times associated with the propagation of electrical self-induced impulses within the two-dimensional domain.
The resulting dynamics highlight the complex interplay between ischemic and non-ischemic zones: the waves originating in the forced pathological region interact with neighboring ischemic patches and propagate into healthy tissue, generating activation patterns shaped by anatomical anisotropy and oxygen-dependent alterations in membrane excitability.

%% file: sample.bib
@article{cangiani_hp-version_2014,
  title     = {hp-version discontinuous {G}alerkin methods on polygonal and polyhedral meshes},
  author    = {Cangiani, Andrea and Georgoulis, Emmanuil H and Houston, Paul},
  journal   = {Mathematical Models and Methods in Applied Sciences},
  volume    = {24},
  number    = {10},
  pages     = {2009--2041},
  year      = {2014},
}

@article{dijkstra2016biophysical,
  title     = {A Biophysical Model for Cytotoxic Cell Swelling},
  author    = {Dijkstra, Kirsten and Hofmeijer, Jeannette and van Gils, Stephan A. and van Putten, Michel J.A.M.},
  journal   = {Journal of Neuroscience},
  volume    = {36},
  year      = {2016},
  publisher = {Society for Neuroscience}
}

@article{leimer2024highorder,
title = {A high-order discontinuous {G}alerkin method for the numerical modeling of epileptic seizures},
journal = {Computers \& Mathematics with Applications},
volume = {205},
pages = {112-131},
year = {2026},
author = {Caterina B. {Leimer Saglio} and Stefano Pagani and Mattia Corti and Paola F. Antonietti},
}

@article{nicholson1998diffusion_tortuosity,
  title={Extracellular space structure revealed by diffusion analysis},
  author={Nicholson, Charles and Sykov{\'a}, Eva},
  journal={Trends in Neurosciences},
  volume={21},
  pages={207--215},
  year={1998},
  publisher={Elsevier}
}

@article{beard2022astrocytes,
  title={Astrocytes as Key Regulators of Brain Energy Metabolism: New Therapeutic Perspectives},
  author={Beard, E. and Lengacher, S. and Dias, S. and Magistretti, P.J. and Finsterwald, C.},
  journal={Frontiers in Physiology},
  volume={12},
  pages={825816},
  year={2022},
  publisher={Frontiers}
}

@article{clarke1999_energy,
  title={Basic neurochemistry: molecular, cellular and medical aspects},
  author={Clarke, Donald D and Sokoloff, L and Siegel, GJ},
  journal={Lippincott-Raven, Philadelphia},
  year={1999}
}

@article{kety1957metabolism_energy,
  title={Metabolism of the nervous system},
  author={Kety, SS},
  journal={Edited by Richter, D. New Y\&k},
  volume={232},
  year={1957}
}

@book{cipolla2009cerebral,
  author    = {Marilyn J. Cipolla},
  title     = {The Cerebral Circulation},
  publisher = {Morgan \& Claypool Life Sciences},
  year      = {2009},
  address   = {San Rafael, CA},
  series    = {Integrated Systems Physiology: From Molecule to Function},
  chapter   = {5}
}

@article{lee2000brain,
  title   = {Brain tissue responses to ischemia},
  author  = {Lee, Jin-Moo and Grabb, Margaret C. and Zipfel, Gregory J. and Choi, Dennis W.},
  journal = {Journal of Clinical Investigation},
  volume  = {106},
  number  = {6},
  pages   = {723--731},
  year    = {2000}
}

@article{chen2022poststroke,
  title   = {Pathogenesis of seizures and epilepsy after stroke},
  author  = {Chen, Jiayu and Ye, Haijiao and Zhang, Jie and Li, Aihong and Ni, Yaohui},
  journal = {Acta Epileptologica},
  volume  = {4},
  year    = {2022}
}

@article{cell_volume_seizures_ullah,
    author = {Ullah, Ghanim AND Wei, Yina AND Dahlem, Markus A AND Wechselberger, Martin AND Schiff, Steven J},
    journal = {PLOS Computational Biology},
    publisher = {Public Library of Science},
    title = {The Role of Cell Volume in the Dynamics of Seizure, Spreading Depression, and Anoxic Depolarization},
    year = {2015}
}

@article{Hellas2021NeuronalSwelling,
  title        = {Neuronal Swelling: A Non-osmotic Consequence of Spreading Depolarization},
  author       = {Hellas, Julia A. and Andrew, R. David},
  journal      = {Frontiers in Cellular Neuroscience},
  year         = {2021},
  volume       = {15},
}

@article{stefanescu2012computational,
  title={Computational models of epilepsy},
  author={Stefanescu, Roxana A and Shivakeshavan, RG and Talathi, Sachin S},
  journal={Seizure},
  volume={21},
  number={10},
  pages={748--759},
  year={2012},
  publisher={Elsevier}
}

@article{saetra2023neural,
  title={Neural activity induces strongly coupled electro-chemo-mechanical interactions and fluid flow in astrocyte networks and extracellular space - {A} computational study},
  author={S{\ae}tra, Marte J and Ellingsrud, Ada J and Rognes, Marie E},
  journal={PLoS Computational Biology},
  volume={19},
  number={7},
  pages={e1010996},
  year={2023},
  publisher={Public Library of Science San Francisco, CA USA}
}

@book{mardal2022mathematical,
  title={Mathematical modeling of the human brain: from magnetic resonance images to finite element simulation},
  author={Mardal, Kent-Andr{\'e} and Rognes, Marie E and Thompson, Travis B and Valnes, Lars Magnus},
  year={2022},
  publisher={Springer Nature}
}

@ARTICLE{Bassi201245,
	author = {Bassi, F. and Botti, L. and Colombo, A. and Di Pietro, D.A. and Tesini, P.},
	title = {{On the flexibility of agglomeration based physical space discontinuous {G}alerkin discretizations}},
	year = {2012},
	journal = {Journal of Computational Physics},
	volume = {231},
	number = {1},
	pages = {45 – 65},
	url = {https://www.scopus.com/inward/record.uri?eid=2-s2.0-80054874713&doi=10.1016%2fj.jcp.2011.08.018&partnerID=40&md5=7d8e9a0256d5b309fcef1ffa374d0212}
}

@article{Potassio_buffer_Contreras2021,
  author    = {Contreras, S. A. and Schleimer, J. H. and Gulledge, A. T. and Schreiber, S.},
  title     = {Activity-mediated accumulation of potassium induces a switch in firing pattern and neuronal excitability type},
  journal   = {PLoS Computational Biology},
  year      = {2021},
  volume    = {17}
}

@article{jaeger2022deriving,
  title={Deriving the bidomain model of cardiac electrophysiology from a cell-based model; properties and comparisons},
  author={J{\ae}ger, Karoline Horgmo and Tveito, Aslak},
  journal={Frontiers in Physiology},
  volume={12},
  pages={811029},
  year={2022},
  publisher={Frontiers Media SA}
}

@article{antonietti2025lymph,
  title={lymph: discontinuous po{LY}topal methods for {M}ulti-{PH}ysics differential problems},
  author={Antonietti, Paola F and Bonetti, Stefano and Botti, Michele and Corti, Mattia and Fumagalli, Ivan and Mazzieri, Ilario},
  journal={ACM Transactions on Mathematical Software},
  volume={51},
  number={1},
  pages={1--22},
  year={2025},
  publisher={ACM New York, NY}
}

@article{cangiani2022hp,
  title={$hp$-version discontinuous {G}alerkin methods on essentially arbitrarily-shaped elements},
  author={Cangiani, Andrea and Dong, Zhaonan and Georgoulis, Emmanuil},
  journal={Mathematics of Computation},
  volume={91},
  number={333},
  pages={1--35},
  year={2022}
}

@article{Antonietti2013A1417,
	author = {Antonietti, Paola F. and Giani, Stefano and Houston, Paul},
	journal = {SIAM Journal on Scientific Computing},
	number = {3},
	pages = {A1417--A1439},
	title = {Hp-Version composite discontinuous {G}alerkin methods for elliptic problems on complicated domains},
	url = {https://www.scopus.com/inward/record.uri?eid=2-s2.0-84884956201&doi=10.1137%2f120877246&partnerID=40&md5=48d57b2f5c5121ad7694653ca5868908},
	volume = {35},
	year = {2013}}

@article{HamBathe2012_FEM_wave,
  author    = {Sanghyun Ham and Klaus-J{\"u}rgen Bathe},
  title     = {A finite element method enriched for wave propagation problems},
  journal   = {Computers \& Structures},
  volume    = {94--95},
  year      = {2012}
}

@article{kager2007seizure,
  title={Seizure-like afterdischarges simulated in a model neuron},
  author={Kager, H and Wadman, WJ and Somjen, GG},
  journal={Journal of Computational Neuroscience},
  volume={22},
  pages={105--128},
  year={2007},
  publisher={Springer}
}

@article{rungta2015cellular_edema,
  title     = {The Cellular Mechanisms of Neuronal Swelling Underlying Cytotoxic Edema},
  author    = {Rungta, Ravi L. and Choi, Hyun B. and Tyson, John R. and Malik, Aqsa and Dissing-Olesen, Lasse and Lin, Paulo J.C. and Cain, Stuart M. and Cullis, Pieter R. and Snutch, Terrance P. and MacVicar, Brian A.},
  journal   = {Cell},
  year      = {2015},
  volume    = {161},
  number    = {3},
  pages     = {610--621}
}

@article{idumah2023spatial_calvetti2,
  author    = {Idumah, Gideon and Somersalo, Erkki and Calvetti, Daniela},
  title     = {A spatially distributed model of brain metabolism highlights the role of diffusion in brain energy metabolism},
  journal   = {Journal of Theoretical Biology},
  year      = {2023},
  volume    = {572}
}

@article{allaman2011astrocyte_neuron,
  title={Astrocyte--neuron metabolic relationships: for better and for worse},
  author={Allaman, Igor and B{\'e}langer, Mireille and Magistretti, Pierre J},
  journal={Trends in neurosciences},
  volume={34},
  pages={76--87},
  year={2011},
  publisher={Elsevier}
}

@article{BONVENTO20211546_neuron_astrocyte,
title = {Astrocyte-neuron metabolic cooperation shapes brain activity},
journal = {Cell Metabolism},
volume = {33},
number = {8},
pages = {1546-1564},
year = {2021},
author = {Gilles Bonvento and Juan P. Bolaños}
}

@article{ShawRudy1997_slow_conduc,
  author       = {Shaw, R. M. and Rudy, Y.},
  title        = {Electrophysiologic effects of acute myocardial ischemia: A mechanistic investigation of action potential conduction and conduction failure},
  journal      = {Circulation Research},
  year         = {1997},
  volume       = {80},
  pages        = {124--138}
}

@article{leimer2025p,
  title={A p-adaptive polytopal discontinuous {G}alerkin method for high-order approximation of brain electrophysiology},
  author={Leimer Saglio, Caterina B. and Pagani, Stefano and Antonietti, Paola F.},
  journal={Computer Methods in Applied Mechanics and Engineering},
  volume={446},
  pages={118249},
  year={2025},
  publisher={Elsevier}
}

@article{hodgkin1990_HH,
  title   = {A Quantitative Description of Membrane Current and its Application to Conduction and Excitation in Nerve},
  author  = {Hodgkin, A. L. and Huxley, A. F.},
  journal = {Bulletin of Mathematical Biology},
  volume  = {52},
  number  = {1-2},
  pages   = {25--71},
  year    = {1990},
  doi     = {10.1007/BF02459568}
}

@article{arnold_unified_2002,
    title       = {Unified analysis of discontinuous {G}alerkin methods for elliptic problems},
    author      = {Arnold, Douglas N and Brezzi, Franco and Cockburn, Bernardo and Marini, L Donatella},
    journal     = {SIAM journal on numerical analysis},
    volume      = {39},
    number      = {5},
    pages       = {1749--1779},
    year        = {2002},
}

@book{cangiani_hp-version_2017,
    title       = {hp-version discontinuous {G}alerkin methods on polygonal and polyhedral meshes},
    author      = {Cangiani, Andrea and Dong, Zhaonan and Georgoulis, Emmanuil H and Houston, Paul},
    year        = {2017},
    publisher   = {Springer},
}

@article{attwell2001energy,
  title={An energy budget for signaling in the grey matter of the brain},
  author={Attwell, David and Laughlin, Simon B.},
  journal={Journal of Cerebral Blood Flow \& Metabolism},
  volume={21},
  number={10},
  pages={1133--1145},
  year={2001},
  publisher={SAGE Publications},
  doi={10.1097/00004647-200110000-00001},
  pmid={11598490}
}

@article{cressman2009_BC1,
  title={The influence of sodium and potassium dynamics on excitability, seizures, and the stability of persistent states: I. Single neuron dynamics},
  author={Cressman, John R. and Ullah, Ghanim and Ziburkus, Jokubas and Schiff, Steven J. and Barreto, Ernest},
  journal={Journal of Computational Neuroscience},
  volume={26},
  number={2},
  pages={159--170},
  year={2009},
  publisher={Springer},
  doi={10.1007/s10827-008-0132-4}
}

@article{CALVETTI2018,
title = {A computational model integrating brain electrophysiology and metabolism highlights the key role of extracellular potassium and oxygen},
journal = {Journal of Theoretical Biology},
volume = {446},
pages = {238-258},
year = {2018},
issn = {0022-5193},
doi = {https://doi.org/10.1016/j.jtbi.2018.02.029},
url = {https://www.sciencedirect.com/science/article/pii/S0022519318300900},
author = {D. Calvetti and G. {Capo Rangel} and L. {Gerardo Giorda} and E. Somersalo}
}

@article{somjen2008computer,
  title={Computer simulations of neuron-glia interactions mediated by ion flux},
  author={Somjen, GG and Kager, H and Wadman, WJ},
  journal={Journal of computational neuroscience},
  volume={25},
  pages={349--365},
  year={2008},
  publisher={Springer}
}

@article{barreto2011ion,
  title={Ion concentration dynamics as a mechanism for neuronal bursting},
  author={Barreto, Ernest and Cressman, John R},
  journal={Journal of biological physics},
  volume={37},
  pages={361--373},
  year={2011},
  publisher={Springer}
}

@article{hodgkin1952propagation,
  title={Propagation of electrical signals along giant nerve fibres},
  author={Hodgkin, Alan Lloyd and Huxley, Andrew Fielding},
  journal={Proceedings of the Royal Society of London. Series B-Biological Sciences},
  volume={140},
  number={899},
  pages={177--183},
  year={1952},
  publisher={The Royal Society London}
}

@article{schreiner2022simulating,
  title={Simulating epileptic seizures using the bidomain model},
  author={Schreiner, Jakob and Mardal, Kent-Andre},
  journal={Scientific Reports},
  volume={12},
  number={1},
  pages={10065},
  year={2022},
  publisher={Nature Publishing Group UK London}
}

@article{CAPORANGEL2020110093,
title = {Metabolism plays a central role in the cortical spreading depression: Evidence from a mathematical model},
journal = {Journal of Theoretical Biology},
volume = {486},
pages = {110093},
year = {2020},
issn = {0022-5193},
author = {G Capo-Rangel and L Gerardo-Giorda and E Somersalo and D Calvetti}
}

@article{koppl20203d,
  title={A 3D-1D coupled blood flow and oxygen transport model to generate microvascular networks},
  author={K{\"o}ppl, Tobias and Vidotto, Ettore and Wohlmuth, Barbara},
  journal={International journal for numerical methods in biomedical engineering},
  volume={36},
  number={10},
  pages={e3386},
  year={2020},
  publisher={Wiley Online Library}
}

@article{goriely2015mechanics,
  title={Mechanics of the brain: perspectives, challenges, and opportunities},
  author={Goriely, Alain and Geers, Marc GD and Holzapfel, Gerhard A and Jayamohan, Jayaratnam and J{\'e}rusalem, Antoine and Sivaloganathan, Sivabal and Squier, Waney and van Dommelen, Johannes AW and Waters, Sarah and Kuhl, Ellen},
  journal={Biomechanics and modeling in mechanobiology},
  volume={14},
  number={5},
  pages={931--965},
  year={2015},
  publisher={Springer}
}

@ARTICLE{muller,
AUTHOR={Müller, Lucas O.  and Watanabe, Sansuke M.  and Toro, Eleuterio F.  and Feijóo, Raúl A.  and Blanco, Pablo J. },   
TITLE={An anatomically detailed arterial-venous network model. Cerebral and coronary circulation},    
JOURNAL={Frontiers in Physiology},      
VOLUME={Volume 14 - 2023},
YEAR={2023},
ISSN={1664-042X}
}

@article{Jiruska2011_high_freq,
  author    = {Premysl Jiruska and Anatol Bragin},
  title     = {High-frequency activity in experimental and clinical epileptic foci},
  journal   = {Epilepsia},
  year      = {2011},
  volume    = {52 Suppl 8},
  pages     = {28--36}
}
